\documentclass[ssy]{imsart}
\usepackage{epsfig,epstopdf,mathtools,amssymb} 
\usepackage[colorlinks = true, citecolor = blue, urlcolor = blue]{hyperref}

\startlocaldefs

\newtheorem{remark}{Remark}[section]

\newcommand{\RR}{{\mathbb R}}
\newcommand{\ZZ}{{\mathbb Z}}

\newcommand{\eeq}{\end{equation}}

\newcommand{\ra}{\rightarrow}
\newcommand{\Ra}{\Rightarrow}

\newcommand{\beql}[1]{\begin{equation}\label{#1}}
\newcommand{\eqn}[1]{(\ref{#1})}
\newcommand{\beq}{\begin{displaymath}}
\newcommand{\eeqno}{\end{displaymath}}

\newcommand{\lm}{\lambda}

\newcommand{\qandq}{\quad\mbox{and}\quad}
\newcommand{\qifq}{\quad\mbox{if}\quad}

\newcommand{\qforq}{\quad\mbox{for}\quad}

\newcommand{\qasq}{\quad\mbox{as ~}}

\newcommand{\qinq}{\quad\mbox{in ~}}
\newcommand{\qforallq}{\quad\mbox{for all}\quad}

\newcommand{\bes}{\begin{equation*}}
\newcommand{\ees}{\end{equation*}}
\newcommand{\bequ}{\begin{equation}}
\newcommand{\ei}{\end{itemize}}
\newcommand{\bsplit}{\begin{split}}
\newcommand{\esplit}{\end{split}}
\newcommand{\bea}{\begin{eqnarray}}
\newcommand{\eea}{\end{eqnarray}}
\newcommand{\beas}{\begin{eqnarray*}}
\newcommand{\eeas}{\end{eqnarray*}}
\newcommand{\btab}{\begin{tabular}}
\newcommand{\etab}{\end{tabular}}

\newcommand{\barq}{\bar{Q}}
\newcommand{\barz}{\bar{Z}}

\newcommand{\barx}{\bar{X}}

\newcommand{\1}{\textbf{1}}

\def\lm{\lambda}
\def\tinf{\rightarrow\infty}

\def\AA{\mathbb{A}}
\def\BB{\mathbb{B}}

\def\mA{\mathcal{A}}

\endlocaldefs

\begin{document}

\begin{frontmatter}

\title{Achieving Rapid Recovery in an Overload Control for Large-Scale Service Systems}
\runtitle{Rapid Recovery in Overload Control}


\author{\fnms{Ohad} \snm{Perry}\ead[label=e1]{ohad.perry@northwestern.edu}}
\address{Department of Industrial Engineering and Management Sciences, Northwestern University, \\
 Evanston, IL 60208 \printead{e1}}
\and
\author{\fnms{Ward} \snm{Whitt}\ead[label=e2]{ww2040@columbia.edu}}
\address{Department of Industrial Engineering and Operations Research, Columbia University \\
 New York, NY, 10027 \printead{e2}}

\runauthor{Perry and Whitt}

\begin{abstract}
We consider an automatic overload control
for two large service systems modeled as multi-server queues, such as call centers.
We assume that the two systems
are designed to operate
independently, but want to help each other respond to unexpected overloads.
The proposed overload control automatically activates
sharing (sending some customers from one system to the other)
once a ratio of the queue lengths in the two systems crosses an activation
threshold (with ratio and activation threshold parameters for
each direction).  To prevent harmful sharing, sharing is allowed
in only one direction at any time.  In this paper, we are primarily concerned
with ensuring that the system recovers rapidly after the
overload is over, either (i) because the two systems return to
normal loading or (ii) because the direction of the overload suddenly shifts in the
opposite direction.
To achieve rapid recovery, we introduce lower thresholds for the queue ratios, below which one-way sharing is released.
As a basis for studying the complex dynamics, we develop a new six-dimensional fluid approximation for a system with time-varying arrival rates,
extending a previous fluid approximation involving a stochastic averaging principle.
We conduct simulations to confirm that the new algorithm is effective for predicting the system performance
and choosing effective control parameters.
The simulation and the algorithm both show that the system can experience an inefficient nearly-periodic behavior, corresponding to an oscillating equilibrium
(congestion collapse), if the sharing is strongly inefficient and the control parameters are set inappropriately.
\end{abstract}



\end{frontmatter}


\section{Introduction}\label{secIntro}

\subsection{An Automatic Overload Control}\label{secControl}

In this paper we study an automatic control to temporarily
activate ``emergency'' measures in an uncertain dynamic
environment to mitigate damage from an unexpected disruption,
and then automatically return to normal operation when the
disruption is over. There are two important questions:  First,
how and when should the control be activated? And, second, how
and when should the control be released?

These issues arise in many contexts and have long been studied
within the discipline of control theory \cite{K02,S98}. A
familiar automatic control is a thermostat, which automatically
turns on and off a heater and/or an air conditioner within a
building.  Since building temperature tends to change slowly
relative to human temperature tolerance, conventional
thermostats operate well with little concern, but special
thermostats are needed for complex environments, such as in
biochemical processes \cite{BFI05}.

Another example of an automated control occurs in a large stock
market exchange, such as the New York Stock Exchange (NYSE). To
respond to the experience of dramatic fluctuations in prices,
in 1988 the NYSE instituted trading curbs called {\em circuit
breakers} or collars, which stop trading for a specified period
in the event of exceptionally large price changes.  With the
increase of high-speed computer trading, these controls have
become even more important and interesting since then
\cite{GK04}.

The specific setting considered here involves two large-scale
telephone call centers (or service pools within the same call
center) that are designed to operate independently, but have
the capability (due to both network technology and agent
training) to respond to calls from the other system, even
though there might be some loss in service effectiveness and
efficiency in doing so. These call centers are designed and
managed to separately respond to uncertain fluctuating demand
and, with good practices, usually can do so effectively; see
\cite{AAM07} for background.  However, these call centers may
occasionally face exceptional unexpected overloads, due to sudden surges in arrivals,
extensive agent absenteeism 
or system malfunction (e.g., due to computer failures).  It
thus might be mutually beneficial for the two systems to agree
to help each other during such overload incidents.  We propose
an automatic control for doing so. We are motivated by this
call-center application, but the insights and methods should be
useful in other service systems.  Since we model the call centers as multi-server queues,
the insights and methods may also be useful for other queueing settings.

In telecommunication systems and the Internet, the standard
overload controls reduce the demand through some form of
admission control (rejecting some arrivals) or otherwise
restricting demand; see \cite{BW98,FF99,LPD02,SKT90,WJLH06} and
references therein.  These controls, that reject or reduce
arrivals, are especially important when the increasing load can
cause the useful throughput (the ``goodput'') not only to reach
its largest possible value, but also to actually decrease. Such
anomalous behavior can occur because some of the customers ``go
bad.''  The classic telephone example is failure during the
call setup process.  The customer might start entering digits
before receiving dial tone or abandon before the call is sent
to the destination.  As a consequence, the vast majority of
system resources may be working on requests that are no longer
active, causing the throughput to actually decrease. In
response, various effective controls have been developed
\cite{DH86,Ko91}.

In contrast, here we assume that {\em no} arrivals will be
directly turned away, although on their own initiative
customers may elect to abandon from queue because they become
impatient.  Instead, we develop a control that automatically
sends some of the arrivals to receive service from the other
service pool when appropriate conditions are met. It is natural
to prefer diverting instead of rejecting arrivals whenever some
response is judged to be better than none at all, even if
delayed. Indeed, diverting instead of rejecting arrivals is the
accepted policy with ambulance diversion in response to
overload in hospital emergency rooms, e.g., see
\cite{BBHA12,DG11,YGGG10} and references therein. The results
here may be useful in that context as well, but then it is
necessary to consider the extra delay for ambulances to reach
alternative hospitals, which has no counterpart in networked
call centers. (We assume that the calls can be transferred
instantaneously.)

\subsection{Congestion Collapse}\label{secCollapse}

An important feature of this kind of sharing, which is captured by our our model, is that the sharing may be
inefficient.  A simple symmetric example that we will consider
in \S \ref{secOSC} has identical service rates for agents
serving their own customers, but identical slower service rates
when serving the other customers.
With such inefficiency,
 the whole system will
necessarily operate inefficiently, with lower throughput of
both classes, if both pools are busy serving the other
customers instead of their own.
Nevertheless, we find that judicious sharing with our proposed overload control
can be effective even with some degree of inefficiency, but care is needed in setting the control parameters.
A major concern with such
inefficient sharing is that the system may possibly experience
{\em congestion collapse}, i.e., the system may reach an
equilibrium with inefficient operation \cite{SW11}.

It is known that control schemes can cause congestion collapse; see, e.g., \cite{EF91}.
 Within telecommunications there is a long history of congestion collapse and its prevention
in the circuit-switched telephone network.  More than 60 years
ago, it was discovered that the capacity and performance of the
network could greatly be expanded by allowing alternative
routing paths \cite{W56}. If a circuit is not available on the
most direct path, then the switch can search for free circuits
on alternative paths. The difficulty is that these alternative
paths may use more links and thus more circuits.  Thus, in
overload situations (the classic example being Mother's Day),
the network can reach a stable inefficient operating regime,
with the system congested, but far less than maximal
throughput. This congestion collapse in the telephone network
was first studied by simulation \cite{W64}. The classical
remedy in such loss networks is trunk reservation control,
where the last few circuits on a link are reserved for direct
traffic; see \cite{FR94}, \S\S 4.3-4.5 of \cite{K91} and
references therein.

Overload controls have also been considered for more general
multi-class loss networks. In the multi-class setting, it may
be desirable to provide different grades of service to
different classes, including protection against overloads
caused by overloads of other classes. Partial sharing controls
achieving these more general goals can be achieved exploiting
upper limit bounds and guaranteed minimum bounds \cite{CLW95}.
Moreover, in \cite{CLW95} algorithms are developed to compute
the performance associated with such complex controls, which
greatly facilitates choosing appropriate control parameters.
For the (different) problem we consider, we also develop a
performance algorithm that can be used to set the control
parameters.

Even though a call center can be regarded as a
telecommunications network, our problem is quite different from
the classical loss network setting discussed above.  By
definition, the loss network has {\em no} queues, so that all
arrivals that cannot immediately enter service are turned away.
In sharp contrast, our system turns {\em no} arrivals away.  As
a consequence, our system is more ``sluggish;'' it responds
more slowly to changes in conditions, and presents new
challenges.

For the model considered here,
we show in \S \ref{secOSC}
that the two call centers can indeed experience behavior that is best described as congestion collapse
if the sharing is strongly inefficient and an inappropriate
control is used.
An unstable oscillating equilibrium is predicted by our numerical algorithm for the approximating fluid model
and confirmed by simulation; see Figures \ref{figZ12oscAb} and \ref{figQ2oscAb} for the simulation and Figures \ref{figZ12oscAbFluid}
and \ref{figQ2oscAbFluid} for the algorithm.

However, this oscillatory phenomenon is far from obvious because the stochastic model
after the overload is over is an ergodic time-homogeneous CTMC with a steady-state limiting distribution.
The situation that we consider in this paper is similar to the nearly periodic behavior of the $G/D/s+GI$ queue exposed in
\cite{LW11}.  In that setting, the actual stochastic system has a well-defined limiting steady-state distribution
and yet the system exhibits nearly periodic behavior over long time periods.  When the scale is large, it turns out that the nearly periodic transient behavior
observed in simulations is
well predicted by a limiting fluid model.  Unlike the stochastic model, the fluid model does not have a unique limiting steady-state.
The reason for this discrepancy is that the two iterated limits (as time gets large and as the scale, determined by the arrival rate, gets large) done in different order are not equal.

In this paper we show the existence of the nearly periodic behavior (with inefficient sharing and inappropriately chosen controls), tantamount to congestion collapse,
with our fluid algorithm and simulation.  We provide additional mathematical support
in \cite{PW14} by proving that unstable oscillating equilibria
can exist for a class of these fluid models.

However, this highly undesirably behavior can be avoided with reasonably chosen controls.
In this paper we develop a model and an algorithm for analyzing that model that can be used
to achieve the benefits of sharing while avoiding such bad behavior.

\subsection{Fixed-Queue-Ratio Controls}

Our overload control is a modification of
the {\em Fixed-Queue-Ratio} (FQR) and more general
Queue-and-Idleness-Ratio (QIR)
 controls proposed for routing and scheduling in a multi-class multi-pool call center under normal operating conditions
in \cite{GW09a,GW09b,GW10}.  For the two-class two-pool $X$
model considered here, the FQR rule sends customers to the
other service pool if the ratio of the queue lengths exceeds a
specified ratio. 
However, the theorems establishing that the FQR control is
effective in \cite{GW09a,GW09b,GW10} have conditions that {\em do not} hold for our networks here,
which has a cyclic routing graph and service rates that depend on the customer class
and service pool.  
Indeed, Example 2 of \cite{PW09} shows
that the $X$ model can experience severe congestion collapse
under normal loading if $FQR$ is used. (The congestion collapse shown in \cite{PW09} is different than the one mentioned above, which is due to the undesired
oscillatory behavior.)

Nevertheless, in \cite{PW09} we showed that the FQR control can
usefully be applied as an overload control for the $X$ model
with inefficient sharing if we introduce additional activation
thresholds. The {\em FQR control with thresholds} (FQR-T) sends
customers to the other service pool if the queue ratio exceeds
the activation threshold. For the $X$ model, the FQR-T control
has four parameters: a target ratio and an activation threshold
for each direction of sharing. The target ratios are chosen to
minimize the long-run average cost during the overload incident
in an approximating stationary deterministic fluid model with a
convex cost function applied to the two queues. To prevent
harmful sharing, we also imposed the condition of {\em one-way
sharing}; i.e., sharing is allowed in only one direction at any
one time.

To better understand the transient behavior of the FQR-T
control, in \cite{PW11a} we developed a deterministic fluid
model to analyze the performance.  That model is challenging
and interesting because it is an {\em ordinary differential
equation} (ODE) involving a stochastic {\em averaging
principle} (AP).  In \cite{PW11b,PW13a,PW13b} we established
supporting mathematical results about the FQR-T control, including a functional weak law of large numbers (FWLLN) and
functional central limit theorem (FCLT) refinement.  The
previous analysis
showed that the FQR-T control can rapidly
respond to and mitigate an unexpected overload, while
preventing sharing under normal conditions.

\subsection{New Contribution:  Rapid Recovery After the Overload Is Over}

In this paper we show that FQR-T needs to be modified in order
to ensure that the system recovers rapidly after an overload is
over, either (i) because the two systems return to normal
loading or (ii) because the direction of the overload suddenly shifts in the opposite
direction. To achieve rapid recovery, we propose additional
release thresholds for the shared-customers processes, below
which one-way sharing is released.  (We had previously recognized that such a modification of FQR-T
was needed, e.g., see paragraph 3 in \S 2.2 of \cite{PW11b} and Remark B.1 in Appendix B of \cite{PW13a}, but we now show for the first time
that the modified control can be analyzed and can be effective.)

As a basis for studying such
more complex dynamics, we extend our previous fluid model
approximation in two ways:  (i)  the new fluid model is $6$-dimensional instead of $3$-dimensional
and (ii) the model is allowed to have time-varying arrival rates and staffing functions.  We also extend
our previous algorithm to numerically compute the fluid
solution to this more complex model.  We implement the new algorithm and conduct simulations to show that the fluid model and the associated algorithm
are effective in predicting system performance.  Finally, we show that the new {\em FQR control with
activation-and-release thresholds} (FQR-ART) can be effective with appropriate control parameters.
We provide guidelines on choosing appropriate control parameters in the paper.  The new model and algorithm can be used to
confirm that a good choice has been made.

As before, for large scale in the model, there is important {\em state space collapse} (SSC) during overload periods with active sharing of
customers.  As a consequence, even though the basic stochastic process is $6$-dimensional, the approximating fluid model
is essentially
$3$-dimensional when there is active sharing instead of $6$-dimensional. One of the
complications in the new setting is to identify if and when SSC
begins and in what direction (which queue is receiving help),
and when it ends.
The stochastic AP determines when SSC occurs, and how the fluid
evolves during periods of SSC.
When we introduce release thresholds, we also discover that to
achieve good robust performance, we also need to increase the
activation thresholds.

In summary, our contribution is fourfold: (i) We continue our
study of the X model and demonstrate how and when it is
beneficial (or harmful) to exploit system flexibility
in response to an overload. (ii) We improve the previous
FQR-T control designed to automatically exploit that system
flexibility when it is beneficial to do so by ensuring rapid
recovery when the overload has ended. (iii) We develop a novel
fluid model to approximate the intractable stochastic system in
the time-varying environment and help determine appropriate
control parameters.
(iv) Finally, we design an efficient algorithm to solve that
fluid model.  (This paper addresses the control from an engineering perspective, as in \cite{PW11a};
we do not focus on underlying mathematics (prove theorems) as in \cite{PW11b,PW13a,PW13b}.)

Simulation also plays an important role in our study.  First,
we use simulation to show that refinements to the FQR-T control
are needed to ensure rapid recovery after the overload is over.
Second, we use simulation to demonstrate that the fluid model
provides a good performance approximation. Finally, we use
simulation to verify that we can indeed gain important insights
into complex system behavior from the fluid model, even for
systems that are not overloaded, as in our examples after the
overload has ended.


\subsection{Other Related Literature}

\paragraph{Time-Varying Models.}  A significant contribution here is extending the analysis of the transient behavior of a stationary fluid model to the analysis
of (the necessarily transient behavior of) a time-varying model.
When the predictable variability captured by time-varying model
parameters dominates the unpredictable stochastic variability,
deterministic fluid models are especially appropriate.
Operationally, the deterministic fluid models tend to capture
the essential performance. Mathematically, the deterministic
fluid models are much easier to analyze than their stochastic
extensions, such as diffusion approximations. The vast majority
of the queueing literature concerns stationary models, but
there have been important exceptions, e.g., \cite{K72,N82}. For
related recent work, see
\cite{HJM09,HMW09,LW11,LW12a,LW12b,LW12c}, and references
therein.

\paragraph{Overloaded Systems and Fluid Models.}

For other work that considers overloaded systems and fluid
models, see \cite{CAB11,GurPer12,HZ05,ST11}. The authors in
\cite{CAB11} suggest using the max-weight policy which, much
like the FQR-ART control here, is easy to implement because it
uses only information on the current state of the system; it
stabilizes the system during normal loads and keeps the queues
at target ratios when the system cannot be stabilized due to
high arrival rates. In \cite{GurPer12} overflow networks in
heavy-traffic are studied in settings of co-sourcing, i.e.,
when firms that operate their own in-house call center overflow
a {\em nonnegligible proportion} of the arrivals to a call
center that is operated by an outsourcer.
Fluid and diffusion limits are obtained via a stochastic
averaging principle.
The authors in \cite{ST11} apply their previously introduced
{\em shadow routing} control to overloaded parallel systems
with unknown arrival rates, and show that it maximizes the
reward rate, assuming a class-dependent reward of each customer
served.

\paragraph{Healthcare Systems.}

System overloads are especially prevalent in healthcare
systems, often even being the ``natural state.'' Some
facilities such as {\em intensive care units} (ICU's) and
equipment such as magnetic resonance imaging (MRI) machines are
so expensive that they are designed to be operated
continuously, and exhibit long lines of waiting patients.
Extreme overloads can occur with mass casualty events.

Hospitals 
have complex queueing dynamics, with multiple internal flows
among its units in addition to exogenous arrival streams. Thus,
overloads in some units of a hospital can ``propagate'' to
other units, creating a system-wide overload. For example, when
{\em inpatient wards} (IW) are overloaded, patients from the
{\em emergency department} (ED) who need to be hospitalized
cannot be transferred to the IW due to the unavailability of
beds, creating the phenomenon of {\em blocked beds} in the ED,
i.e., beds that are occupied by patients who finished their
treatment in the ED; see, e.g., \cite{BBHA12} for a current
review.
See \cite{ArmonyPatient} for a data-based study of queueing aspects in hospital settings as well as an extensive literature review. 

In \cite{ChanGalit}, a fluid approximation of an ICU
experiencing periods of overload periods is studied, in which
the service rate of current ICU patients increases (is
``sped-up'') if the number of patients that are waiting to be
admitted to the ICU exceeds a certain threshold. In turn, the
sped-up patients have an increased probability of readmission
to the ICU, so that alleviating overloads by employing speedup
increases future overloads.
The fluid model in \cite{ChanGalit} exploits an averaging
principle in the spirit of \cite{PW11a}.



\subsection{Organization of the Rest of the Paper}

In \S \ref{secModel} we define the stochastic X model
and the FQR-T and FQR-ART controls. Building on simple fluid
considerations, In \S \ref{secRelaxOneWay} and \S \ref{secOSC}
we demonstrate the need to modify FQR-T in order to rapidly
recover after the overload is over. In \S \ref{secRelaxOneWay}
we show why release thresholds are needed.  In \S \ref{secOSC}
we show that, unless precaution is taken, the release
thresholds can cause congestion collapse when the system
recovers from an overload. To avoid that bad behavior, the
activation thresholds need to be increased beyond the FQR-T
values. In \S \ref{secFluidMain} we develop the fluid
approximation and in
\S \ref{secPis} we develop an efficient algorithm to
numerically solve it. In \S \ref{secNum} we provide numerical
examples, demonstrating the effectiveness of both the FQR-ART
control and the fluid model by comparing the results of the
numerical algorithm for the ODE to the results of simulation
experiments.
Finally, in \S \ref{secConclude} we draw conclusions and
suggest directions for further research.

\section{The Time-Varying X Model} \label{secModel}

\begin{figure}[h!]
\begin{center}
\includegraphics[width=5cm]{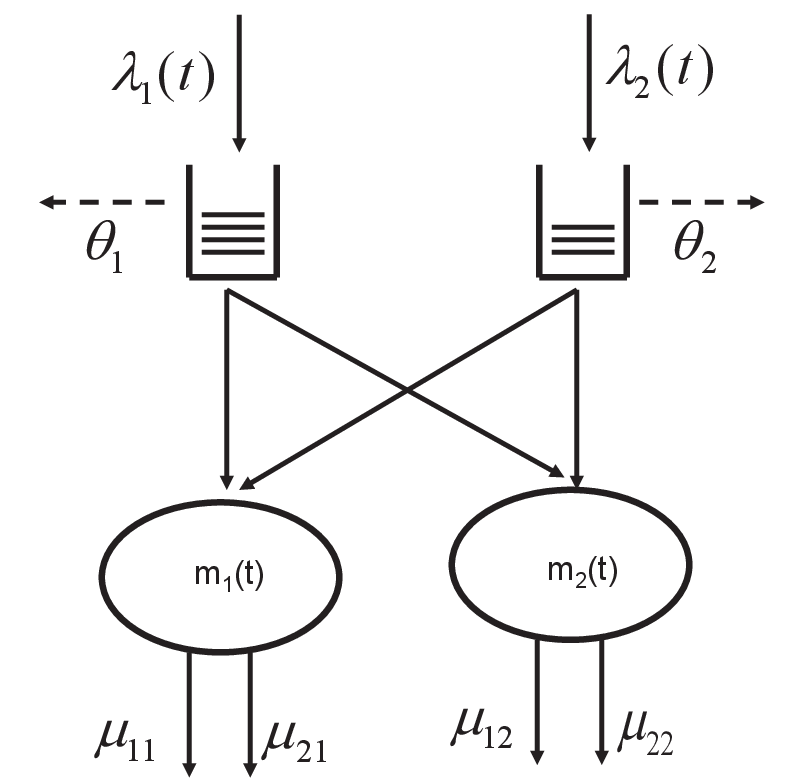}
\caption{The $X$ model}
\end{center}
\label{figX}
\end{figure}

As depicted in Figure \ref{figX}, the $X$ model has two
customer classes and two agent pools, each with many
homogeneous agents working in parallel.
We assume that each customer class has a service pool primarily
dedicated to it, but all agents are cross-trained so that they
can handle calls from the other class, even though they may do
so inefficiently, i.e., customers may be served at a slower
rate when served in the other class pool. We assume that the
service times are independent exponential random variables,
with $1/\mu_{i,j}$ being the expected time for a class $i$ customer
to be served in service pool $j$. Each class has a buffer with
unlimited capacity where customers who are not routed
immediately into service upon arrival wait to be served. Within
each class, customers enter service according to the
first-come-first-served discipline. Customers have limited
patience, so that they may abandon from the queue. The
successive patience times of class $i$ customers are i.i.d.
exponential variables with mean $1/\theta_i$.

We assume that customers arrive according to independent {\em
nonhomogeneous} Poisson processes, one for each class, with
time-varying deterministic rate functions. The staffing levels
are assumed to be time dependent as well, usually chosen to
respond to anticipated changes in the arrival rates; see
\cite{LW12b} and references therein.
As discussed in \S1 of \cite{LW12a}, it is necessary to specify
how the system responds when the staffing level of a service
pool is scheduled to decrease.  As in \cite{LW12a}, we allow
server switching (an agent can take over service from an agent
scheduled to leave).  Since service times are exponential, it
thus suffices to let idle agents leave when staffing decreases,
and the first agent to become idle leave when all agents are
busy when staffing is scheduled to decrease.


Even though we do not prove any limit theorems as we did in
\cite{PW13a,PW13b}, and instead develop direct fluid models to
approximate the stochastic system, we will use asymptotic
considerations in our analysis. We therefore consider a
sequence of X systems, as just described, indexed by a
superscript $n$. As is standard for many-server heavy-traffic
limits \cite{PW13a}, the service rates and abandonment rates are
independent of $n$, but the arrival rates and staffing levels
increase. Specifically, for each $n \ge 1$, let $\lm^n_i(t)$ be
the arrival rate to pool $i$ and let $m^n_j(t)$ be the number
of agents in pool $j$ at time $t$. For the fluid approximation,
we assume that \bequ \label{ParamLim} \lm^n_i(t)/n \ra \lm_i(t)
\qandq m^n_j(t)/n \ra m_j(t) \qasq n \tinf, \eeq uniformly in
$t$ over each bounded time interval.

As in \cite{LW12a}, we assume that the limit functions $\lm_i$
and $m_j$ in \eqn{ParamLim} are piecewise-smooth, by which we
mean that they have only finitely many discontinuities in any
finite interval, have limits from the left and right at each
discontinuity point and are differentiable at all continuity
points.  That assumption is not restrictive for applications
and supports analysis of the approximating fluid model by
differential equations. For call-center applications, it
usually suffices to consider piecewise-constant functions, but
we allow greater generality because our methods can be applied
in other settings.

Let $Q^n_i(t)$ be the number of customers waiting in the
class-$i$ buffer and $Z^n_{i,j}(t)$ be the number of class-$i$
customers in service pool $j$ at time $t$ in system $n$. Let
the associated six-dimensional vector process be \bequ
\label{Xn} X^n \equiv X^n(t) \equiv (Q^n_i(t), Z^n_{i,j}(t) :
i,j = 1,2), \quad t \ge 0. \eeq We consider controls that are
functions of $X^n(t)$ at each $t$, making $X^n$ a
nonhomogeneous CTMC.

To define asymptotic regimes, let $\rho^n_i(t) :=
\lm^n_i(t)/(\mu_{i,i} m^n_i(t))$ be the instantaneous
traffic-intensity function of class $i$ (and pool $i$) alone in
system $n$ at time $t$. By \eqn{ParamLim}, \bequ \label{beta}
\rho^n_i(t)-1 \ra \beta_i(t) \qasq n \ra \infty, \eeq uniformly
in $t$ over each bounded time interval. We say that class $i$
(and pool $i$) is {\em underloaded } at time $t$ if $\beta_i(t)
< 0$, {\em overloaded} at time $t$ if $\beta(t) > 0$ and {\em
normally loaded} at time $t$ if $\beta_i(t) = 0$.

The generality we have introduced allows for many possible
scenarios, but here we restrict attention to an unexpected
overload incident followed by a subsequent instantaneous switch
in state, either (i) a return to normal loading or (ii) a
switch in the direction of overloading. Thus, now there are
three intervals:  first normally loaded, then overloaded and
then a final new regime, which is either normal loading for
both classes or an overload in the opposite direction. During
each of these three intervals, the arrival rates and staffing
functions are allowed to change.

As before, we consider the system starting at the unanticipated
time when the first overload incident begins. However, now the
arrival rates and staffing functions no longer need to be
constant within each interval. By assumption, they have
discontinuities at the beginning of the first overload incident
and at the subsequent time when the overload is over. For the
generality that we do consider, we exploit the fact that we
know how to staff to stabilize the system in face of
time-varying arrival rates under normal loading; see
\cite{LW12a,LW12b} and references therein.

\subsection{The Initial FQR-T Control} \label{secFQR-T}

For each $n \ge 1$, the FQR-T control is based on two positive (activation) thresholds, $k^n_{1,2}$ and $k^n_{2,1}$ 
and the two queue-ratio parameters, $r_{1,2}$ and $r_{2,1}$
(which are chosen independent of $n$ under \eqn{ParamLim}). We
define two (centered) queue-difference stochastic processes
\bequ \label{DiffProc}
\bsplit
D^n_{1,2}(t) & \equiv Q^n_1(t) - k^n_{1,2} - r_{1,2}Q^n_2(t) \qandq \\
D^n_{2,1}(t) & \equiv r_{2,1}Q^n_2(t) - k^n_{2,1} - Q^n_1(t),\quad t \ge 0.
\end{split}
\eeq

As long as $D^n_{1,2}(t) < 0$ and $D^n_{2,1}(t) < 0$ we
consider the system to be not overloaded so that no customers
are routed to be served in the other class pool. Once one of
these inequalities is violated, the system is considered to be
overloaded, and sharing is initiated. For example, if
$D^n_{1,2}(t) \ge 0$, then class $1$ is judged to be overloaded
(because then $Q^n_1 - r_{1,2} Q^n_2 \ge k^n_{1,2}$), and it is
desirable to send class-$1$ customers to be served in pool $2$.
Note that $D^n_{1,2}(t) \ge 0$ does not exclude the case that
class $2$ is also overloaded; we can have $\beta_i(t) > 0$ for
{\em both} $i$. However, once one of the thresholds is crossed,
its corresponding class is considered to be ``more overloaded''
than the other class. (We refer to this situation as {\em
unbalanced overloads}.) We call $k^n_{1,2}$ and $k^n_{2,1}$
{\em activation thresholds}, because exceeding one of these
thresholds activates sharing (and not exceeding prevents
sharing when it is not desired).

The behavior of $X^n$ in \eqref{Xn} depends on the choice of
the thresholds $k^n_{i,j}$. In particular, we want the
thresholds to be large enough so that sharing will not take
place if both service pools are normally loaded, and to be
small enough to detect any overload quickly, and start sharing
in the correct direction once the overload begins. Note that
without sharing, the two pools operate as two independent
$M_t/M/m_t^n+M$ (time-varying Erlang-A) models. The familiar fluid and diffusion
limits for the stationary Erlang-A model give insight as to how to choose
these thresholds; e.g., see \cite{GMR02,PTW07}. In Assumption 2.4 of \cite{PW13a} and Assumption 3 of \cite{PW13b}
we assumed that the activation thresholds are chosen to
satisfy:
\bequ \label{thresholds}
k^n_{i,j}/n \ra 0 \qandq k^n_{i,j}/\sqrt{n} \ra \infty \qasq n \tinf, \quad i,j = 1,2
\quad \mbox{with}\, \, i \not= j.
\eeq
The first limit in
\eqref{thresholds} ensures that overloads are detected quickly
(immediately in the fluid model obtained as $n \tinf$), whereas
the second limit in \eqref{thresholds} ensures that stochastic
fluctuations of normally-loaded pools will not cause undesired
sharing, since the diffusion-scaled queue in that case are of
order $\sqrt{n}$.

Given that the system is designed so that sharing of customers
takes place only during overloads, it is reasonable to assume
that agents serve the other class customers (the so-called
``shared customers'') at a slower rate than they serve their own
designated customers. Thus, substantial sharing is likely to
reduce the effective service rate of the helping pool. In our
previous work we took measures to avoid sharing in both
directions simultaneously. In particular, we imposed the
one-way sharing rule described in \S \ref{secIntro}. However,
it is evident that the one-way sharing rule may considerably
slow the recovery after the overload is over. We elaborate in
\S \ref{secRelaxOneWay} below.

To remedy this problem, we could consider removing the one-way
sharing rule altogether and rely solely on the activation
thresholds to avoid undesired sharing. However, removing the one-way sharing rule makes it necessary to increase the activation thresholds
substantially, increasing the time until overloads
are detected. Moreover, if these thresholds are too large, then
some overloads may not be detected at all, because
abandonment keeps the queues from increasing indefinitely.
(While there is also a need to increase the activation thresholds in our setting here, that increase is less than
would be required if the one-way sharing was completely removed.)
Moreover, if sharing is taking place in one direction and then
immediately starts in the other direction in response to a
switch in the overload, then the combined service capacity of
both pools may be reduced significantly, creating a period of
severe congestion in both directions. Hence, it is beneficial
to avoid too much simultaneous two-way sharing. We again refer to Example 2 in \cite{PW09}.
Therefore, our new control relaxes the one-way sharing rule by introducing the release thresholds alluded to above.
We elaborate in the following subsection.

\subsection{The Proposed FQR-ART Control}

For the reasons discussed above, we suggest a modification of
the one-way sharing rule by introducing {\em release
thresholds} (RT). For each $n \ge 1$, we introduce two strictly
positive numbers $\tau^n_{1,2}$ and $\tau^n_{2,1}$. A newly
available type-$2$ agent is allowed to take a class-$1$
customer at time $t$ only if $Z^n_{2,1}(t) \le \tau^n_{2,1}$,
i.e., if the number of type-$1$ agents serving class-$2$
customers at the same time $t$ is below $\tau^n_{2,1}$ (and of
course $D^n_{1,2}(t) \ge 0$), and similarly in the other
direction.  (Ways to choose the parameters $\tau^n_{1,2}$ and $\tau^n_{2,1}$
will be discussed later.)

However, the new release thresholds allow small simultaneous sharing in both directions, which can slightly increase the overload.
In some cases, this slight increase in the overload is sufficient to cause the system to spin out of control and start to oscillate, as we
demonstrate in \S \ref{secOSC} below.
In particular, the new release thresholds make
activation thresholds satisfying \eqref{thresholds} unsuitable.
We therefore conclude that these activation thresholds should be positive in ``fluid
scale'', i.e., they should be chosen so as to satisfy
\bequ
\label{ThreshQ} \lim_{n\tinf} k^n_{i,j}/n = k_{i,j} > 0, \quad
i,j = 1,2.
\eeq

Thus, the FQR-ART control is specified by the parameter
six-tuple
$$(r_{1,2}, r_{2,1}, k^n_{1,2}, k^n_{2,1},\tau^n_{1,2}, \tau^n_{2,1})$$
and the routing and scheduling
rules which depend on the values of the two processes
$D^n_{i,j}$ and $Z^n_{i,j}$, $i \ne j$, in the manner described
above. Note that FQR-T requires knowing only the queue lengths
$Q^n_i (t)$ at each time $t$ (specifically, the values of the
two difference processes \eqref{DiffProc}), whereas FQR-ART
also requires knowledge of $Z^n_{1,2}$ and $Z^n_{2,1}$. Under
either control, the X model is a (possible inhomogeneous) CTMC.

\subsection{Analysis Via Fluid Approximations} \label{secFluidApprox}

Since the stochastic process $X^n$ in \eqref{Xn} under FQR-ART
is evidently too difficult to analyze exactly, we will employ a
deterministic dynamical-system approximation, and refer to that
approximation as ``fluid approximation'' or ``fluid model''
interchangeably. The main idea in using fluid approximations is
that, for large $n$, $\barx^n \approx x$, for some
deterministic function $x$ that is easier to analyze than the
untractable stochastic process $X^n$. (We use the `bar'
notation throughout to denote fluid scaled processes, e.g.,
$\barx^n \equiv X^n/n$.) In particular, the fluid counterpart
of $X^n$ in \eqref{Xn} is the six-dimensional deterministic
function \bequ \label{x} x \equiv x(t) \equiv (q_i(t),
z_{i,j}(t) : i, j = 1,2), \quad t \ge 0, \eeq where $q_i$ and
$z_{i,j}$ are the fluid approximations for the stochastic
processes $Q^n_i$ and $Z^n_{i,j}$, $i,j=1,2$. The approximation
$\barx^n \approx x$ should be supported by a {\em functional
law of large numbers} (FLLN), stating that $\barx^n \Rightarrow
x$ as $n \ra \infty$, extending \cite{PW13a}, but that remains
to be established.  (However, the FWLLN has been established for the FQR-T model in
\cite{PW12a}.)

In the stochastic system, customer routing depends on the
values of the difference processes in \eqref{DiffProc}. For
example, if sharing is taking place with pool $2$ helping class
$1$, and assuming $Z^n_{2,1} \le \tau^n_{2,1}$, the process
$D^n_{1,2}$ determines which customer class a newly available
type-$2$ agent will take.  as in \cite{PW11a,PW11b,PW13a}, that implies that the resulting fluid
model is much more complicated than most fluid models in the literature.
In particular,
 in the fluid system we
cannot simply replace the process $D^n_{1,2}$ with a process
\bes d_{1,2}(t) \equiv q_1(t) - k_{1,2} - r_{1,2}q_2(t), \quad
t \ge 0. \ees
In fact, the purpose of the control is to keep
$d_{1,2}(t) = 0$ during the overload. Hence, as in
\cite{PW11a,PW11b,PW13a}, a refined asymptotic analysis of the behavior of
$D^n_{1,2}$ (or $D^n_{2,1}$ during overloads in the other
direction) is required. That refined analysis can be carried out thanks to a stochastic averaging
principle, which replaces the processes $D^n_{i,j}$, $i,j = 1,2$, with the long-run average behavior of corresponding limiting stochastic processes.
In turn, those deterministic long-run averages determine the evolution of the fluid model;
see \S \ref{secFluidMain} below, where the fluid equations are developed.


\section{The Need to Relax the One-Way Sharing Rule} \label{secRelaxOneWay}

Relying on the fluid approximation, we now demonstrate why the
one-way sharing rule impedes recovery after the overload
incident is over. The simple fluid analysis suggests that
release thresholds provide a good remedy, and helps indicate
how they should be chosen.

\subsection{The Recovery Time With One-Way Sharing} \label{secDrawback}

We consider two consecutive time intervals $I_1 = [t_0, t_1)$
and $I_2 = [t_1, t_2)$ with $0 \le t_0 < t_1 < t_2 \le \infty$,
with the system being overloaded in opposite direction over
each interval. Suppose that class $2$ is overloaded over the
time interval $I_1$ and that sharing is taking place with pool
$1$ helping class $2$. Then, at time $t_1$ the loads suddenly
change in such a way that sharing is required in the other
direction. In particular, we assume that $\beta_1(t) \le 0$ and
$\beta_2(t) >  0$ for $t \in I_1$, whereas $\beta_1(t) > 0$ and
$\beta_2(t) \le  0$ for $t \in I_2$. We also assume that
$z_{2,1}(t_1) > 0$.

We do two different mathematical analyses.  We first consider a direct fluid model analysis,
and then afterwards we consider the stochastic system.
A fluid approximation for the evolution of $Z^n_{1,2}$ (which we refer to as $z_{1,2} (t)$) can
easily be derived using rate considerations. Since every
type-$1$ agent who is helping a class-$2$ customer at time $t >
t_1$ will finish service immediately after time $t$ at a rate
$\mu_{2,1}$, regardless of the value of $t$, due to the
memoryless property, and since there are no more class $2$
customers routed to pool $1$ after time $t_1$, we expect that
$z_{2,1}$ will satisfy the ODE
\bes \dot{z}_{2,1}(t) = -\mu_{2,1} z_{2,1}(t), \quad t \in I_2,
\ees
whose unique
solution is
\bequ \label{z21NoEmpty} z_{2,1} (t) = z_{2,1}(t_1)
e^{-\mu_{2,1}t}, \quad t \in [t_1,t_2). \eeq
As a consequence, for
the fluid model, if $z_{2,1}(t_1) > 0$, then pool $1$ will {\em
never} empty, so that sharing can {\em never} begin in the opposite
direction.

We now characterize the random time $T^n$ after the time $t_1$ in the stochastic system with scale $n$
for $Z^n_{2,1}(t)$ to first hit $0$.  The time required for all these customers to complete service is
the maximum of $Z^n_{2,1}(t_1)$ i.i.d. exponential random variables.
It is well known that the maximum of $n$ i.i.d. exponential random variables with mean $1$ is
the harmonic sum $H_n \equiv\sum_{j=1}^{n} (1/j)$.  Moreover, it is well known that
$H_n - \log_{e}{n} \ra \gamma$ as $n \ra \infty$, where $\gamma \equiv 0.57721 \ldots$ is the
Euler-Mascheroni constant.  This limit is relevant for us, because
from the established FWLLN in \cite{PW13a}, we know that having
$z_{2,1}(t_1) > 0$ implies that $Z^n_{2,1}(t_1) \approx z_{2,1}
(t_1) n$.

Hence, given $Z^n_{2,1}(t_1)$ and its approximate value, for large $n$,
\bequ \label{Z12ln}
E[T^n] = \sum_{j=1}^{Z^n_{2,1}(t_1)}{\frac{1}{j\cdot\mu_{2,1}}} \approx
\frac{\log_{e}{(Z^n_{2,1}(t_1))}}{\mu_{2,1}} \approx \frac{\log_{e}{(n
z_{2,1} (t_1))}}{\mu_{2,1}}.
\eeq
We thus see that
the expected time required for a pool to empty its shared
customers after an overload is over, and no new shared
customers are routed to that pool,
 is of order $\log_{e}(n)$ as $n \ra \infty$.


\subsection{Choosing Appropriate Release Thresholds} \label{secLower}

The simple considerations leading to \eqref{z21NoEmpty} and \eqref{Z12ln}
 show that a large system will be slow to
recover after an overload is over. That analysis also helps
choose appropriate release thresholds. Indeed, the fluid model
easily generates an approximate recovery time. In particular,
if a release threshold of $\tau_{2,1}$ is used in the fluid
model starting with $z_{2,1}(t_1)$ at time $t_1$, where
$z_{2,1}(t_1) > \tau_{2,1} > 0$, then the release threshold
will be hit at time \bes T \equiv
\frac{1}{\mu_{2,1}}\log_{e}\left(\frac{z_{2,1}(t_1)}{\tau_{2,1}}\right).
\ees

The analysis above indicates that the release thresholds in
stochastic system $n$ should be of order $O(n)$ as $n$
increases. It suffices to pick two strictly positive numbers
$\tau_{1,2}$ and $\tau_{2,1}$ and let \bequ
\label{ReleaseScaling} \tau^n_{1,2} \equiv n \tau_{1,2} \qandq
\tau^n_{2,1} \equiv n \tau_{2,1}. \eeq With the scaling in
\eqn{ReleaseScaling}, the recovery time $T^n$ in system $n$
should be approximately a constant, independent of $n$.

In summary, with FQR-ART, an available type-$2$ agent is
allowed to serve a class-$1$ customer only if $Z^n_{2,1}(t) \le
\tau^n_{2,1}$ (or, equivalently, only if $\barz^n_{2,1} (t) \le
\tau_{2,1}$), and of course $D^n_{1,2}(t) \ge 0$, and similarly
in the other direction. The choice in \eqn{ReleaseScaling}
shows that the release thresholds should be proportional to
$n$, but does not determine the proportionality constants
$\tau_{1,2}$ and $\tau_{1,2}$. Further analysis shows that
these can be quite small, as we show next.

\subsection{Simulation Experiments} \label{secSimRelease}

\begin{figure}[h!]
  \hfill
  \begin{minipage}[t]{.4\textwidth}
    \begin{center}
      \includegraphics[width=6cm]{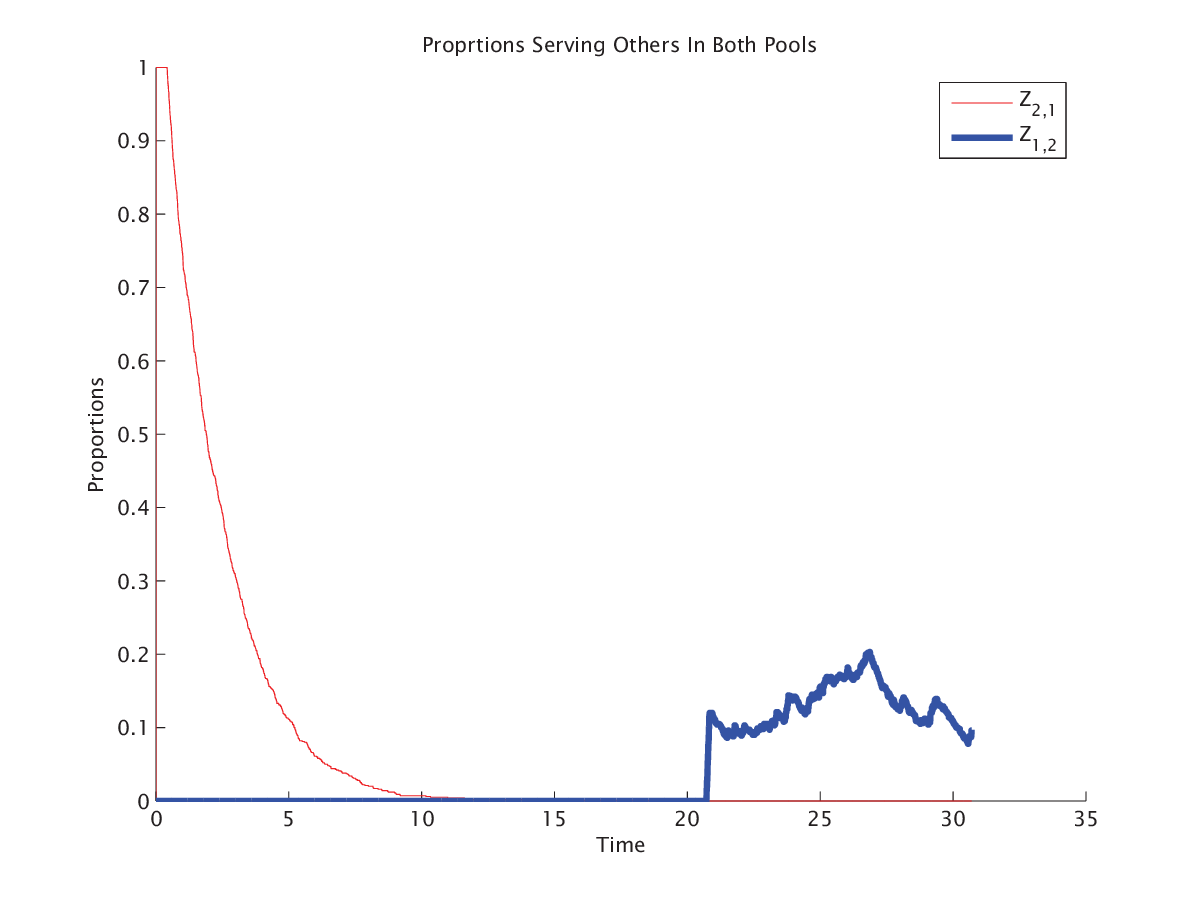}
      \caption{Sample paths of $\barz^n_{1,2} (t)$ and $\barz^n_{2,1} (t)$ initialized incorrectly, without release thresholds.}
      \label{fig2}
    \end{center}
  \end{minipage}
  \hfill
  \begin{minipage}[t]{.45\textwidth}
    \begin{center}
      \includegraphics[width=6cm]{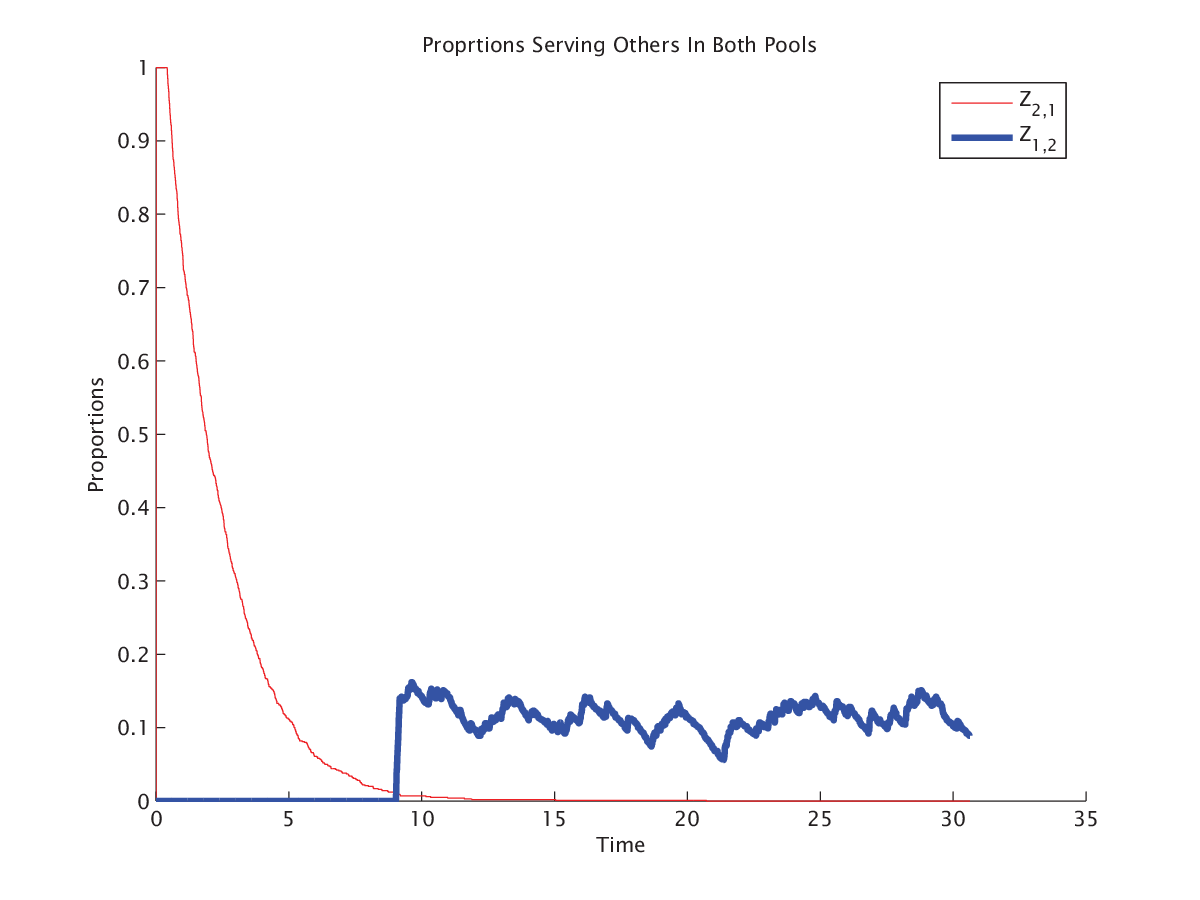}
      \caption{Sample paths of $\barz^n_{1,2} (t)$ and $\barz^n_{2,1} (t)$ initialized incorrectly, with release thresholds $\tau_{1,2} = \tau_{2,1} = 0.01$.}
      \label{fig3}
    \end{center}
  \end{minipage}
  \hfill
\end{figure}

To illustrate the importance of the release thresholds for
stochastic systems, we conducted simulation experiments,
comparing the performance of a system with and without release
thresholds. The results can be seen in Figures \ref{fig2} and
\ref{fig3}.

The (fixed) parameters for this simulation are
\beas
m^n_1 & = & m^n_2 = 1000, \quad \lm^n_1 = 1200, \quad \lm^n_2 = 990, \quad \mu_{1,1} = \mu_{2,2} = 1, \\
\mu_{1,2} & = & \mu_{2,1} = 0.5, \quad \kappa^n_{1,2} = \kappa^n_{2,1} = 100, \qandq r_{1,2} =
r_{2,1} = 1.
\eeas
(Here, we can think of $n$ as being fixed
and equal to $1000$.) With these parameters, $\rho^n_1 = 1.2$
and $\rho^n_2 = 0.99$, where $\rho^n_i \equiv
\lm^n_i/(m_i^n\mu_{i,i})$, so that class 1 may be regarded as
overloaded, whereas class $2$ may be regarded as normally
loaded (recall \eqref{beta}).

To respond to that unbalanced overload by having pool $2$ help
class $1$, we should have $Z^n_{1,2} > 0$ and $Z^n_{2,1} = 0$
if one-way sharing is employed. However, we initialize the
system at time $0$ sharing in the opposite direction, with {\em
all} pool 1 agents serving class 2 customers. We are interested
in the time it takes the stochastic process $Z^n_{2,1}$ to
reach $0$, so that the desired sharing can begin. Without
release thresholds, the required recovery time is quite long,
approximately $21$ (mean service times, of their own type). In
contrast, with release thresholds of only $\tau^n_{1,2} =
\tau^n_{2,1} = 0.01n = 10$, that time is reduced from about
$21$ to about $9$ service times. Thus, clearing the last $1\%$
of the class-$2$ customers in pool $1$ without release
thresholds takes more than half the total clearing time!

We hasten to admit that we just considered an extreme example
in which {\em all} of service pool $1$ is initially busy with
customers from class $2$.
 We did so in order to convey the message that
{\em it is the last few agents working with class $1$ that
cause the largest part of the delayed response}. In particular,
the $Z^n_{2,1}$ process decreases fast at the beginning, but
then the decrease rate slows down considerably.

From Figures \ref{fig2} and \ref{fig3}, it is also easy to see
what happens in less extreme cases, when $0 < Z_{2,1}(0) <
m_1$.
For example, if we initialize with $20\%$ sharing in the wrong direction, 
we see that, without a release threshold, the time to activate
sharing in the right direction is about $21 - 4 = 17$ time
units. In contrast, with release thresholds, it is about $9 - 4
= 5$ time units.  (Figures \ref{fig2} and \ref{fig3} show that the common value $4$ in these calculations is the time to go from $100\%$ sharing in the wrong direction to
only $20\%$ sharing in the wrong direction, which would be the same in the two cases.)
When we start with a lower percentage of
agents sharing the wrong way, the difference becomes even more
dramatic, because we eliminate a common initial period (here of
length $4$ time units).


\section{Congestion Collapse Due to Oscillations}\label{secOSC}

The previous section dramatically showed the need for the
release thresholds when the direction of the overload suddenly
shifts. However, a more common case is for the two systems to
simply return to normal loading, after which no sharing in
either direction is desired. We now show that the release
thresholds can cause serious problems when the system returns
to normal loading after an overload incident
if the activation thresholds are too small. 
In this case, there is a potential difficulty when the
inefficient sharing condition holds, i.e., when
$\mu_{1,1} > \mu_{2,1}$ 
and $\mu_{2,2} > \mu_{1,2}$,
which is what we now assume.  We show that, with inefficient
sharing, the release thresholds combined with small activation
thresholds can lead to oscillatory poor performance. We
emphasize that, even though the performance is oscillatory, the
model after the overload is over is a (necessarily aperiodic)
positive-recurrent and stationary time-homogeneous CTMC when there is abandonment (as discussed in \S \ref{secCollapse}).
In particular, our examples below are time-homogeneous CTMCs because the arrival rates and staffing levels are kept fixed.

\subsection{Simulations of Oscillating Systems with Inefficient Sharing} \label{appOSC}

The oscillatory behavior is more evident when there is no
abandonment, so we start by considering a system without
abandonment. We start with an extreme case having very
inefficient sharing; i.e., we let $\mu_{1,1} = \mu_{2,2} = 1$,
but $\mu_{1,2} = \mu_{2,1} = 0.1$. Afterwards we consider a
more realistic example with customer abandonment and less
efficiency loss from sharing. We consider a relatively heavily
loaded symmetric system. In particular, let there be $m_1 = m_2
= 100$ agents in each pool and let the arrival rates be $\lm_1
= \lm_2 = 98$. Thus each class alone is stable, but if all the agents are busy in one pool with
$n$ serving the other class, then the total service rate out is
$0.1 k + (100 -k) = 100 - 0.9k$.  When $k \ge 3$, the maximum service rate is
less than the arrival rate $98$, so that the rate in exceeds
the maximum rate out at that instant.

We now illustrate bad behavior for poorly chosen thresholds.
We make both the activation thresho.lds and the release thresholds be too small.
In particular, we use
ratio parameters $r_{1,2} = r_{2,1} = 1$, activation
thresholds $k^n_{i,j} = 10$ and release thresholds
$\tau^n_{i,j} = 1$ for $i,j = 1,2$ and $i \not= j$. We start
the system with both pools busy serving their own class, but no
queues, i.e., $Z^n_{1,1} (0) = Z^n_{2,2} (0) = 100$ and $Q^n_1
(0) = Q^n_2 (0) = 0$.
 The symmetry implies that both pools and
queues exhibit symmetric behavior.
\begin{figure}[h!]
  \hfill
  \begin{minipage}[t]{.4\textwidth}
    \begin{center}
      \includegraphics[width=6cm]{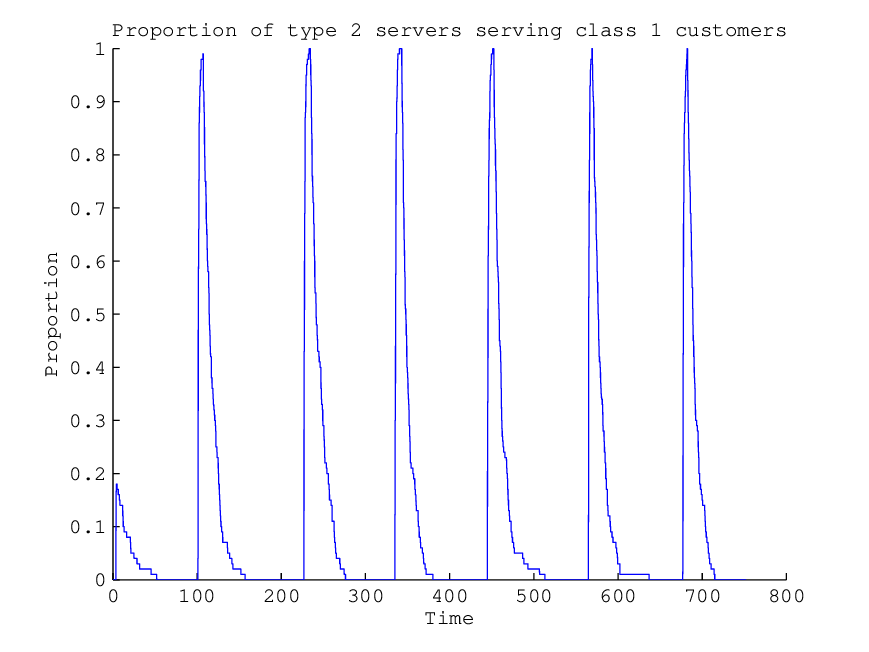}
      \caption{Oscillations of $\barz^n_{1,2}$ in the extreme symmetric example with $\tau^n_{i,j} = 1$, $k^n_{i,j} = 10$ and no abandonment.}
      \label{figZ12oscNo1}
    \end{center}
  \end{minipage}
  \hfill
  \begin{minipage}[t]{.4\textwidth}
    \begin{center}
      \includegraphics[width=6cm]{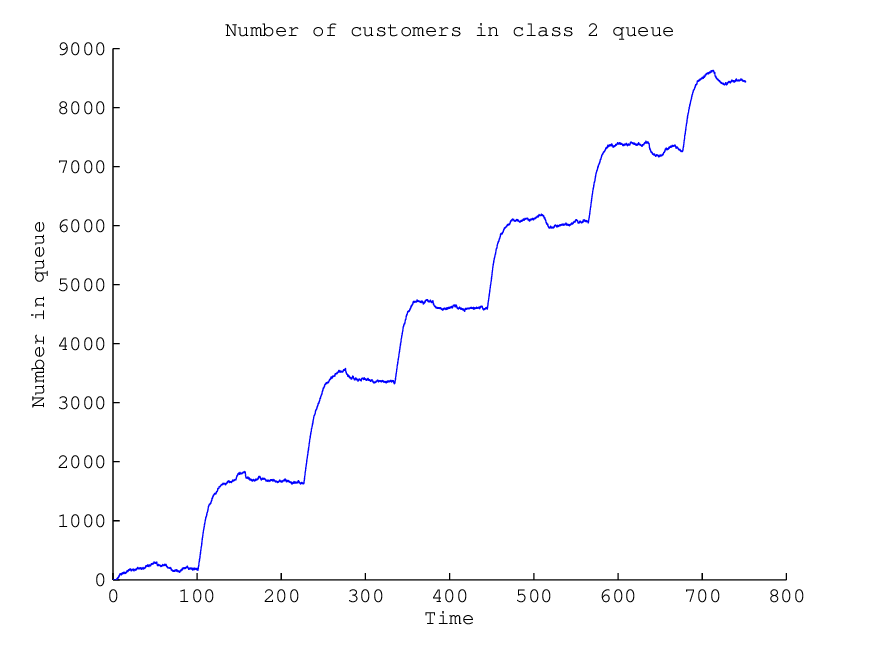}
      \caption{Oscillating growth of $\barq^n_2$ in the extreme symmetric example with $\tau^n_{i,j} = 1$, $k^n_{i,j} = 10$ and no abandonment.}
      \label{figQ2oscNo1}
    \end{center}
  \end{minipage}
  \hfill
\end{figure}

Figures \ref{figZ12oscNo1} and \ref{figQ2oscNo1} show a single
simulated sample path.
Figure \ref{figZ12oscNo1} shows that the proportion of agents
serving customers from the other pool oscillates between $0$
and $1$, alternating between these two extremes over this
horizon.
 The oscillatory behavior is occurring despite the fact
that there is no sharing initially. Figure \ref{figQ2oscNo1}
shows that the queue lengths are growing in an oscillating
manner over the time interval $[0,800]$ at an average rate of
$10$.

The oscillatory behavior also occurs for systems with
abandonment, but it is often hard to detect, because the
abandonment ensures that the stationary stochastic system after
the overload has ended is stable and it dampens any oscillatory
behavior. Nevertheless, the difficulty highlighted above
remains with abandonment.

\begin{figure}[h!]
  \hfill
  \begin{minipage}[t]{.4\textwidth}
    \begin{center}
      \includegraphics[width=6cm]{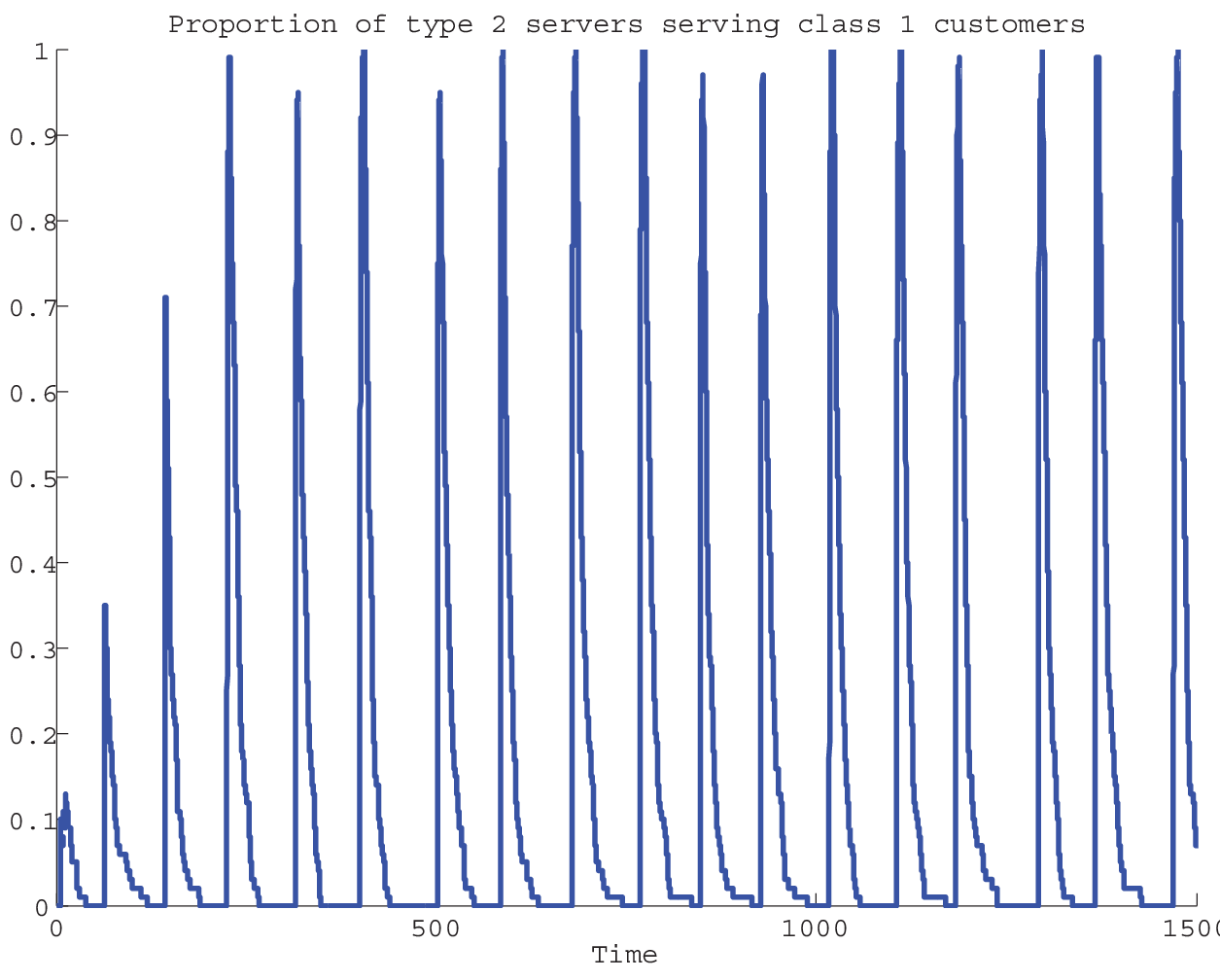}
      \caption{Oscillations of $\barz^n_{1,2}$ in the extreme symmetric example with $\tau^n_{i,j} = 1$, $k^n_{i,j} = 10$ with abandonment.}
      \label{figZ12oscAb}
    \end{center}
  \end{minipage}
  \hfill
  \begin{minipage}[t]{.4\textwidth}
    \begin{center}
      \includegraphics[width=6cm]{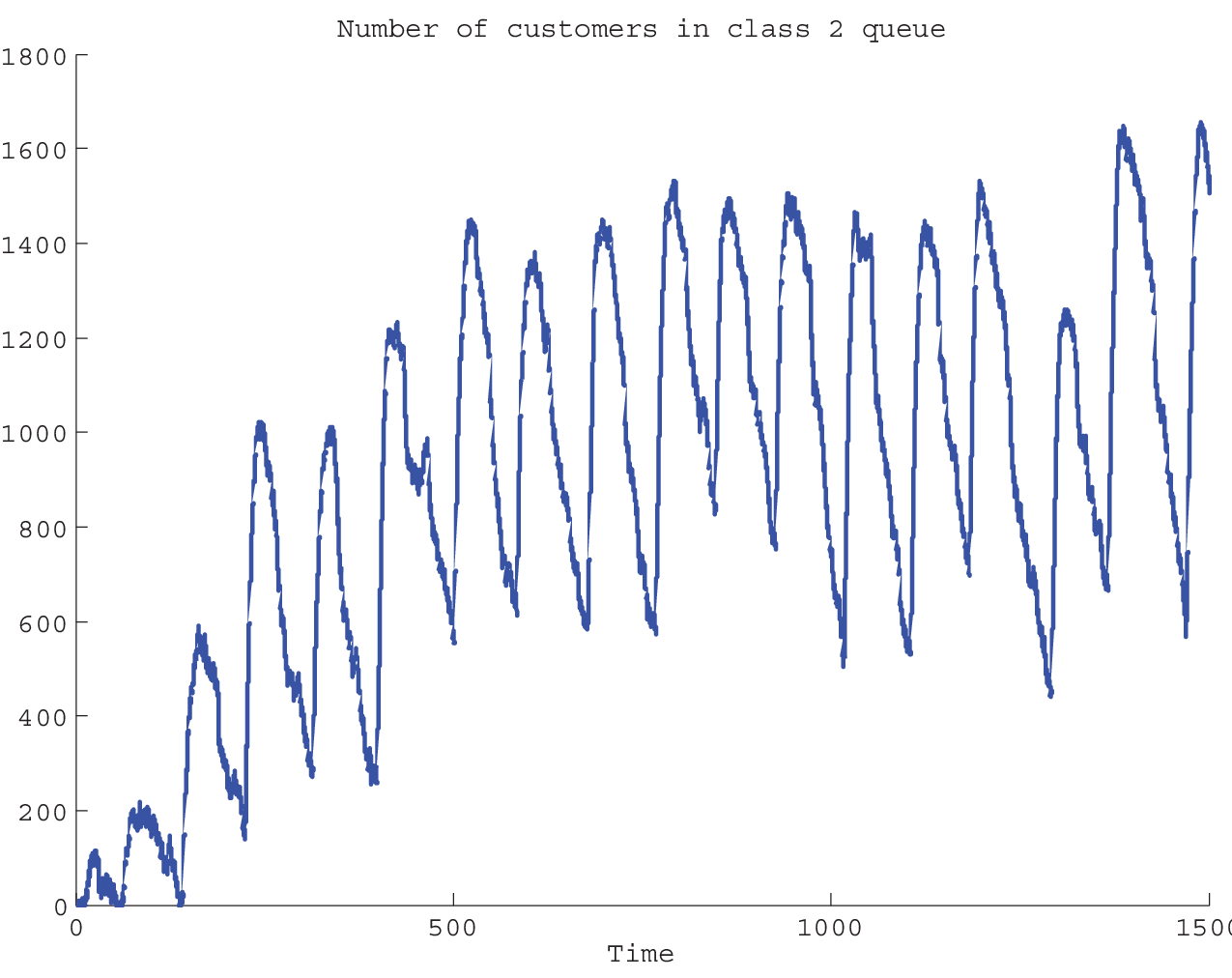}
      \caption{Oscillating stable behavior of $\barq^n_2$ in the extreme symmetric example with $\tau^n_{i,j} = 1$, $k^n_{i,j} = 10$ with abandonment.}
      \label{figQ2oscAb}
    \end{center}
  \end{minipage}
  \hfill
\end{figure}

To demonstrate dramatically, we simulated the same system
considered in the previous example, but now with the low
positive abandonment rates $\theta_1 = \theta_2 = 0.01$.
Figures \ref{figZ12oscAb} and \ref{figQ2oscAb} show that the
oscillatory behavior remains.
Moreover, Figure \ref{figQ2oscAb} suggests that $Q^n_2$ (and,
by symmetry, also $Q^n_1$) stabilizes at an overloaded
oscillatory equilibrium. The oscillatory behavior in Figures
\ref{figZ12oscAb}--\ref{figQ2oscAb} may be surprising at first,
because the underlying (time-homogeneous) CTMC after the
overload has ended is ergodic, as we mentioned above.
Fortunately, the fluid model provides valuable insight, as we
explain in \S \ref{secInsight}.

\begin{figure}[h!]
  \hfill
  \begin{minipage}[t]{.4\textwidth}
    \begin{center}
      \includegraphics[width=6cm]{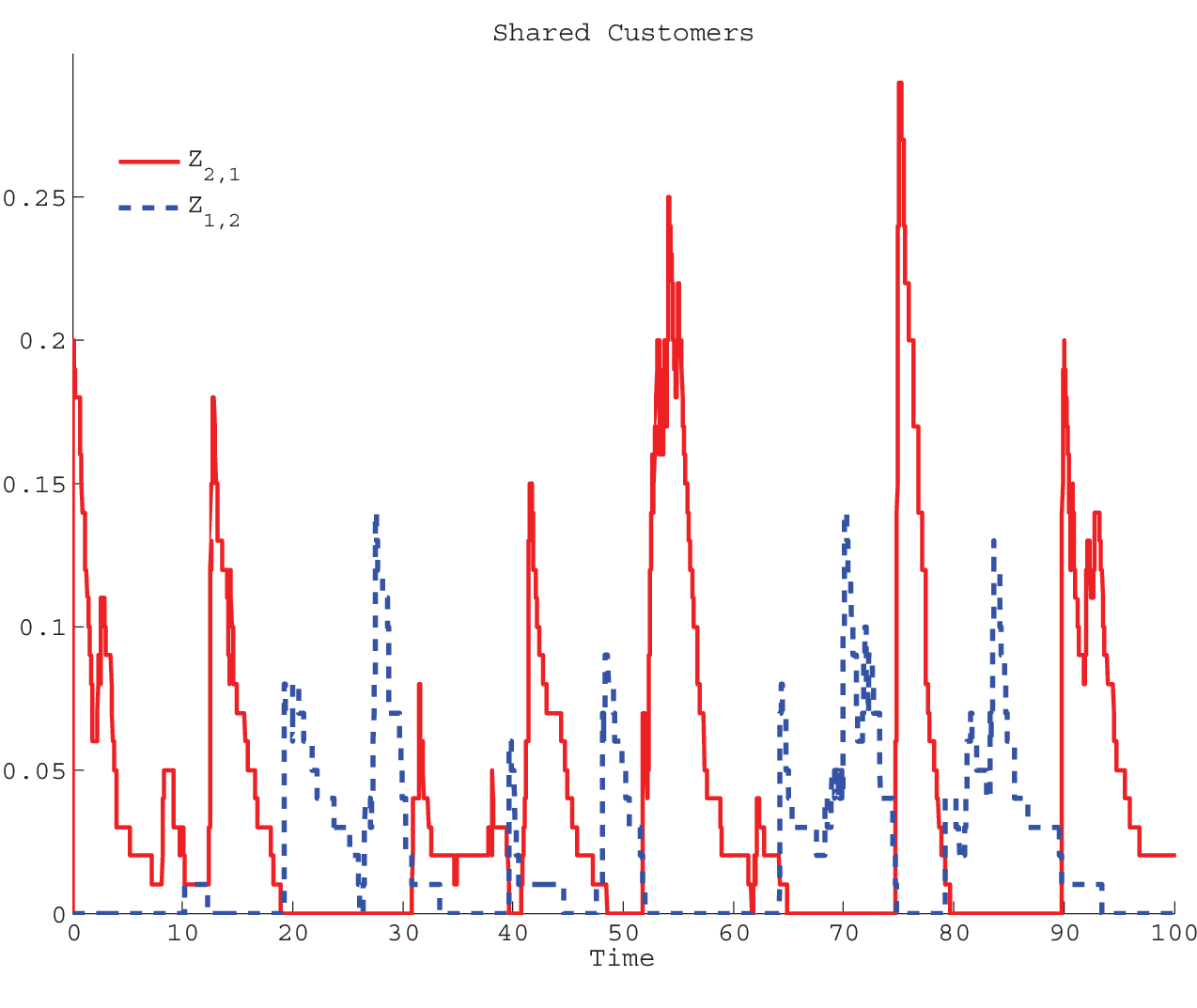}
      \caption{Oscillations of $\barz^n_{1,2}$ in the more realistic symmetric example with abandonment: $\mu_{1,1} = 1$, $\mu_{1,2} = 0.5$,
      $\theta_{1} = 0.5$, $\tau^n_{1,2} = 1$ and $k^n_{i,j} = 10$.}
      \label{figZosc}
    \end{center}
  \end{minipage}
  \hfill
   \begin{minipage}[t]{.4\textwidth}
    \begin{center}
      \includegraphics[width=6cm]{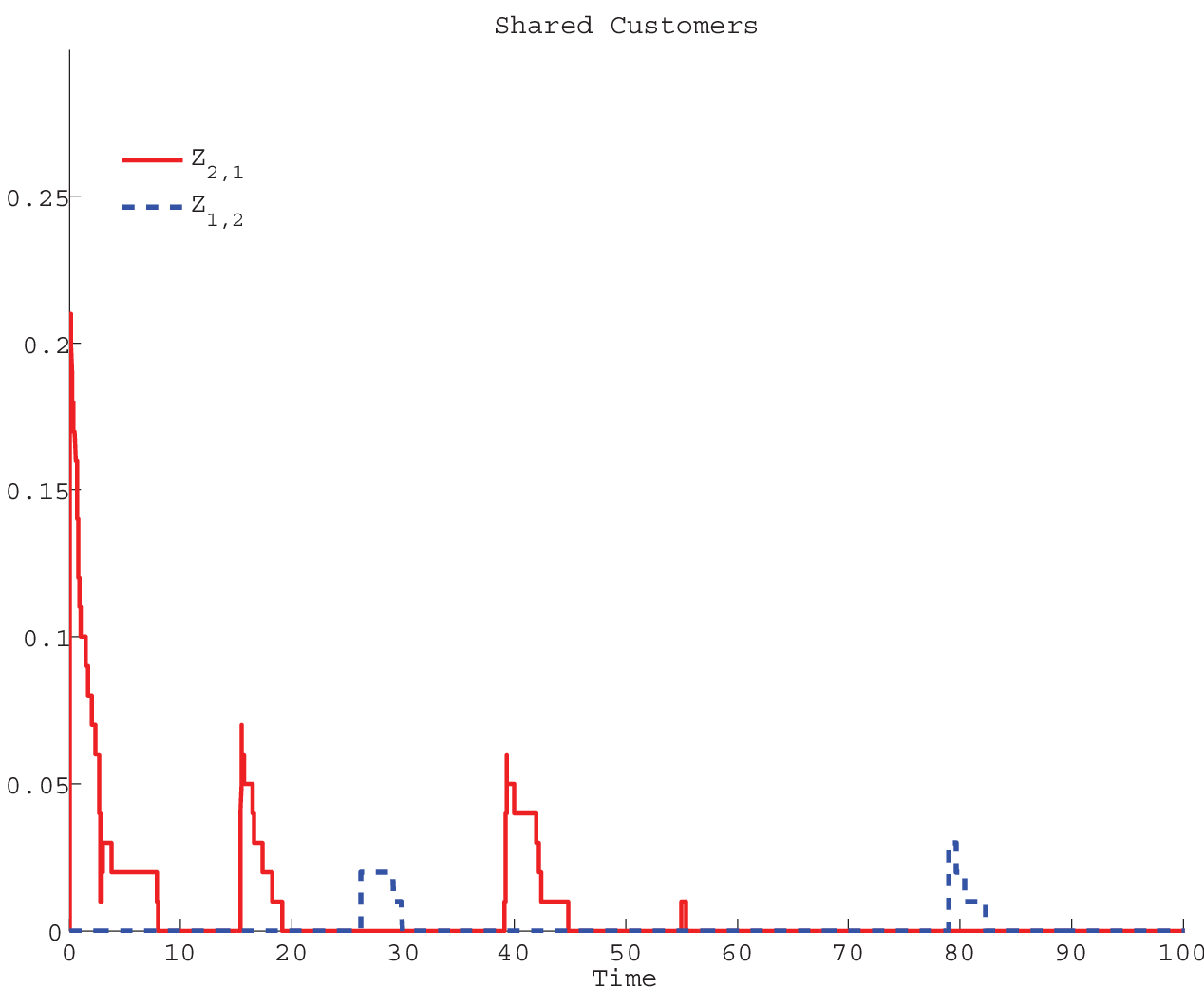}
      \caption{$\barz^n_{1,2}$ in the more realistic example 
      but with higher activation thresholds $k^n_{i,j} = 35$.}
      \label{figZnoOsc}
    \end{center}
  \end{minipage}
  \hfill
\end{figure}

We now consider a less-extreme more realistic example, in which
the sharing service rates and abandonment rates are changed to
$\mu_{1,2} = \mu_{2,1} = \theta_1 = \theta_2 = 0.5$.
First, Figure \ref{figZosc} shows the proportion of shared
customers over time with the previously specified activation
thresholds of $k^n_{i,j} = 10$, but we now consider a system
that is recovering from an overload in which pool $1$ was
helping class $2$ customers. In particular, there are initially
$20$ type-$1$ agents helping class-$2$ customers. By taking
this initial condition, we are considering a system that starts
``worse off'' than before, because it is initially overloaded.
(In the other two examples, the systems were initialized
empty.) We consider the time interval $[0,100]$ to make the
figures clear, but the behavior shown in the figures below
remained for the whole duration of the simulation (which lasted
for 1500 time units).

In this case, substantial customer abandonment significantly
dampens the sharing oscillations seen previously. Nevertheless,
Figure \ref{figZosc} shows that the pools share repeatedly in
an oscillating manner over the time interval $[0,100]$.
Although the long-run average number of agents that are helping
the other class is not significant, this oscillatory behavior,
is clearly undesirable. We do not show figures of the queues
because they are uninformative (the oscillations are
insignificant). Hence, the bad behavior in a system with a
relative substantial customer abandonment may be hard to detect
by only observing the queues, so that a system with no
abandonment, or low abandonment rate, gives important insights.

To remedy the problem in Figure \ref{figZosc}, we propose
increasing the activation thresholds. To illustrate the
potential benefit,
 Figure \ref{figZnoOsc} shows the sharing when the activation thresholds are increased to $k^n_{i,j} = 35$, $i,j = 1,2$,
with all other parameters kept the same. Even though some
customers are shared occasionally, especially just after the
overload is over, the oscillatory behavior is minimal and
decays quickly.


\subsection{Insight from the Deterministic Fluid Model}\label{secInsight}

In the examples we have just considered, the six-dimensional stochastic process $X^n$ in \eqn{Xn} describing the system performance
after the overload incident has ended is a stationary CTMC.
With customer abandonment,
that CTMC is necessarily stable, so that with FQR-ART and any parameter setting,
the stochastic process $X^n$ in \eqn{Xn} necessarily has a unique
steady-state distribution.  Nevertheless, we have just seen that the system can exhibit quite complex undesirable behavior for some initial conditions
if the control parameters are not set properly.

Fortunately, the fluid model we develop provides an effective means to
study the complex system performance and set the control parameters.
The oscillating behavior we see in the simulations looks periodic, but it is not quite; it is nearly periodic, just as in \cite{LW11}.
The system becomes more nearly periodic as the scale increases.
In the many-server heavy-traffic limit, the stochastic process $X^n$ approaches the deterministic solution of the fluid model we introduce next
to serve as an approximation.  From the algorithm for that fluid model, we see that it possesses a periodic equilibrium for
some initial conditions.

As a consequence, the fluid model can be bistable; it can have a periodic equilibrium in addition to a
stable equilibrium, depending on the initial conditions.  Consequently, the order in which two different limits occur
leads to different stories.  As time increases, for any fixed scale, the stochastic process approaches its unique steady-state
distribution.  In contrast, as the scale increases, a properly-scaled version of the stochastic process approaches a deterministic function, which can be periodic.
Thus, the fluid model provides important insight: an oscillatory fluid approximation implies that
a corresponding large system experiences oscillatory behavior for prohibitively large time intervals, even though it is essentially a stationary CTMC.
In \cite{PW14} we prove that the fluid models can exhibit this bi-stability;
Here it is verified numerically by applying the fluid algorithm.

\section{The Fluid Model} \label{secFluidMain}

The fluid model approximating the stochastic system $X^n$ under FQR-ART is described as the solution to an ordinary differential equation (ODE),
but that ODE depends on a stochastic averaging principle (AP).
In this section we derive that ODE via a heuristic representation of the inhomogeneous CTMC in \eqref{Xn}.
The reasoning in the justification of the fluid model approximation parallels the heuristic engineering discussion in \cite{PW11a}, to which we refer for more discussion.
For mathematical support for that reasoning, see \cite{PW11b,PW13a}.

\subsection{Representation of the Stochastic System During Overloads} \label{secRep}

The sample paths of the queueing system can be represented in terms of its primitive processes, i.e., the arrival, abandonment
and service processes, as a function of the control.
Unlike traditional fluid models, in which the primitive stochastic processes are replaced by their long-run rates,
the deterministic fluid model here is more involved and includes a stochastic ingredient in the form of a stochastic AP, which we describe in detail
in \S \ref{secAP} below.

Even though we are not proving that the fluid model arises as a weak limit of the fluid-scaled stochastic system, we need to take asymptotic
considerations in order to develop the fluid approximation.
We thus start with a representation of the stochastic system during overloads, assuming that both service pools are full over an interval $[0, T]$,
i.e.,
\bequ \label{eqZinitial}
Z^n_{1,1}(t) + Z^n_{2,1}(t) = m^n_1(t) \qandq Z^n_{2,2}(t) + Z^n_{1,2}(t) = m^n_2(t), \quad t \in [0, T].
\eeq
During the time interval $[0,T]$ no customers can enter service immediately upon arrival, and so all customers are delayed in queue.
For simplicity, we first consider intervals over which the staffing functions are continuous and differentiable everywhere.
In \S \ref{secNum} we give an example of a staffing function with discontinuity; see Figure \ref{figVaryPool} below.

We represent the sample paths of $X^n$ as random time-changes of independent unit-rate Poisson processes,
as reviewed in \cite{PTW07}; see Equations (41)-(43) in \cite{PW13a} for such a representation applied to the X model operating under FQR-T.
Let
\bequ \label{mA}
\bsplit
\mA^n_{1,2}(s) & \equiv \{\{D^n_{1,2}(s) > 0\} \cap \{Z^n_{2,1}(s) \le \tau^n_{2,1}\}\} \qandq \\
\mA^n_{2,1}(s) & \equiv \{\{D^n_{2,1}(s) > 0\} \cap \{Z^n_{1,2}(s) \le \tau^n_{1,2}\}\},
\end{split}
\eeq
the representation of $Q^n_1$ over $[0, T]$ is
\bes
\bsplit
Q^n_1(t) & = N^a_{1} \left(\int_0^t \lm^n_1(s) ds\right) - N^u_1 \left(\theta_1 \int_0^t Q^n_1(s) ds \right) \\
& \quad - N^+_{1} \left(\int_0^t \1_{\mA^n_{1,2}(s)} \left(\mu_{1,1} Z^n_{1,1}(s) + \mu_{1,2} Z^n_{1,2}(s) + \mu_{2,1} Z^n_{2,1}(s) + \mu_{2,2} Z^n_{2,2}(s)\right) ds \right) \\
& \quad - N^-_1 \left(\int_0^t (1 - \1_{\mA^n_{1,2}(s)} - \1_{\mA^n_{2,1}(s)}) \left(\mu_{1,1} Z^n_{1,1}(s) + \mu_{2,1} Z^n_{2,1}(s) \right) ds\right),
\end{split}
\ees
where $N^a_1, N^u_1, N^+_1$ and $N^-_1$ are mutually independent unit rate (homogeneous) Poisson processes,
and $\1_{A}$ is the indicator function that is equal to $1$ if event $A$ occurs, and to $0$ otherwise.

Note that the representation of $Q^n_1$ is essentially a flow conservation equation (based on the memoryless property of the exponential distribution).
That is, the queue at time $t$ is all those customers who arrived by that time, captured by the Poisson process $N^a_1$, minus all the customers that abandoned,
captured by the Poisson process $N^u_1$, minus all those who were routed into service, as captured by the last two Poisson processes in the expression.
Similar expressions hold for the other processes in $X^n$.

We elaborate on how the intensities of the last two Poisson processes in the right-hand side (RHS) of the representation were obtained.
First, if at time $s \in [0,T]$ the event $\mA^n_{1,2}(s)$ in \eqref{mA} holds,
then any newly available agent in the system will take his next customer
from the head of queue $1$. Since agents become available at an instantaneous rate
$\sum_{i,j} \mu_{i,j}Z^n_{i,j}(s)$ at time $s$, we get the third component in the RHS of $Q^n_1(t)$.
Next we recall that, by the routing rule of FQR-ART, if
at a time $s \in [0, T]$ $\mA^n_{2,1}(s)$ in \eqref{mA} holds,
then any newly available agent takes his next customer from queue $2$, in which case queue $1$ will not decrease due to a service completion.
If neither of the events $\mA^n_{1,2}(s)$ or $\mA^n_{2,1}(s)$ holds at a time $s$, then only service completions at pool $1$ will cause a decrease
at queue $1$ due to a customer from that queue being routed to service. That explains the last term in the RHS of the representation.

Next, we exploit the fact that each of the Poisson processes in the representation minus its random intensity constitutes a martingale
(again, see \cite{PTW07,PW13a}), e.g.,
$$M^{n,u}_1 \equiv N^u_1 \left(\theta_1 \int_0^t Q^n_1(s) ds \right) - \theta_1 \int_0^t Q^n_1(s) ds$$ is a martingale.
Thus, subtracting and then adding all the random intensities, and using the fact that a sum of martingales is again a martingale,
we get the following
representation for the processes $Q^n_1, Q^n_2, Z^n_{1,2}, Z^n_{2,1}$
(the remaining two processes $Z^n_{1,1}$ and $Z^n_{2,2}$ are determined by \eqref{eqZinitial}):

\bequ \label{RepXstoc}
\bsplit
Q^n_1(t) & = M^n_1(t) + \int_0^t \lm^n_1(s) ds - \int_0^t\theta_1 Q^n_1(s)ds \\
& \quad - \int_0^t \1_{\mA^n_{1,2}(s)}
\left(\mu_{1,1} Z^n_{1,1}(s) + \mu_{1,2} Z^n_{1,2}(s) + \mu_{2,1} Z^n_{2,1}(s) + \mu_{2,2} Z^n_{2,2}(s)\right)ds \\
& \quad - \int_0^t (1 - \1_{\mA^n_{1,2}(s)} - \1_{\mA^n_{2,1}(s)})
\left(\mu_{1,1} Z^n_{1,1}(s) + \mu_{2,1} Z^n_{2,1}(s) \right) ds, \\
Q^n_2(t) & = M^n_2(t) + \int_0^t \lm^n_2(s) ds - \int_0^t \theta_2 Q^n_2(s) ds \\
& \quad - \int_0^t \1_{\mA^n_{2,1}(s)}
\left(\mu_{1,1} Z^n_{1,1}(s) + \mu_{1,2} Z^n_{1,2}(s) + \mu_{2,1} Z^n_{2,1}(s) + \mu_{2,2} Z^n_{2,2}(s)\right) ds \\
& \quad - \int_0^t (1 - \1_{\mA^n_{1,2}(s)} - \1_{\mA^n_{2,1}(s)})
\left(\mu_{2,2} Z^n_{2,2}(s) + \mu_{1,2} Z^n_{1,2}(s) \right) ds, \\
Z^n_{1,2}(t) & = M^n_{1,2}(t) + \int_0^t \1_{\mA^n_{1,2}(s)} \mu_{2,2}Z^n_{2,2}(s) ds - \int_0^t (1 - \1_{\mA^n_{1,2}(s)})\mu_{1,2}Z^n_{1,2}(s) ds, \\
Z^n_{2,1}(t) & = M^n_{2,1}(t) + \int_0^t 1_{\mA^n_{1,2}(s)} \mu_{1,1}Z^n_{1,1}(s) ds - \int_0^t (1 - \1_{\mA^n_{2,1}(s)}) Z^n_{2,1}(s)) ds, \\
\end{split}
\eeq
where $M^n_1, M^n_2, M^n_{1,2}$ and $M^n_{2,1}$ are the martingale terms alluded to above.
It is not hard to show that those martingales are negligible in the fluid scaling
(divided by $n$, e.g., $\bar{M}^n_i \equiv n^{-1} M^n_i$), i.e., that $\bar{M}^n_i \Ra 0$ and $\bar{M}^n_{i,j} \Ra 0$ as $n \tinf$,
uniformly over $[0, T]$, $i,j = 1,2$; see, e.g., Lemma 6.1 in \cite{PW13a}.
Hence, we consider those martingales as a negligible stochastic noise that can be ignored for the purpose of developing the fluid approximation for \eqref{RepXstoc}.

To replace the stochastic integral representation in \eqref{RepXstoc} with a deterministic one,
we need to replace the indicator functions with smooth functions.
This is where the AP comes in.
What we do is replace the term $1_{\{D^n_{1,2} (t) > 0\}}$ by the steady state probability that
an associated fast-time scale process (FTSP) is greater than or equal to $0$, denoted by $\pi_{1,2} (x(t))$, which is a function of the fluid state at time $t$, $x(t)$.
Both $x(t)$ and $\pi(x(t))$ turn out to be a continuous function of $t$.  This complicated step requires more explanation and justification,
which again was the subject of \cite{PW11a,PW11b,PW13a}.

We give a brief account in the rest of this section and the following one.
We start by assuming that there is a fluid counterpart $x$ for $X^n$ in \eqref{RepXstoc} which is {\em continuous and differentiable}.
(This fact can be shown to hold by a minor modification of Corollary 5.1 in \cite{PW13a}).
For any fluid point $x(t)$, let
\bequ \label{dij}
d_{1,2}(x(t)) \equiv q_1(t) - r_{1,2}q_2(t) - k_{1,2} \qandq d_{2,1}(x(t)) \equiv r_{2,1}q_2(t) - q_1(t) - k_{2,1}.
\eeq
We first observe that, if $d_{i,j}(x(t)) > 0$ then, since $d_{i,j}(\cdot)$ is a continuous function, $d_{i,j}$ is strictly positive over an interval,
and similarly if $d_{i,j} < 0$, $i,j = 1,2$. In such cases the indicator functions are easy to deal with because each is a constant over the
interval, and equals either $1$ or $0$. For example, if $d_{1,2}(x(t)) > 0$ for $t \in [s_1, s_2)$, for some $0 \le s_1 < s_2 < \infty$,
and in addition, $Z^n_{2,1}(t) \le \tau^n_{2,1}$ over that interval for all $n$ large enough, then
\bes
1_{\mA^n_{1,2}(t)} \equiv \1_{\{\{D^n_{1,2}(t) > 0\} \cap \{Z^n_{2,1}(t) \le \tau^n_{2,1}\}\}} = \1_{\{[s_2, s_2)\}}(t) \quad \mbox{for all $n$ large enough}.
\ees
Hence, a careful study is required for all $x(t) = \gamma$ in the {\em boundary sets} defined by
\bequ \label{setsBB}
\BB_{1,2} \equiv \{\gamma \in \RR_6 : d_{1,2}(\gamma) = 0\} \qandq \BB_{2,1} \equiv \{\gamma \in \RR_6 : d_{2,1}(\gamma) = 0\}
\eeq
FQR-ART aims to ``pull'' the fluid model to one of these two boundary sets during overloads, when sharing is actively taking place,
i.e., $\BB_{i,j}$ is the region of the state space where we aim the fluid model to be when pool $j$ helps class $i$, $i,j = 1,2$.

Unfortunately, there is no straightforward fluid counterpart to the stochastic processes $D^n_{1,2}$ and $D^n_{2,1}$
when the fluid is in the boundary sets. However, there are two related stochastic processes, operating in an infinitely faster time scale,
whose behavior determines the evolution of the fluid model, as we now explain.

\subsection{A Stochastic Averaging Principle} \label{secAP}

For the discussion now, assume that $x(t) \in \BB_{1,2}$ and consider $D^n_{1,2}$.
To be able to apply the results in \cite{PW13a}, we assume (for now) that the arrival rates are fixed (the arrival processes are
homogeneous Poisson processes) and that
$Z^n_{2,1} < \tau_{2,1}$, so that routing is determined solely on the value of $D^n_{1,2}$.
In particular, sharing can take place if $D^n_{1,2}(t) > 0$. Then, by Theorem 4.5 in \cite{PW13a},
\bequ \label{APlocal}
D^n_{1,2}(t) \Ra D_{1,2}(x(t), \infty) \qinq \RR \qasq n\tinf,
\eeq
where $D_{1,2}(\gamma, \cdot) \equiv \{D_{1,2}(\gamma, s) : s \ge 0\}$ is a CTMC associated with $\gamma \in \RR_6$
whose distribution is determined by the value $\gamma$. (There is a different process for each $\gamma$.)

An analogous result holds for $D^n_{2,1}$ when $x(t) \in \BB_{2,1}$.
The notation $D_{i,j}(\gamma, \infty)$ stands for
a random variable that has the steady-state distribution of the CTMC $D_{i,j}(\gamma, \cdot)$.
Loosely speaking, $D^n_{i,j}$ moves so fast when $x(t)$ is in $\BB_{i,j}$, that it reaches its steady state instantaneously as $n \tinf$.
Hence, we call $D_{i,j}(\gamma, \cdot)$ the {\em fast-time-scale process} (FTSP) associated with the point $\gamma$, or simply the FTSP.

Since we are interested in analyzing the indicator functions in \eqref{RepXstoc}, we first define for all $\gamma \in \RR_6$
$$D_{i,j}(\gamma, \cdot) \equiv + \infty \qifq d_{i,j}(\gamma) > 0 \qandq D_{i,j}(\gamma, \cdot) \equiv - \infty \qifq d_{i,j}(\gamma) < 0.$$
Next, we define
\bequ \label{pi's}
\bsplit
\pi_{1,2}(\gamma) & \equiv P(D_{1,2}(\gamma, \infty) > 0), \qforq \gamma \in \BB_{1,2}  \qandq \\
\pi_{2,1}(\gamma) & \equiv P(D_{2,1}(\gamma, \infty) > 0), \qforq \gamma \in \BB_{2,1}.
\end{split}
\eeq
Now, by Theorem 4.1 in \cite{PW13a}, which was proved for the process $D^n_{1,2}$ when $x \in \BB_{1,2}$, and assuming that
$Z^n_{2,1}(s) \le \tau^n_{2,1}$ over $[t_1, t_2]$ for all $n$ large enough, we have that, as $n \tinf$,
\bes
\int_{t_1}^{t_2} \1_{\mA^n_{1,2}(s)} ds \equiv \int_{t_1}^{t_2} \1_{\{\{D^n_{1,2}(s) > 0\} \cap \{Z^n_{2,1}(s) \le \tau^n_{2,1}\}\}} ds \Ra \int_{t_1}^{t_2} \pi_{1,2}(x(s)) ds.
\ees
Similarly, if $x \in \BB_{2,1}$ over an interval $[t_3, t_4]$, and $Z^n_{1,2}(s) \le \tau^n_{1,2}$ for all $n$ large enough over that interval, we have
\bes
\int_{t_3}^{t_4} \1_{\mA^n_{2,1}(s)} ds \equiv \int_{t_3}^{t_4} \1_{\{\{D^n_{2,1}(s) > 0\} \cap \{Z^n_{1,2}(s) \le \tau^n_{1,2}\}\}} ds \Ra \int_{t_3}^{t_4} \pi_{2,1}(x(s)) ds.
\ees
The convergence in both equations above holds uniformly.

We called these limits a ``stochastic averaging principle'', or simply an {\em averaging principle} (AP), since the process $D^n_{i,j}(t)$
is replaced by the {\em long-run average} behavior of the corresponding FTSP $D_{i,j}(x(t), \cdot)$ for each time $t$
over the appropriate interval.

In the FQR-ART settings, the AP holds under the assumption that $Z^n_{i,j}$ lies below the appropriate release threshold over the interval $[t_1, t_2]$
for all $n$ large enough (i.e., with probability converging to $1$ as $n \tinf$).
If $Z^n_{i,j}$ is larger than the appropriate release threshold for all $n$ large enough (again, with probability converging to $1$) over $[t_1, t_2]$,
then the limit of the integral considered above is clearly the $0$ function.
It remains to rigorously prove convergence theorems at points at which $Z^n_{i,j}(t) = \tau^n_{i,j} + o_P(n)$,
where $o_P(n)$ denotes a random variable satisfying $o_P(n)/n \Ra 0$ as $n \tinf$.
However, it is not hard to determine what the dynamics of the limit should be at such points if the limit exists.
That is the basis for our heuristic fluid model approximation below.

\subsection{Representation via an ODE} \label{secFluidLimits}

The heuristic limiting arguments above lead to the following fluid approximation for the X system under FQR-ART during overload periods.
Considering an interval $[0, T]$ for which
\bequ \label{Full}
z_{1,1}(t) + z_{2,1}(t) = m_1(t) \qandq z_{2,2}(t) + z_{2,2}(t) = m_2(t) \qforallq t \in [0, T],
\eeq
together with an initial condition $x(0)$, the fluid model of $X^n$ is the solution $x \equiv \{x(t) : t \ge 0\}$ over $[0, T]$ to the ODE:

\bequ \label{ODE-6}
\bsplit
\dot{q}_1(t) & = \lm_1(t) - \theta_1 q_1(t) - \Pi_{1,2}(x(t)) \left(\mu_{1,1} z_{1,1}(t) + \mu_{1,2} z_{1,2}(t) + \mu_{2,1} z_{2,1}(t) + \mu_{2,2} z_{2,2}(t)\right) \\
& \quad - (1 - \Pi_{1,2}(x(t)) - \Pi_{2,1}(x(t))) \left(\mu_{1,1} z_{1,1}(t) + \mu_{2,1} z_{2,1}(t) \right), \\
\dot{q}_2(t) & = \lm_2(t) - \theta_2 q_2(t) - \Pi_{2,1}(x(t)) \left(\mu_{1,1} z_{1,1}(t) + \mu_{1,2} z_{1,2}(t) + \mu_{2,1} z_{2,1}(t) + \mu_{2,2} z_{2,2}(t)\right) \\
& \quad - (1 - \Pi_{1,2}(x(t)) - \Pi_{2,1}(x(t))) \left(\mu_{2,2} z_{2,2}(t) + \mu_{1,2} z_{1,2}(t) \right), \\
\dot{z}_{1,2}(t) & = \Pi_{1,2}(x(t)) \mu_{2,2}z_{2,2}(t) - (1 - \Pi_{1,2}(x(t)))\mu_{1,2}z_{1,2}(t), \\
\dot{z}_{2,1}(t) & = \Pi_{2,1}(x(t)) \mu_{1,1}z_{1,1}(t) - (1 - \Pi_{2,1}(x(t)))\mu_{2,1}z_{2,1}(t), \\
\dot{m}_1(t) & = \dot{z}_{1,1}(t) + \dot{z}_{2,1}(t), \\
\dot{m}_2(t) & = \dot{z}_{2,2}(t) + \dot{z}_{1,2}(t),
\end{split}
\eeq
where, for $\pi_{i,j}(x(t))$ in \eqref{pi's}, $i,j = 1,2$,
\bes
\Pi_{i,j}(x(t)) := \left\{\begin{array}{ll}
\pi_{i,j}(x(t)) & \mbox{if } z_{j,i}(t) < \tau_{j,i}, \\
0               & \mbox{otherwise.}
\end{array}\right.
\ees

We remark that the ODE \eqref{ODE-6} can be equivalently represented by an integral equation resembling \eqref{RepXstoc},
but with the negligible martingale terms omitted,
all the stochastic processes replaced by their fluid counterparts, and the indicator functions replaced by the appropriate $\Pi_{i,j}$ functions.

In practice we do not a-priori know the value of $T$, and there is a need to make sure that the ODE is a valid approximation for the stochastic system.
We consider the ODE \eqref{ODE-6} valid (i.e., a legitimate representation of the evolution of the system) as long as the following two conditions are satisfied:
(i) the two queues are strictly positive; (ii) if a queue is equal to $0$
at some time $t \ge 0$, then the derivative of that queue is nonnegative at time $t$ (so that the queue is nondecreasing
at this time).
When the ODE \eqref{ODE-6} is not valid, then other fluid models should be employed to approximate the system.
We discuss such scenarios in \S \ref{secFluidNoShare} below.

We elaborate on Condition (ii). Consider, for example, the ODE for $q_1$ and assume that $q_1(t) = 0$ and $\dot{q}_1(t) < 0$ for some $t \ge 0$.
Necessarily $\Pi_{1,2}(x(t)) = 0$, because $d_{1,2}(x(t)) \le 0$, and the assumption that $\dot{q}_1(t) < 0$ implies that
\bequ \label{q1Valid}
\lm_1(t) - (1-\Pi_{2,1}(x(t))) (\mu_{1,1}z_{1,1}(t) + \mu_{2,1} z_{2,1}(t)) < 0.
\eeq
In addition, since all the class-$1$ arrivals must immediately enter service (for otherwise, the queue will be increasing), it also holds that
\bes
\dot{z}_{1,1}(t) = \lm_1(t) - \mu_{1,1}z_{1,1}(t).
\ees
Hence,
\bequ \label{m1Valid}
\bsplit
\dot{z}_{1,1}(t) + \dot{z}_{2,1}(t) & = \lm_1(t) - \mu_{1,1}z_{1,1}(t) \\
& \quad +  \Pi_{2,1}(x(t)) \mu_{1,1}z_{1,1}(t) - (1 - \Pi_{2,1}(x(t)))\mu_{2,1}z_{2,1}(t) \\
\quad & = \lm_1(t) - (1 - \Pi_{2,1}(x(t))) (\mu_{1,1}z_{1,1}(t) + \mu_{2,1} z_{2,1}(t)) \\
& \quad < 0,
\end{split}
\eeq
where the inequality follows from \eqref{q1Valid}.

Now, since $\dot{m}_1(t) = \dot{z}_{1,1}(t) + \dot{z}_{2,1}(t)$,
we see that pool $1$ can remain full just after time $t$ only if $m_1(t)$ happens to decrease exactly as in \eqref{m1Valid}.
However, $q_1$ is becoming negative, so that the ODE is not valid.
On the other hand, if \eqref{m1Valid} holds (which ODE \eqref{ODE-6} enforces to be equal to $\dot{m}_1(t)$) and $q_1(t) = 0$,
then necessarily $\dot{q}_1(t) < 0$, so that the queue is becoming negative.
In either case, we see that the ODE is valid as an approximation for the stochastic system when $q_1(t) = 0$ only if pool $1$
can be kept full without enforcing $q_1$ to become negative.
Similar reasonings hold for the $q_2$ and $m_2$ processes.

\subsection{The Fluid Model When There is No Active Sharing} \label{secFluidNoShare}

The ODE for the fluid model above was developed for all cases for which both pools are full, i.e., \eqref{Full} holds.
This is the main case because systems are typically designed to operate with very little extra service capacity (if any),
and is the primary case when overloads occur.
Nevertheless, the system may go through periods in which at least one of the pools is underloaded.  Hence, we now briefly describe the fluid models for
underloaded pools.

Consider an interval $I \subset [0, \infty)$.
If no sharing takes place and $z_{1,2}(t) = z_{2,1}(t) = 0$ for all $t \in I$,
then the two classes operate as two independent single-pool models (with time-varying parameters and staffing) over that interval $I$,
to which fluid limits are easy to establish. Specifically, assuming without loss of generality, that $I = [0, s)$
for some $0 < s < \infty$, the fluid dynamics of both classes obey the ODE
\bequ \label{ODE-NoShare}
\bsplit
\dot{q}_i(t) & = (\lm_i(t) - \mu_{i,i}z_{i,i}(t) - \theta_i q_i(t))\1_{\{q_i(t) \ge 0\}} \\
\dot{z}_{i,i}(t) & = \left\{\begin{array}{lr}
\dot{m}_i (t) & \mbox{if } q_i(t) > 0, \\
\lm_i(t) \1_{\{z_{i,i}(t) \le m_i(t)\}} - \mu_{i,i} z_{i,i}(t) & \mbox{if } q_i(t) = 0.
\end{array} \right.
\end{split}
\eeq
In the {\em time-invariant case}, when the arrival rates and staffing functions are fixed constants,
the unique solution for a given initial condition to the ODE in \eqref{ODE-NoShare} is easily seen to be
\bequ \label{fluidNoShare}
\bsplit
q_i(t) & = \left( \frac{\lm_i - \mu_{i,i}m_i}{\theta_i} + \left(q_i(0) - \frac{\lm_i - \mu_{i,i}m_i}{\theta_i}\right) e^{-\theta_i t} \right) \vee 0, \\
z_{i,i}(t) & = \left\{\begin{array}{lr}
m_{i,i} & \mbox{if } q_i(t) > 0, \\
\frac{\lm_i}{\mu_{i,i}} + \left( z_{i,i}(0) - \frac{\lm_i}{\mu_{i,i}}\right) e^{-\mu_{i,i}t} & \mbox{if } q_i(t) = 0.
\end{array} \right.
\end{split}
\eeq
where $a \vee b \equiv \max\{a,b\}$ and $(q_1(0), q_2(0), z_{1,1}(0), z_{2,2}(0))$ is a deterministic vector in $[0,\infty)^2 \times [0,m_1] \times [0, m_2]$.

If $z_{1,2}(s_0) > 0$ (or $z_{2,1}(s_0) > 0$) for some $s_0 \ge 0$ and there is no active sharing over the interval $[s_0, s_1)$, then $z_{1,2}$ ($z_{2,1}$)
is strictly decreasing over that interval. Then $z_{i,j}$, $i \ne j$, satisfies the ODE
\bes
\dot{z}_{i,j}(t) = - \mu_{i,j}z_{i,j}(t), \quad s_0 \le t < s_1
\ees
which is the same as the ODE for $z_{i,j}$ in \eqref{ODE-6} with $\Pi_{i,j} = 0$.

\begin{remark} \label{remODEsol}
{\em
A proof of existence of a unique solution to the ODE \eqref{ODE-6} following the lines of \cite{PW11b} requires showing that the RHS is a local Lipschitz continuous function of $x$
and is piecewise continuous in $t$. We do not prove such a result here, but it is important to consider arrival rates and staffing functions
that ensure that the right side of the ODE satisfies the piecewise continuity condition in the time argument.
}
\end{remark}

\section{Solving the ODE} \label{secPis}

To appreciate that the algorithm cannot be a routine solution of an ODE, observe that
computing the solution to \eqref{ODE-6} requires computing the two steady-state probabilities $\pi_{1,2}(x(t))$ and $\pi_{2,1}(x(t))$
for all times $t$ and states $x(t) \in \RR_6$. Simplification is achieved when $r_{1,2} = r_{2,1} = 1$, because the FTSP's $D_{i,j}(x(t), \cdot)$, $i,j = 1,2$,
become simple {\em birth-and-death} (BD) processes.
To facilitate the discussion we thus consider this simpler case and refer to \S 6.2 in \cite{PW11b}
for the treatment of the FTSP $D_{1,2}$ as a {\em quasi-birth-and-death process} (QBD) when the ratio parameters are not equal to $1$.
(In \cite{PW11b} FQR-T is studied with one overload incident, with pool $1$ receiving help,
but the same method can be applied to $D_{2,1}$ with sharing in the opposite direction.)

For simplicity, we again start by assuming that the arrival processes are homogenous Poisson processes,
having constant arrival rates $\lm_1$ and $\lm_2$ over $[0, T]$, and that the staffing functions are also fixed over that time interval at $m_1$ and $m_2$.
Recall that $D_{i,j}(\gamma, \cdot) \equiv \infty$ if $d_{i,j}(\gamma) > 0$ and $D_{i,j}(\gamma, \cdot) \equiv - \infty$ if $d_{i,j}(\gamma) < 0$,
and let $\AA_{1,2}$ and $\AA_{2,1}$ be the subsets of $\RR_6$ in which the FTSP's $D_{1,2}(\gamma, \cdot)$ and $D_{2,1}(\gamma, \cdot)$
are positive recurrent, i.e.,
\bequ \label{Asets}
\AA_{1,2} \equiv \{\gamma \in \BB_{1,2}: 0 < \pi_{1,2}(\gamma) < 1\} \qandq \AA_{2,1} \equiv \{\gamma \in \BB_{2,1} : 0 < \pi_{2,1}(\gamma) < 1\}.
\eeq
By definition, if the fluid model at time $t$ is in $\AA_{i,j}$, i.e., $x(t) \in \AA_{i,j}$, then $d_{i,j}(x(t)) = 0$.
However, if $d_{i,j}(x(t)) = 0$, then $x(t)$ is not necessarily in $\AA_{i,j}$, because the FTSP $D_{i,j}(x(t), \cdot)$
may be transient (drift to $+\infty$ or $-\infty$) or null recurrent; in particular,
{\em The evolution of the fluid model is determined by the distributional characteristics of the FTSP's $D_{1,2}$ and $D_{2,1}$}.
Hence, even before we try to compute $\pi_{i,j}(x(t))$, which is necessary in order to solve the ODE \eqref{ODE-6}, there is a need to
determine whether $x(t)$ is in one of the sets $\AA_{1,2}$ or $\AA_{2,1}$.
We focus on $D_{1,2}$, with the analysis of $D_{2,1}$ being similar.

To determine the behavior of the FTSP $D_{1,2}$ it is again helpful to think of $x$ as a fluid limit of the fluid-scaled sequence $\{\barx^n : n \ge 1\}$
and to recall that $D_{1,2}$ was achieved as a limit of $D^n_{1,2}$ without any scaling; see \eqref{APlocal}.
(See also Theorem 4.4 in \cite{PW13a} which provides a process-level limit relating $D_{1,2}$ and $D^n_{1,2}$.)
Hence, both processes are defined on the same state space,
which, for $r_{1,2} = 1$, is $\ZZ \equiv \{\dots, -1,0,1, \dots\}$.

Now, for a fixed $x(t)$, when $D_{1,2}(x(t), \cdot) = m > 0$, the birth and death rates of the FTSP are, respectively,
\bes
\bsplit
\lm^+(x(t),m) & \equiv \lm_1 + \theta_2 q_2(t), \\
\mu^+(x(t),m) & \equiv \lm_2 + \mu_{1,1}z_{1,1}(t) + \mu_{1,2}z_{1,2}(t) + \mu_{2,1}z_{2,1}(t) + \mu_{2,2}z_{2,2}(t) + \theta_1 q_1(t).
\end{split}
\ees
In analogy to the (non-Markov) process $D^n_{1,2} = Q^n_1 - Q^n_2 - k^n_{1,2}$,
$\lm_+(x(t),m)$ corresponds to an increase of $D_{1,2}$ due to arrival to queue $1$
plus an abandonment from queue $2$ (since either one of these two events cause an increase by $1$ of $D^n_{1,2}$ in the stochastic system).
Since any other event causes $D^n_{1,2}$ to decrease by $1$, due to the scheduling rules of FQR-ART, we get the expression for $\mu^+(x(t),m)$.

Next, if $D_{1,2}(x(t), m) = m \le 0$, the birth and death rates are, respectively,
\bes
\bsplit
\lm^-(x(t), m) & \equiv \lm_1 + \mu_{2,2}z_{2,2}(t) + \mu_{1,2}z_{1,2}(t) + \theta_2 q_2(t), \\
\mu^-(x(t), m) & \equiv \lm_2 + \mu_{1,1}z_{1,1}(t) + \mu_{2,1}z_{2,1}(t) + \theta_1 q_1(t).
\end{split}
\ees
Again, whenever $D^n_{1,2}$ is non-positive and sharing is taking place with pool $2$ helping class $1$,
a ``birth'' occurs if there is an arrival to queue $1$ or an abandonment from queue $2$, or if there is a service completion in pool $2$
(since then a newly available type-$2$ agent takes his next customer from queue $2$).
Similarly, a ``death'' occurs if there is an arrival to class $2$, an abandonment from queue $1$, or a service completion in pool $1$.

We see that the FTSP $D_{1,2}(x(t), \cdot)$ is a two-sided $M/M/1$ queue, i.e., it behaves like an $M/M/1$ queue
with ``arrival rate'' $\lm^+(x(t), m)$ and ``service rate'' $\mu^+(x(t), m)$ for all $m > 0$, and behaves like a different $M/M/1$
queue with ``arrival rate'' $\mu^-(x(t), m)$ and ``service rate'' $\lm^-(x(t), m)$, for all $m \le 0$. Thus, for
\bes
\delta^+(\gamma) \equiv \lm^+(\gamma, \cdot) - \mu^+(\gamma, \cdot) \qandq \delta^-(\gamma) \equiv \lm^-(\gamma, \cdot) - \mu^-(\gamma, \cdot), \quad \gamma \in \BB_{1,2},
\ees
the set $\AA_{1,2}$ can be characterized via
\bes
\AA_{1,2} \equiv \{\gamma \in \BB_{1,2}: \delta^+(\gamma) < 0 < \delta^-(\gamma)\}.
\ees
Next, letting $T^+(\gamma)$ and $T^-(\gamma)$ denote, respectively, the busy period of the $M/M/1$ in the positive region
and the busy period of the $M/M/1$ in the negative region,
and using simple alternating renewal arguments for the renewal process $D_{1,2}(\gamma, \cdot)$, we have
\bequ \label{piCalc}
\pi_{1,2}(\gamma) = \frac{E[T^+(\gamma)]}{E[T^+(\gamma)] + E[T^-(\gamma)]},
\eeq
where, from basic $M/M/1$ theory,
\bes
E[T^{\pm}(\gamma)] = \frac{1}{\mu^{\pm}(\gamma) - \lm^\pm(\gamma)}.
\ees
Note that if $d_{1,2}(\gamma) = 0$ but $\gamma \notin \AA_{1,2}$, then $\pi_{1,2}(\gamma)$ is equal to either $1$ or $0$.
In particular,
\bequ \label{piNotAA}
\mbox{if~ $\delta^+(\gamma) \ge 0$,~ then ~$\pi_{1,2}(\gamma) = 1$~ and if~ $\delta^-(\gamma) \le 0$ ~then ~$\pi_{1,2}(\gamma) = 0$.}
\eeq
There are no other options, since for any $\gamma = x(t)$ for which both pools are full (as is required for the
ODE \eqref{ODE-6} to be valid), it holds that
\bes
\delta^-(x(t)) - \delta^+(x(t)) = 2(\mu_{1,2}z_{1,2}(t) + \mu_{2,2}z_{2,2}(t)) > 0, \quad
\ees
where the inequality above follows from the fact that $z_{1,2}(t) + z_{2,2}(t) = m_2(t) > 0$.

We see that the sets $\AA_{i,j}$ and the computation of $\pi_{i,j}(\cdot)$ are completely determined by the staffing, arrival rates, service and abandonment rates
for any given point $\gamma \in \RR_6$, where the only points that require careful analysis are those in one of the two sets $\BB_{i,j}$.
However, recall that we have assumed for simplicity that the arrival rates and staffing functions are not time dependent.
If, instead, the arrival rates or the staffing functions are time dependent, then the distribution of the FTSP $D_{i,j}(x(t), \cdot)$
is also time dependent. In particular, given a $\gamma \in \RR_6$ we cannot determine whether
$D_{1,2}(\gamma, \cdot)$ is positive recurrent or not, since that may depend on the time $t \in [0, T]$.
Thus, the sets at which the FTSP's are ergodic are themselves time dependent.
Hence, for a full analysis, we would need
 to consider sets of the form $\{\AA_{i,j}(t) : t \in [0, T]\}$, where
\bequ \label{PosRec}
\AA_{i,j}(t) \equiv \{(\gamma, t) \in \BB_{i,j} \times \RR_+ : \delta^+(\gamma, t) < 0 < \delta^-(\gamma, t)\},
\eeq
where $\delta^+(\gamma, t)$ and $\delta^-(\gamma, t)$ are the drifts of the FTSP $D_{1,2}(\gamma, \cdot)$ at the point $\gamma$ at time $t$.
Fortunately, for the purpose of solving the ODE, we do not actually need to characterize the sets $\{\AA_{i,j}(t) : t \in [0, T]\}$,
because we can determine whether $D_{i,j}(x(t), \cdot)$ is ergodic at each time $t$ as we solve the ODE.

\subsection{A Numerical Algorithm to Solve the ODE} \label{secAlg}

Given the ODE in \eqref{ODE-6} with a fully specified RHS at each $t$,
we compute the solution $x$ over an interval $[0, T]$ by employing the classical Euler method, combined with the AP.
Given a step size $h$ and the time $T$, the number of iterations needed is $N \equiv T/h$.
Let $\dot{x} = \Psi(x)$, where $\Psi(x)$ is the RHS of the appropriate ODE, e.g., if both pools are full, then $\Psi(x)$ is the RHS of \eqref{ODE-6}.
Given $x(0)$, we can compute $x(h)$ using the first Euler step: $x(h) = x(0) + h\Psi(x(0))$.
Given $x(h)$ we can compute $\Pi_{1,2}(x(h))$ and $\Pi_{2,1}(x(h))$, if needed, and then compute $x(2h)$ using the second Euler step.
In general, the solution to the ODE is computed via
\bes
x((k+1)h) = x(kh) + h \Psi(x(kh)), \quad 0 \le k \le N,
\ees
where at each step, if $x(kh) \in \BB_{1,2}$ or $x(kh) \in \BB_{2,1}$, we can compute $\Pi_{1,2}(kh)$ and $\Pi_{2,1}(kh)$ as explained above.

The algorithm just described remains unchanged when the ratio parameters are general (not equal to $1$), except that the sets $\AA_{i,j}$
and the computations of $\pi_{i,j}$ are more complicated (the FTSP's are no longer BD processes).
We refer to \cite{PW11b} for these more complicated settings.

To evaluate the RHS in each step, we
use the analysis in \S \ref{secPis}, starting at a given initial condition $x(0)$,
since we can now determine the value of $\Pi_{i,j}(x(t))$ for each $t \ge 0$. For example, if at a time $t \ge 0$ $d_{1,2}(x(t)) = 0$,
then we check whether \eqref{PosRec} holds, so that $x(t) \in \AA_{1,2}(t)$. If $z_{2,1}(t) \le \tau_{2,1}$, then $\Pi_{1,2}(x(t)) = \pi_{1,2}(x(t))$
and it can be computed using \eqref{piCalc}. If $z_{2,1}(t) > \tau_{2,1}$, then $\Pi_{2,1}(t) = 0$.
If $d_{1,2}(x(t)) = 0$ but $x(t) \notin \AA_{1,2}(t)$, i.e., if \eqref{PosRec} does not hold,
then we can determine the value of $\pi_{1,2}(x(t))$,
and thus of $\Pi_{1,2}(x(t))$, by computing the drifts of the FTSP and employing \eqref{piNotAA}
(replacing the drifts in \eqref{piNotAA} with the time dependent drifts as in \eqref{PosRec}).
Similarly we can compute the value of $\Pi_{2,1}(x(t))$ whenever $d_{2,1}(x(t)) = 0$.

In all other regions of the state space for which both pools are full, i.e., $z_{i,j}(t) + z_{j,i}(t) = m_j(t)$, $i \ne j$,
we can easily determine the value of $\pi_{1,2}(x(t))$ by considering whether $d_{i,j}(x(t))$ is bigger or smaller than $0$.
For example, if at time $t \ge 0$ $d_{1,2}(x(t)) > 0$, then $\pi_{1,2}(x(t)) = 1$ and if $d_{1,2}(x(t)) < 0$, then $\pi_{1,2}(x(t)) = 0$.
This, together with the value of $z_{2,1}(t)$, immediately gives the value of $\Pi_{1,2}(x(t))$.

We need to use other fluid equations when at least one of the two pools is not full.
If, for example $z_{1,1}(t) + z_{2,1}(t) < m_1(t)$, then necessarily $q_1(t) = 0 < k_{1,2}$, so that
\bes
\dot{z}_{1,2}(t) = -\mu_{1,2} z_{1,2}(t) \qandq \dot{z}_{1,1}(t) = \lm_1(t) - \mu_{1,1}z_{1,1}(t).
\ees
The evolution of $z_{2,1}$ in this case is determined by whether $q_{2}(t) < k_{2,1}$ or $q_{2}(t) \ge k_{2,1}$.
In the first case $z_{2,1}(t)$ must be strictly decreasing at time $t$ if it is positive, or remain at $0$ otherwise.
In the latter case, when $q_{2}(t) \ge k_{2,1}$, the excess fluid - that is not routed to pool $2$ and does not abandon, if such excess fluid exists -
is flowing to pool $1$. We thus have $\dot{z}_{2,1}(t)$ is equal to
\bequ \label{whenInk}
\begin{array}{ll}
-\mu_{2,1}z_{2,1}(t) & \mbox{if $q_2(t) < k_{2,1}$} \\
-\mu_{2,1}z_{2,1}(t) + (\lm_2(t) - \mu_{2,2}z_{2,2}(t) - \mu_{1,2}z_{1,2}(t) - \theta_2 k_{2,1})^+ & \mbox{if $q_2(t) = k_{2,1}$}
\end{array}
\eeq

Similar reasonings lead to the fluid model of $z_{1,2}$ when pool $1$ is full, but pool $2$ has spare capacity.

If both pools have spare capacity at time $t$, then $q_1(t) = q_2(t) = 0$ and
\bes
\dot{z}_{i,j}(t) = -\mu_{i,j}z_{i,j}(t) \qandq \dot{z}_{i,i}(t) = \lm_i - \mu_{i,i}z_{i,i}(t), \quad i,j = 1,2, \quad i \ne j.
\ees

\begin{remark} \label{RemAlgo}{\em
If at iteration $k \ge 0$ the solution lies outside the set $\BB_{1,2} \cup \BB_{2,1}$, then
due to the discreteness of the algorithm, there is a need to ensure that the boundary is not missed in the following iterations.
Hence, if in the $k^{th}$ iteration $d_{1,2}(x(kh)) > 0$ ($< 0$) and in the $(k+1)^{st}$ iteration $d_{1,2}((x(k+1)h)) < 0$ ($> 0$),
then the boundary $d_{1,2}$ necessarily was missed, because the fluid is continuous, and
so we set $d_{1,2}((x(k+1)h)) = 0$. We then check whether $x((k+1)h) \in \AA_{1,2}((k+1)h)$, compute $\pi_{1,2}(x(k+1)h)$
and use its value to compute the value in the $(k+2)^{nd}$ iteration.
It is significant that {\em we do not force the solution to be on the boundary},
e.g., we do not compute $q_1((k+1) h)$ and use its value to compute $q_2((k+1) h)$ via
\bequ \label{InBB}
q_2((k+1) h) = q_1((k+1) h) - k_{1,2}.
\eeq
We solve the six-dimensional ODE in \eqref{ODE-6}, and if indeed \eqref{InBB} holds whenever it should,
then we have a good indication that the algorithm works.
That is, we can check at which iteration the boundary $\BB_{1,2}$ was hit, and then observe if $q_1(t) - q_2(t) = k_{1,2}$ over an interval
for which we have indication that this should hold.
(Of course, the solution to the algorithm might leave the boundary for legitimate reasons, i.e., because the fluid model leaves it.)
}
\end{remark}

\section{Numerical Examples} \label{secNum}

We now study three examples. The first two are piecewise-continuous models, whereas the third is for a general time-varying model.
In all three examples the system starts empty, so that we also check the numerical algorithm in
periods when \eqref{Full} does not hold, as in \S \ref{secFluidNoShare}.

We compare the numerical solutions to the ODE to simulations, to see how well the fluid model approximates stochastic systems.
In the first two examples we simulate three systems, each can be considered as a component in a sequence $\{\barx^n : n \ge 1\}$. In the smallest system we take
$50$ agents in each service pool, in the middle one there are $100$ agents in a pool, and the largest has $400$ agents in each pool, i.e.,
we simulate $\barx^n$ for $n = 50, 100, 400$. That allows us to observe the ``convergence'' of the stochastic system to the fluid approximation.
We plot the fluid and simulation results together, normalized to $n=10$.
(E.g., for the system with $400$ agents in each pool we divide all processes by $40$.)

The following parameters are used for all three simulations:

$\mu_{1,1} = \mu_{2,2} = 1;$ $\mu_{1,2} = \mu_{2,1} = 0.8$, $\theta_1 = \theta_2 = 0.5$. In addition, we take $r_{1,2} = r_{2,1} = 1$.
We take $k^n_{1,2} = k^n_{2,1} = 0.3n$; $\tau^n_{1,2} = \tau^n_{2,1} = 0.02n$, so that, for $n = 50, 100, 400$, we have
$k^n_{1,2} = k^n_{2,1} = 15, 30, 120$ and $\tau^n_{1,2} = \tau^n_{2,1} = 1, 2, 8$, respectively.

\subsection{A Single Overload Incident} \label{secSimSingle}

The first example aims to check whether FQR-ART detects overloads automatically when they occur and starts sharing in the right direction,
and whether, once an overload incident is over, FQR-ART avoids oscillations.
In particular, over the time interval $[0, 60]$ the arrival rates are as follows: $\lm^n_2 = n$ throughout that time interval.
Over $[0, 20)$ and $[40, 60]$ the arrival rate to pool $1$ is $\lm^n_1 = n$. Hence, both pools are normally loaded during these two subintervals.
However, during the interval $[20, 40)$ the arrival rate of class $1$ changes to $\lm^n_1 = 1.4n$,
so that, during $[20, 40)$ the system is overloaded, and pool $2$ should be helping class $1$.

We compare the solution to the fluid equations, solved using the algorithm, to an average of $1000$ independent simulation runs for the three cases
$n = 50, 100, 400$.
The results are shown in Figures \ref{figQ1_1OL}-\ref{figZ_1OL} below.
In addition Figure \ref{figd12_1OL} plots $q_1 - r_{1,2} q_2 - k_{1,2}$. Since shortly after time
$20$ the value is $0$ in Figure \ref{figd12_1OL}, we have a strong indication that the numerical solution is correct, because during most of the overload period,
when sharing takes place, it should hold that $d_{1,2}(x(t)) = 0$.

\begin{figure}[h!]
  \hfill
  \begin{minipage}[t]{.4\textwidth}
    \begin{center}
      \includegraphics[width=6cm]{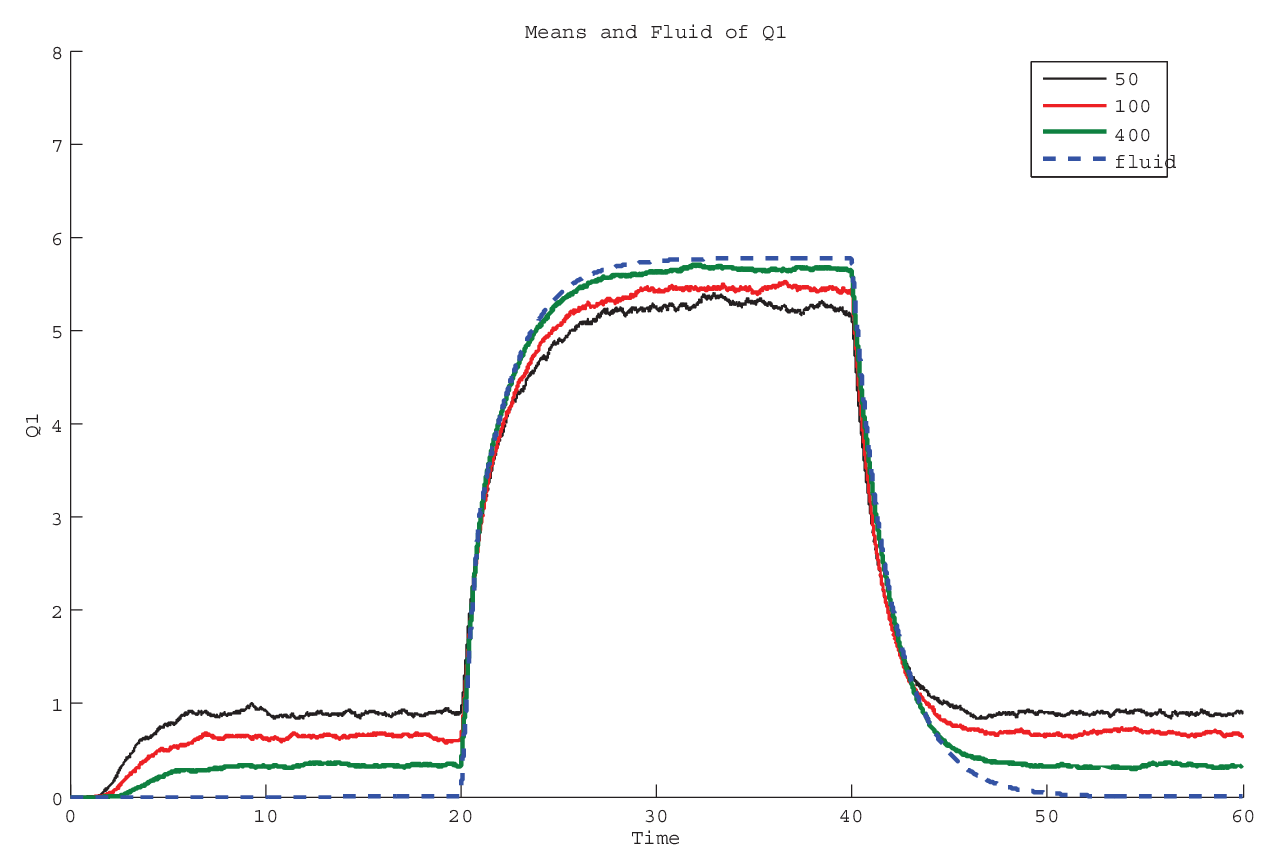}
      \caption{comparison of the fluid model to simulations of $10\barq^n_1$ for $n = 50$, $100$ and $400$ with a single overload}
      \label{figQ1_1OL}
    \end{center}
  \end{minipage}
  \hfill
  \begin{minipage}[t]{.4\textwidth}
    \begin{center}
      \includegraphics[width=6cm]{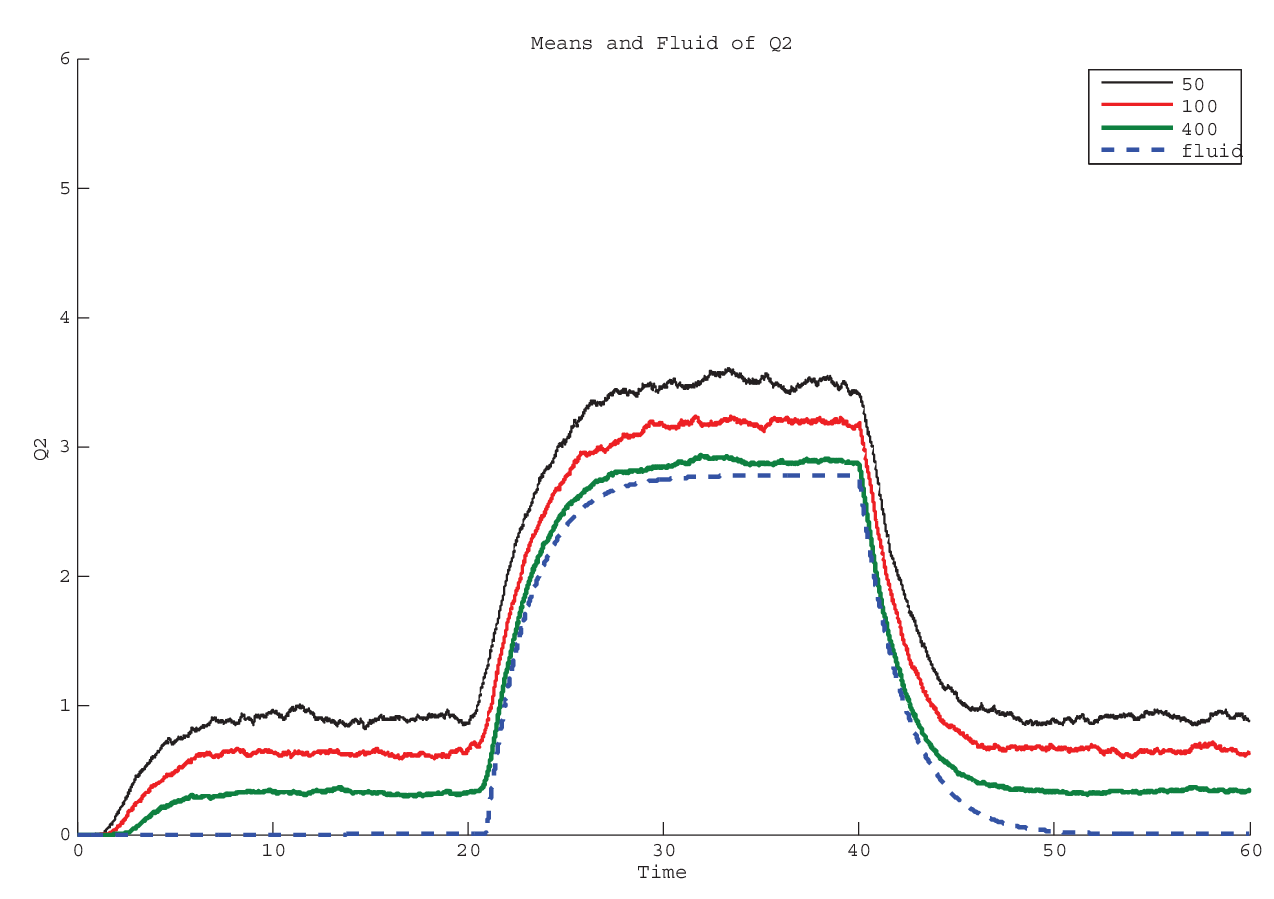}
      \caption{comparison of the fluid model to simulations of $10\barq^n_2$ for $n = 50$, $100$ and $400$ with a single overload}
      \label{figQ2_1OL}
    \end{center}
  \end{minipage}
  \hfill
  \begin{minipage}[t]{.4\textwidth}
    \begin{center}
      \includegraphics[width=6cm]{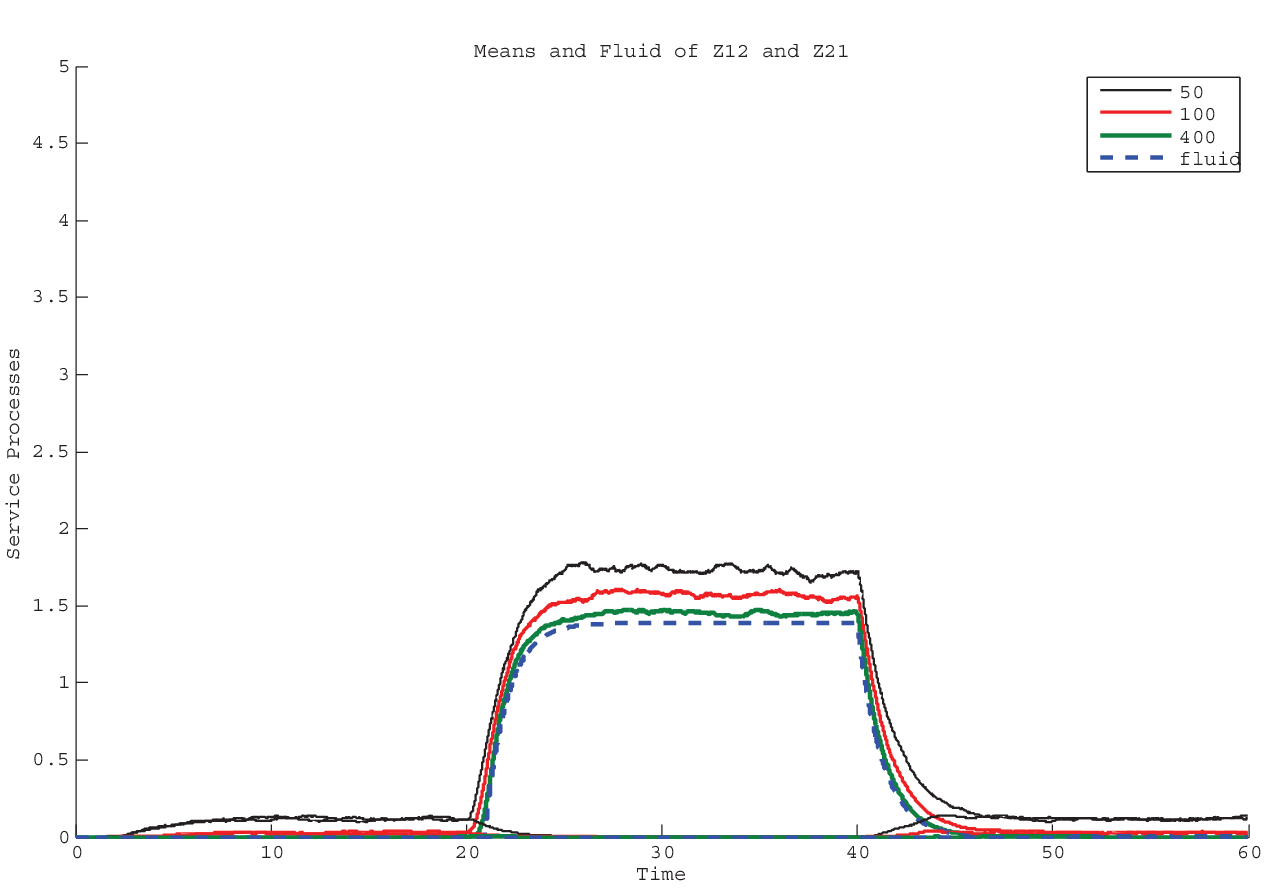}
      \caption{comparison of the fluid model to simulations of $10\barz^n_{1,2}$ for $n = 50$, $100$ and $400$ with a single overload}
      \label{figZ_1OL}
    \end{center}
  \end{minipage}
\hfill
  \begin{minipage}[t]{.4\textwidth}
    \begin{center}
      \includegraphics[width=6cm]{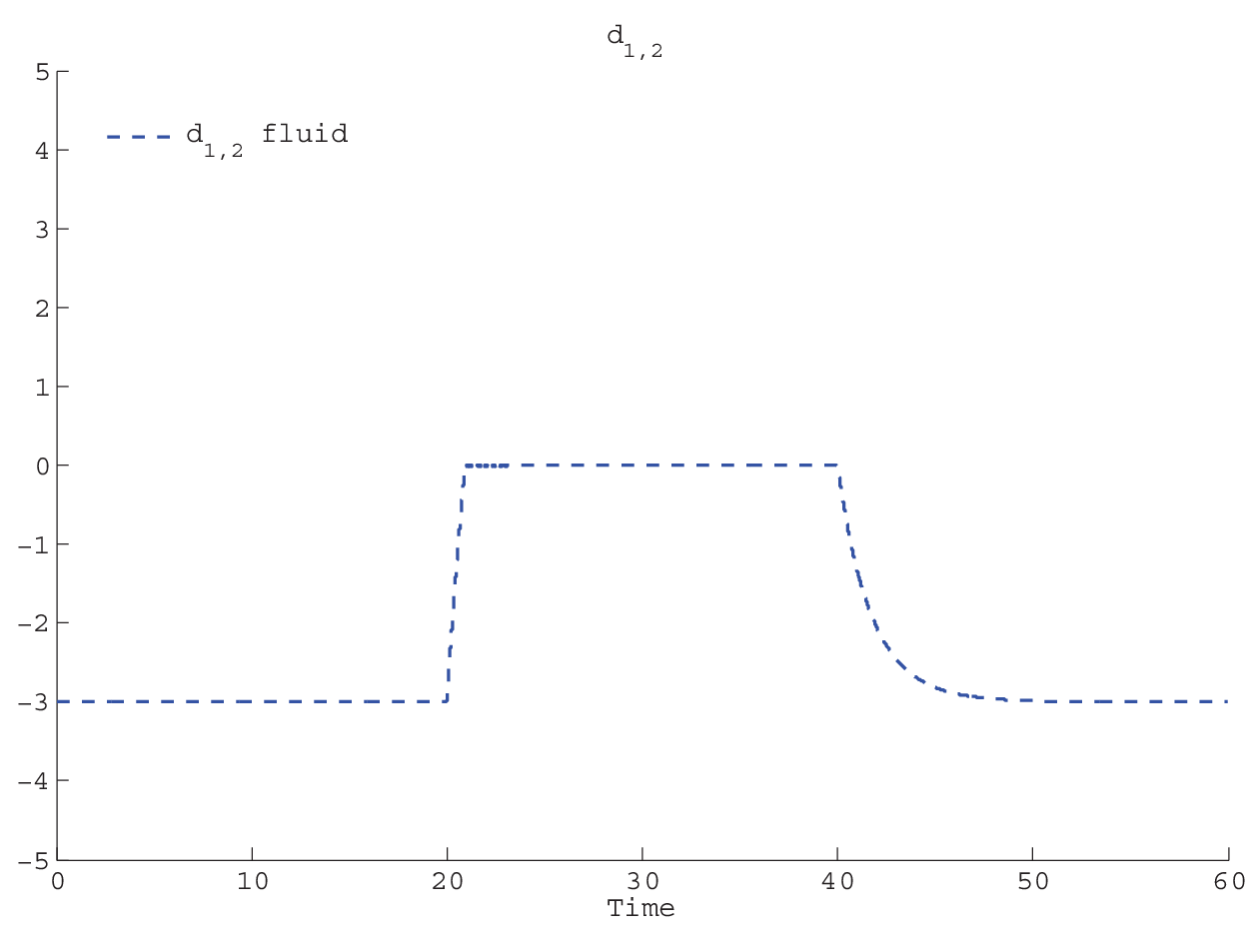}
      \caption{plot of $q_{1} - r_{1,2} q_2 - k_{1,2}$ in the single overload example}
      \label{figd12_1OL}
    \end{center}
  \end{minipage}
  \hfill
\end{figure}

The simulation experiments indicate that the fluid model approximates well the mean behavior of the system even for relatively small systems, e.g., when
$n = 50$. Of course, the accuracy of the approximation grows as $n$ becomes larger.
The simulation experiments show that FQR-ART quickly detects the overload and the correct direction of sharing.
Moreover, the control ensures that there are no oscillations, as in \S \ref{secOSC}.

Another observation is that when the system is normally loaded and there is no sharing, the fluid model, which has null queues,
does not describe the queues well. In those cases there is an increased importance to stochastic refinements for the queues.
If there is only negligible sharing, as FQR-ART ensures, then such stochastic refinements are well approximated by diffusion limits
for the Erlang A model, as in \cite{GMR02}.

\subsection{Switching Overloads} \label{secSimSwitch}

\begin{figure}[h!]
  \hfill
  \begin{minipage}[t]{.4\textwidth}
    \begin{center}
      \includegraphics[width=6cm]{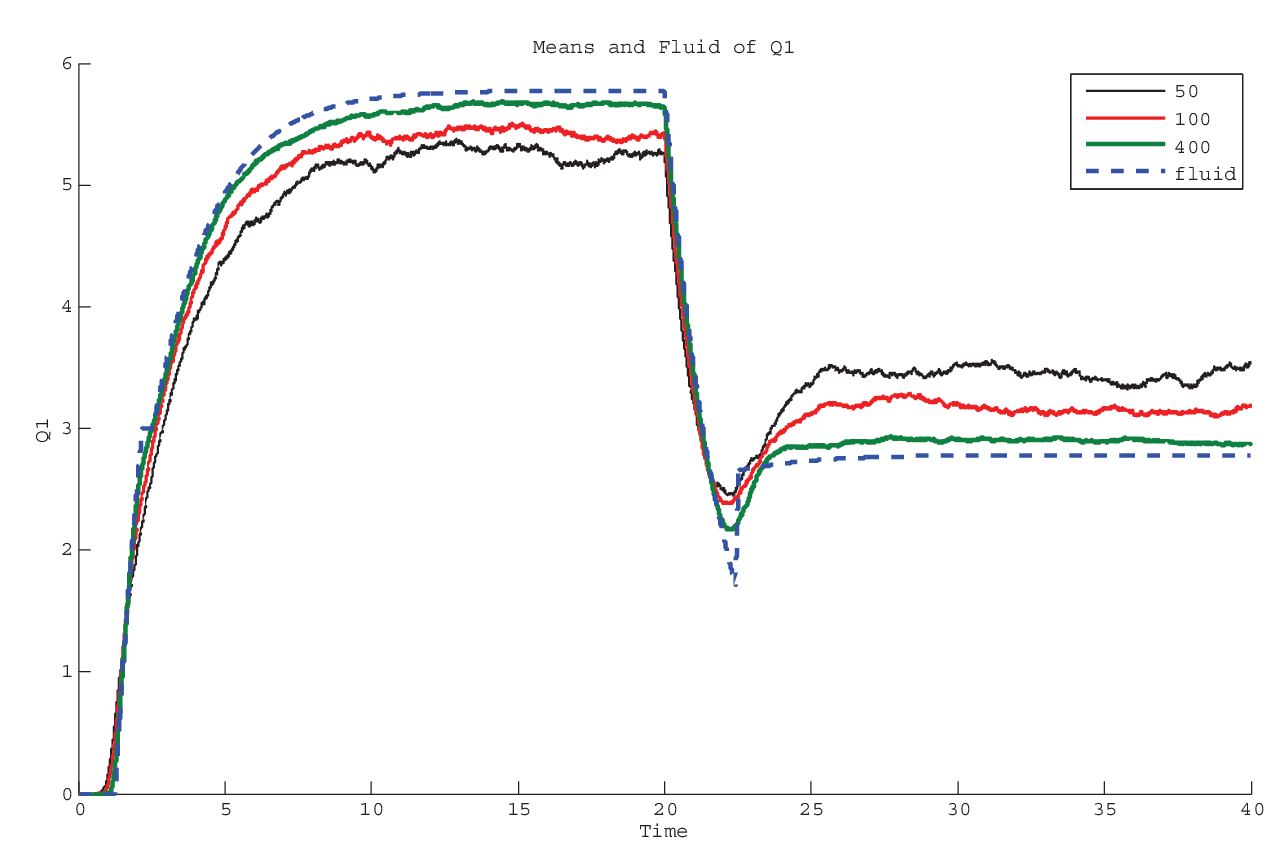}
      \caption{comparison of the fluid model to simulations of $10\barq^n_1$ for $n = 50$, $100$ and $400$ with the switching overloads}
      \label{figQ1switch}
    \end{center}
  \end{minipage}
  \hfill
  \begin{minipage}[t]{.4\textwidth}
    \begin{center}
      \includegraphics[width=6cm]{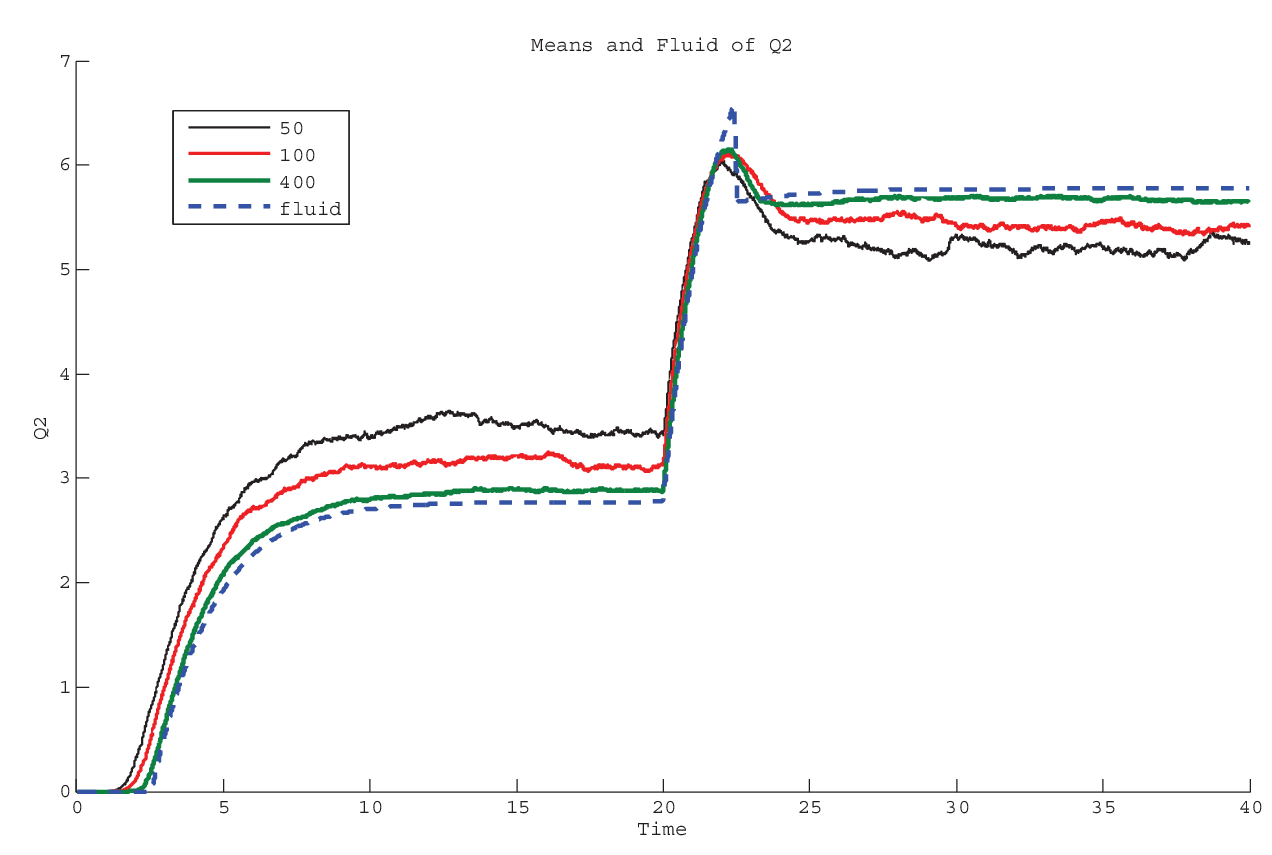}
      \caption{comparison of the fluid model to simulations of $10\barq^n_2$ for $n = 50$, $100$ and $400$ with the switching overloads}
      \label{figQ2switch}
    \end{center}
  \end{minipage}
  \hfill
\end{figure}

\begin{figure}
  \begin{minipage}[t]{.4\textwidth}
    \begin{center}
      \includegraphics[width=6cm]{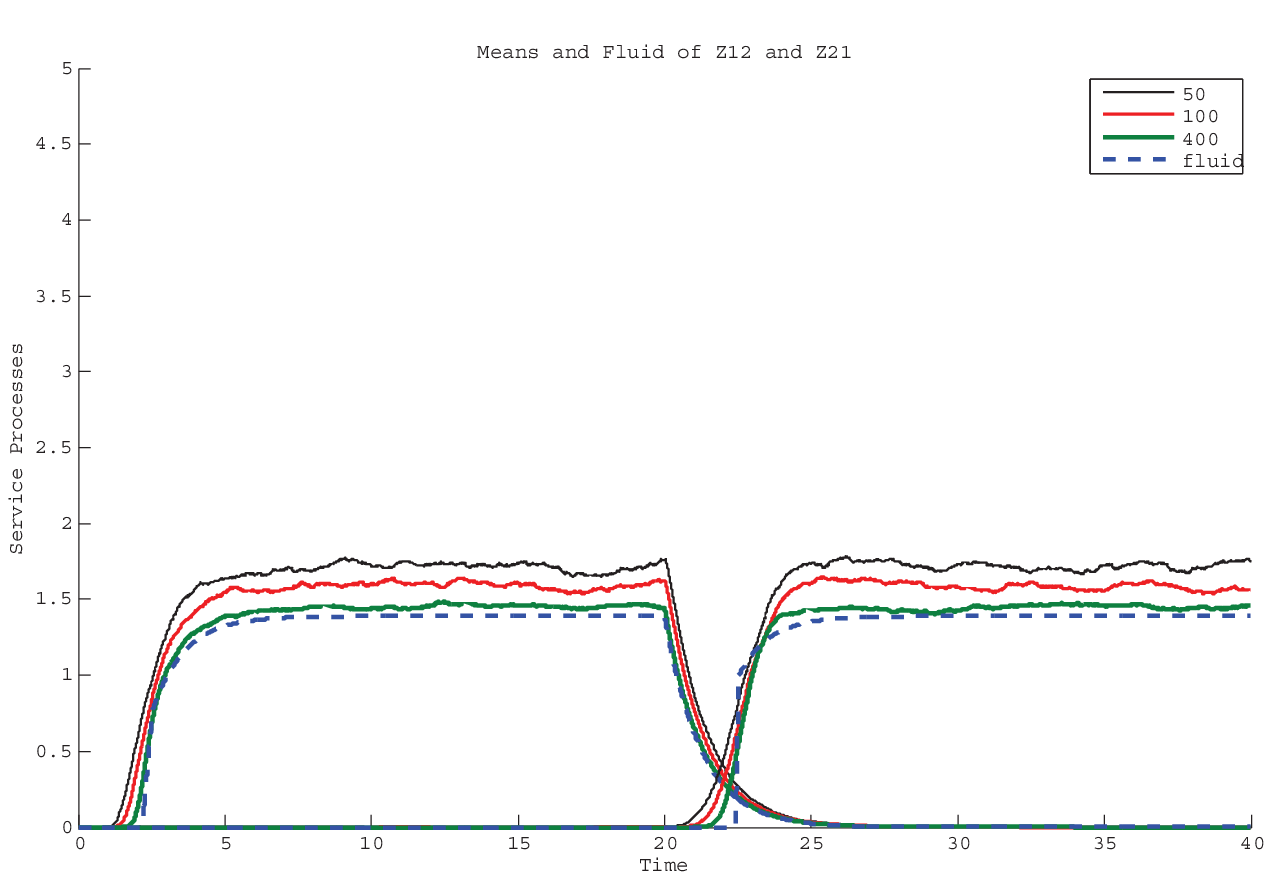}
      \caption{comparison of the fluid model to simulations of $10\barz^n_{1,2}$ for $n = 50$, $100$ and $400$ with the switching overloads}
      \label{figServSwitch}
    \end{center}
  \end{minipage}
  \hfill
  \begin{minipage}[t]{.4\textwidth}
    \begin{center}
      \includegraphics[width=6cm]{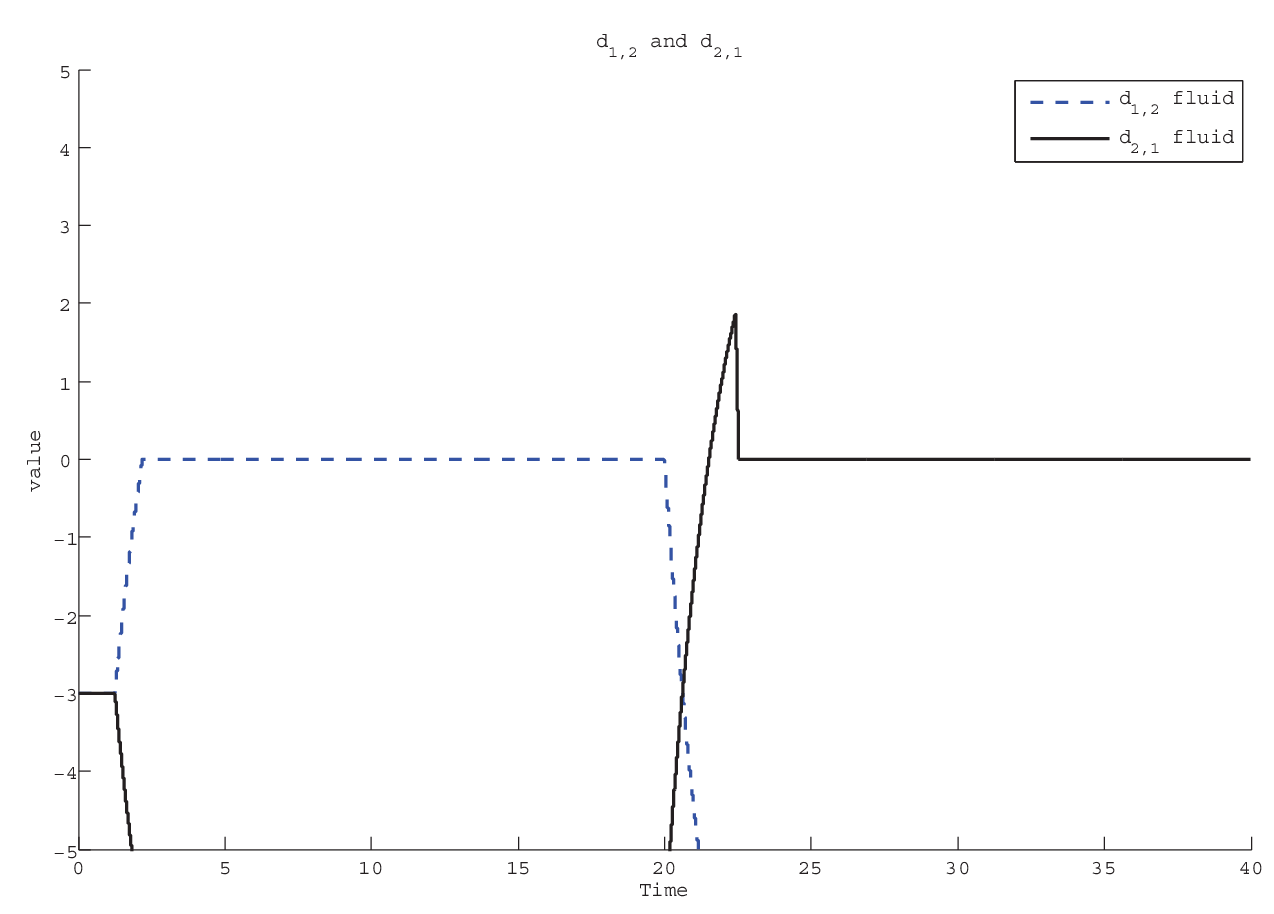}
      \caption{the two fluid difference processes with the switching overloads}
      \label{figd}
    \end{center}
  \end{minipage}
  \hfill
\end{figure}

In the second example we consider an overloaded system, with pool $1$ being overloaded initially,
and with the direction of overload switching after some time, making pool $2$ overloaded.
Specifically, we let the arrival rates be $\lm^n_1 = 1.4n$ and $\lm^n_2 = n$ over $[0, 20)$,
and $\lm^n_1 = n$, $\lm^n_2 = 1.4n$ on $[20, 40]$.
The results are plotted in Figures \ref{figQ1switch}-\ref{figServSwitch}.
Figure \ref{figd} plots $q_1 - r_{1,2}q_2 - k_{1,2}$ and $r_{2,1}q_2 - q_1 - k_{2,1}$.

Once again, the fact that the appropriate difference process equals to $0$
shortly after the corresponding overload begins is an indication that the solution to the ODE is correct, since each queue is calculated
via the averaging principle, without forcing the relations $d_{1,2}(x(t)) = 0$ and $d_{2,1}(x(t)) = 0$.


As in the figures in \S \ref{secSimSingle},
it is easily seen from the figures above that the fluid model approaches a fixed point, so long as the arrival rates are fixed.
Then, once a change in the rates occurs, the fluid goes through a new transient period until it relaxes in a new fixed point.

\subsection{General Non-stationary Model with Switching Overloads}\label{secSimNone}

We next test our algorithm in a more challenging time-varying
example. This example is unrealistic in call-center setting,
because the arrival rates and staffing functions are not likely
to change so drastically, but it demonstrates the robustness of
our fluid model and of the algorithm.

We assume that the arrival rate to pool $1$ over the time
period $[0, 20)$ is sinusoidal. We further assume that
management anticipated the basic sinusoidal pattern of the
arrival rate, but did not anticipated the magnitude, so that
pool $1$ is overloaded. To specify the staffing with the sinusoidal arrival rate,
we assume that staffing follows the appropriate {\em
infinite-server} approximation; see, e.g., Equation (9) in
\cite{FMMW08}.  The purpose of that staffing rule in our
setting, is to stabilize the system at a fixed point
eventually, as in the examples above. In particular, for $t \in
[0, 20]$, we let

\vspace{0.05in}
$\lm^n_1(t) = 1.3n + 0.1n \sin(t)$\: and\:
$m^n_1(t) = n + 0.05n[\sin(t) - \cos(t)]$;

$\lm^n_2(t) = n$\: and\: $m^n_2(t) = n$.

\vspace{0.05in} Then, on the time interval $[20, 40]$ the
overload switches, with pool $2$ becoming overloaded and
experiencing a sinusoidal arrival rate. However, we now take
fixed staffing in both service pools. In particular, the
parameters over the second overload interval $[20, 40]$ are

\vspace{0.05in} $\lm^n_1(t) = n$\: and\: $m^n_1(t) = n$; \quad
$\lm^n_2(t) = 1.1n + 0.1n \sin(t)$\: and\: $m^n_2(t) = n$.

\vspace{0.05in} \noindent Thus, we test two overload settings
in this example. In the first interval, we can see whether the
fluid approximation stabilizes. Since there is sharing of
class-$1$ customers, previous results such as in \cite{LW12a}
do not apply directly to our case. In the second interval, we
expect to see a sinusoidal behavior of the system, because the
staffing in both pools is fixed. In particular, the fluid model
should not approach a fixed point after the switch at time $t =
20$.

\begin{figure}[h!]
  \hfill
  \begin{minipage}[t]{.4\textwidth}
    \begin{center}
      \includegraphics[width=6cm]{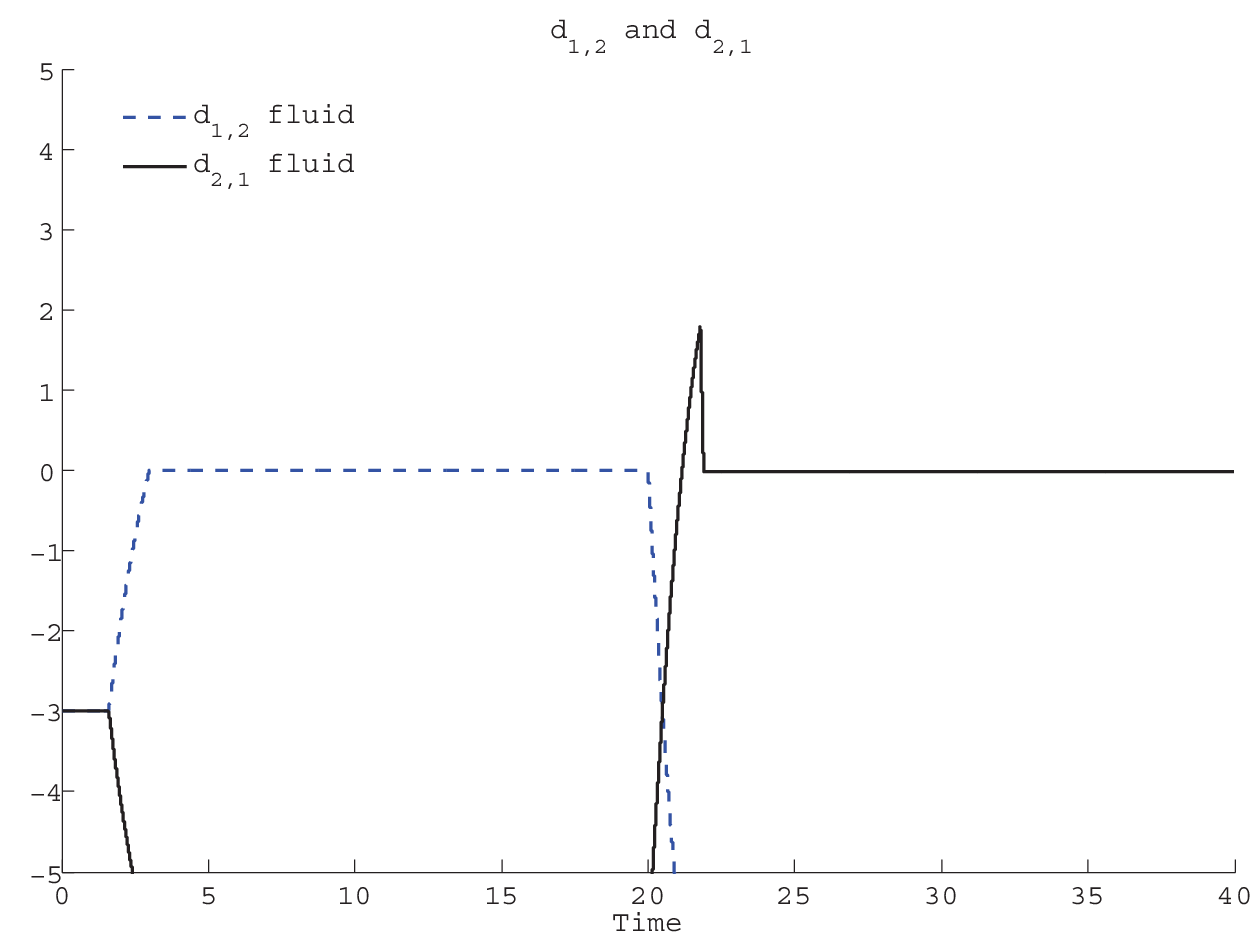}
      \caption{the two fluid difference functions $d_{1,2}$ and $d_{2,1}$ with the switching sinusoidal overloads}
      \label{figVaryDiff}
    \end{center}
  \end{minipage}
  \hfill
  \begin{minipage}[t]{.4\textwidth}
    \begin{center}
      \includegraphics[width=6cm]{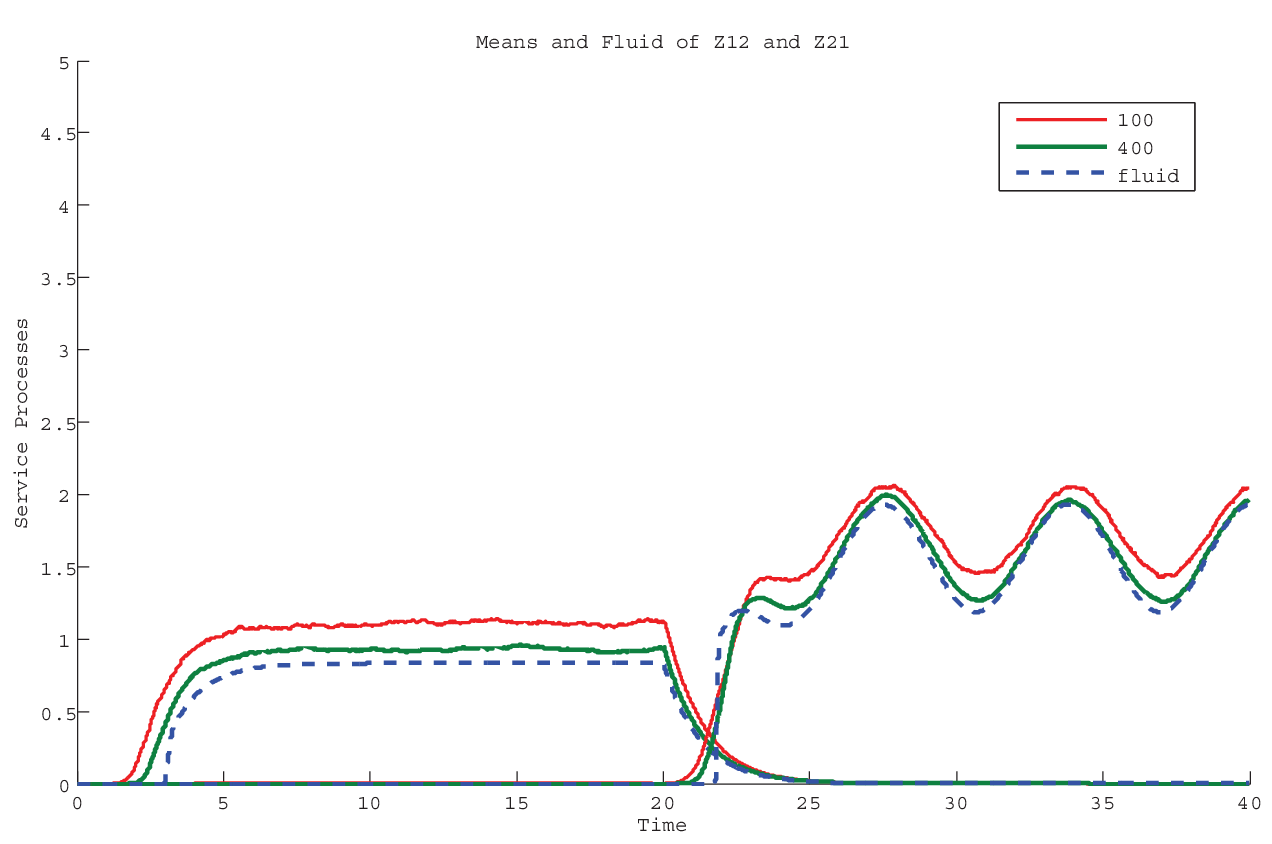}
      \caption{comparison of the fluid model to simulations of $10\barz^n_{1,2}$ and $10\barz^n_{2,1}$
      for $n = 100$ and $400$ with the switching sinusoidal overloads}
      \label{figVaryZ}
    \end{center}
  \end{minipage}
  \hfill
\end{figure}

\begin{figure}
  \begin{minipage}[t]{.4\textwidth}
    \begin{center}
      \includegraphics[width=6cm]{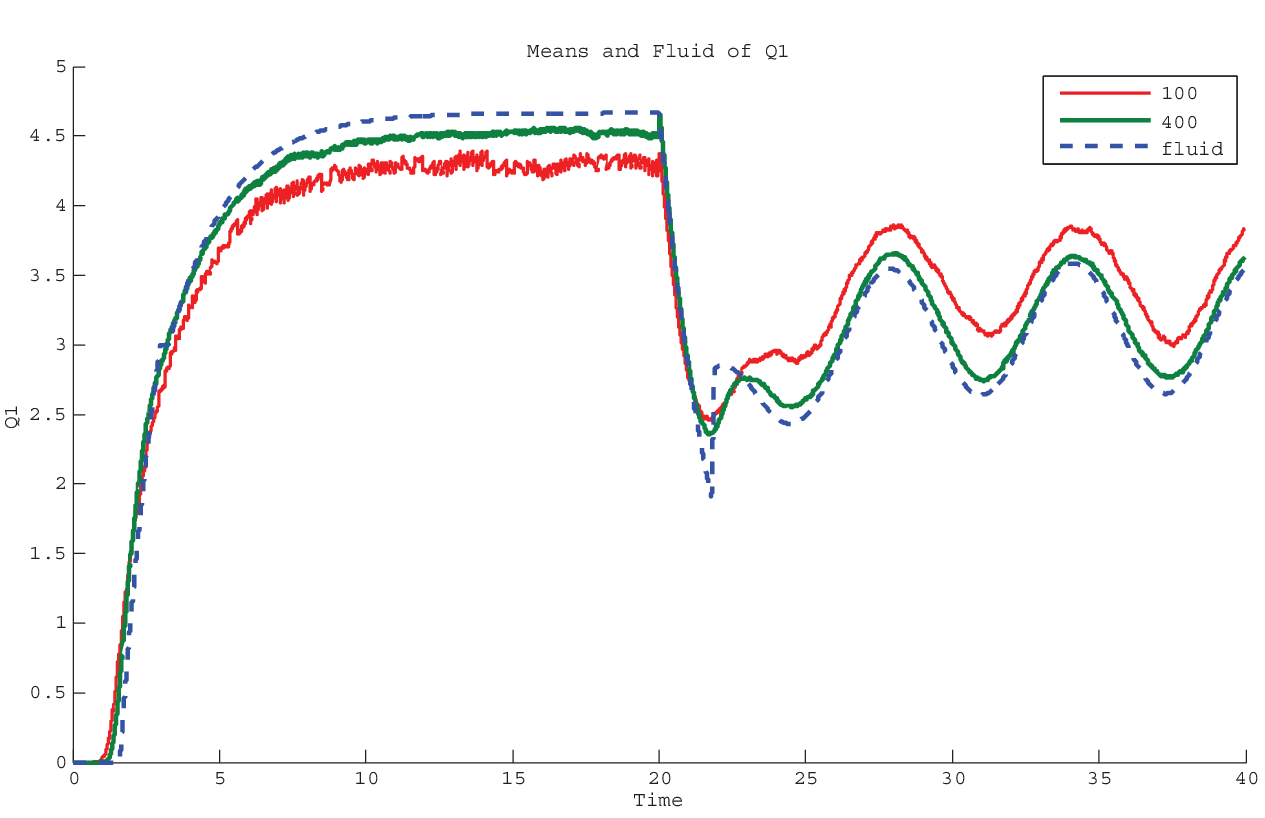}
      \caption{comparison of the fluid model to simulations of $10\barq^n_1$ for $n = 100$ and $400$ with the switching sinusoidal overloads}
      \label{figVaryQ1}
    \end{center}
  \end{minipage}
  \hfill
  \begin{minipage}[t]{.4\textwidth}
    \begin{center}
      \includegraphics[width=6cm]{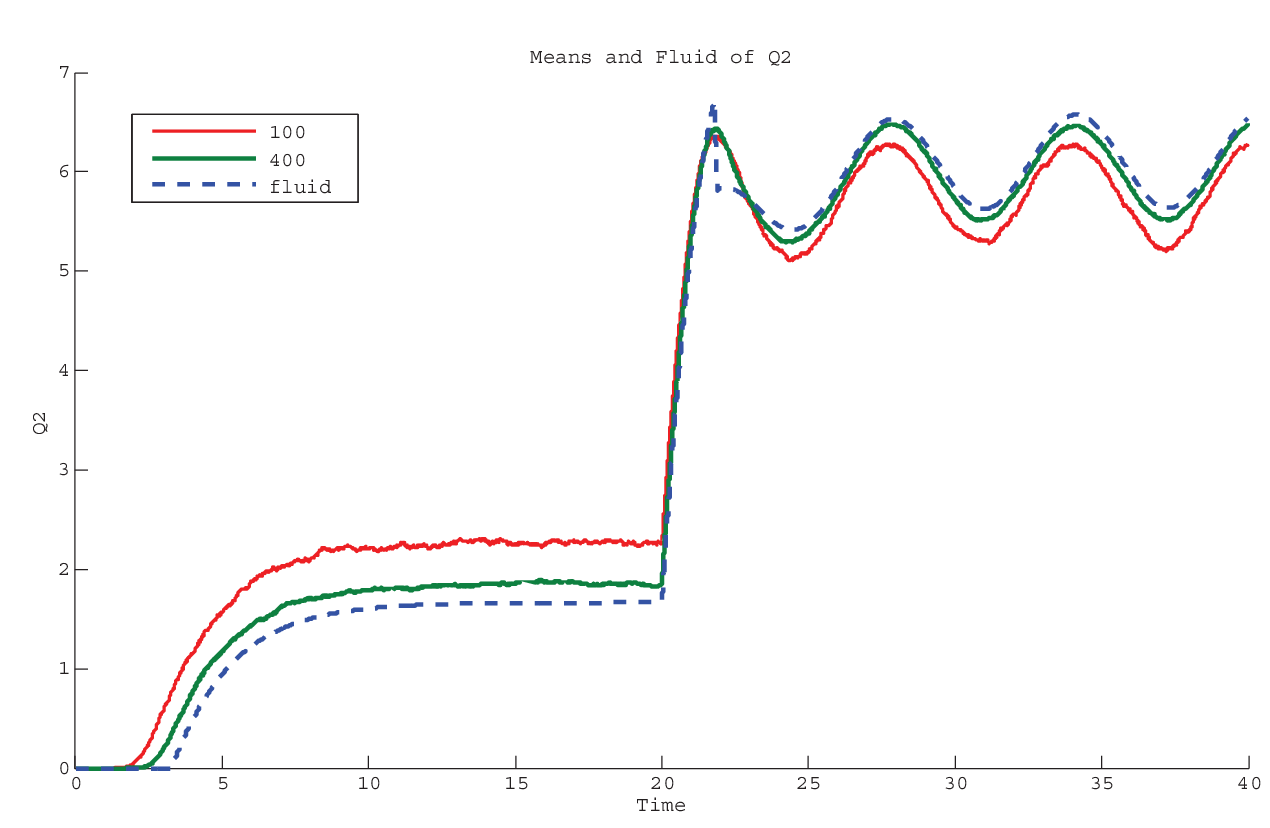}
      \caption{comparison of the fluid model to simulations of $10\barq^n_2$ for $n = 100$ and $400$ with the switching sinusoidal overloads}
      \label{figVaryQ2}
    \end{center}
  \end{minipage}
  \hfill
\end{figure}

We compare the fluid approximation to simulations for $n= 100$
and $n = 400$. Figures \ref{figVaryDiff}--\ref{figVaryQ2}
demonstrate the effectiveness of the fluid model and the
numerical algorithm.
As expected, the fluid over $[0, 20)$ approaches a fixed point,
and exhibits a sinusoidal behavior after $t = 20$, with the
accuracy of the fluid approximation increasing in the scale
parameter $n$.


\begin{figure}[h!]
\begin{center}
      \includegraphics[width=6cm]{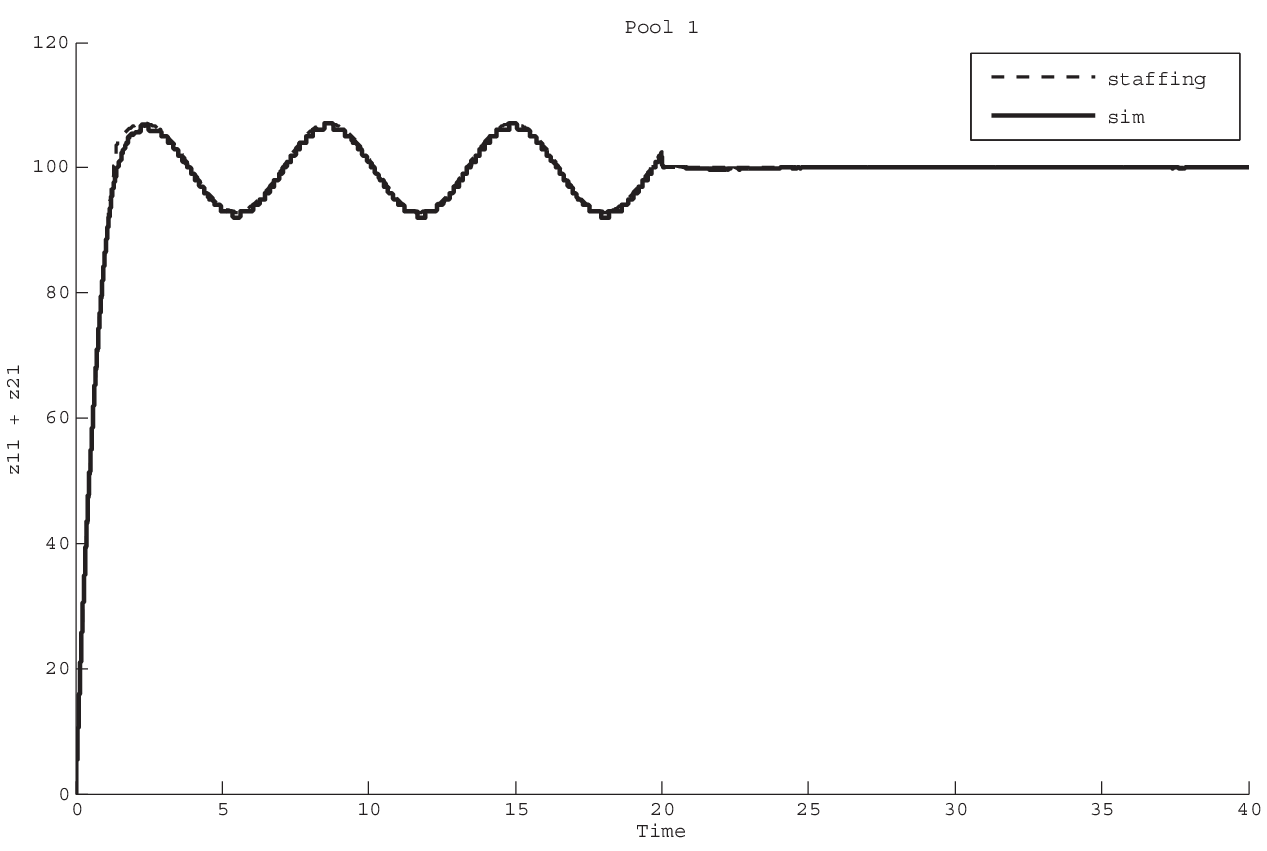}
      \caption{Fluid vs. simulations: Number of agents in pool $1$}
      \label{figVaryPool}
    \end{center}
\end{figure}

As was mentioned above, the fluid model requires special care
when the staffing functions are decreasing; we refer again to
\cite{LW12a}. Figure \ref{figVaryPool} shows the actual number
of agents in Pool $1$ for the case $n=100$ (the average of the
$1000$ simulations), and the staffing function $m^n_1(t)$ given
above. Clearly, the fluid model follows the actual staffing
closely. We further note that there is a downward jump in the
staffing function at time $t=20$. In the fluid model, we simply
eliminated the appropriate amount of staffing from the pool,
together with the fluid that was processed with that removed
capacity (this fluid in service is lost). However, in the
simulation, agents are removed only when they are done serving,
so there is no jump in the actual staffing at $t=20$, and no
customer in service is lost. Nevertheless, the fluid model with
the jump is clearly a good approximation for the stochastic
model with no jump. This behavior is to be expected, since
there are many service completions over short time intervals in
large systems.

\subsection{The Oscillatory Model} \label{secOscFluid}

\begin{figure}[h!]
  \hfill
   \begin{minipage}[t]{.4\textwidth}
    \begin{center}
      \includegraphics[width=6cm]{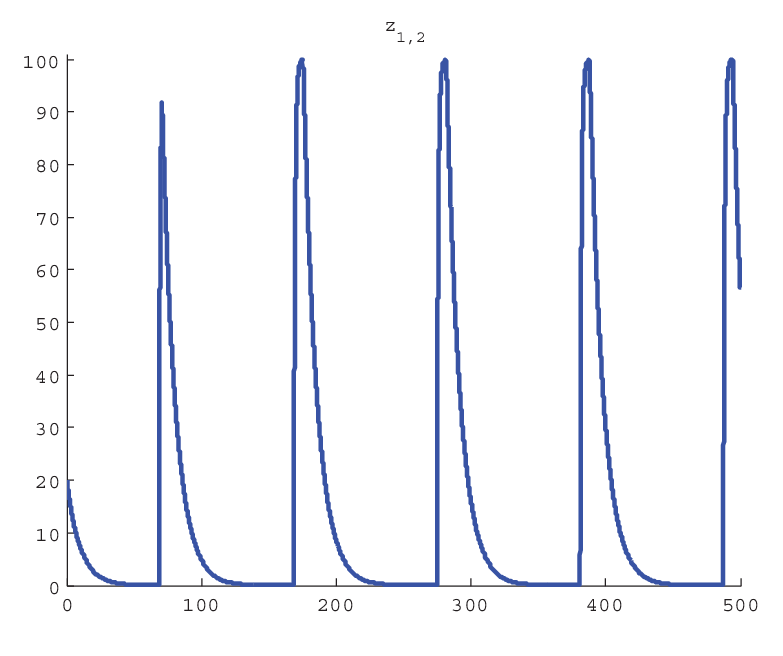}
      \caption{Oscillations $z_{1,2} (t)$ in the fluid model of the extreme example with no abandonment.}
      \label{figZ12oscF}
    \end{center}
  \end{minipage}
  \hfill
   \begin{minipage}[t]{.4\textwidth}
    \begin{center}
      \includegraphics[width=6cm]{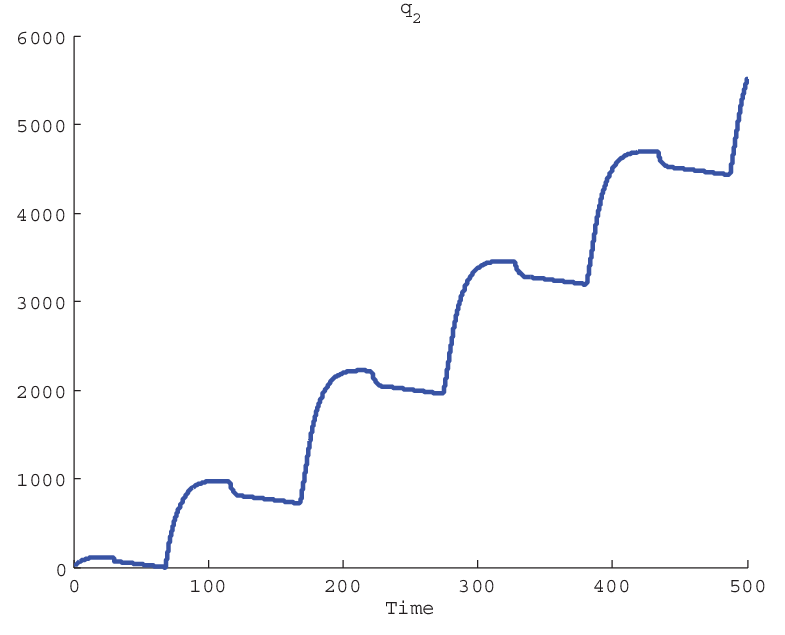}
      \caption{Oscillating growth of the content $q_{2} (t)$ in the fluid model of the extreme example with no abandonment.}
      \label{figQ2oscF}
    \end{center}
  \end{minipage}
  \hfill
\end{figure}

\begin{figure}[h!]
  \hfill
   \begin{minipage}[t]{.4\textwidth}
    \begin{center}
      \includegraphics[width=6cm]{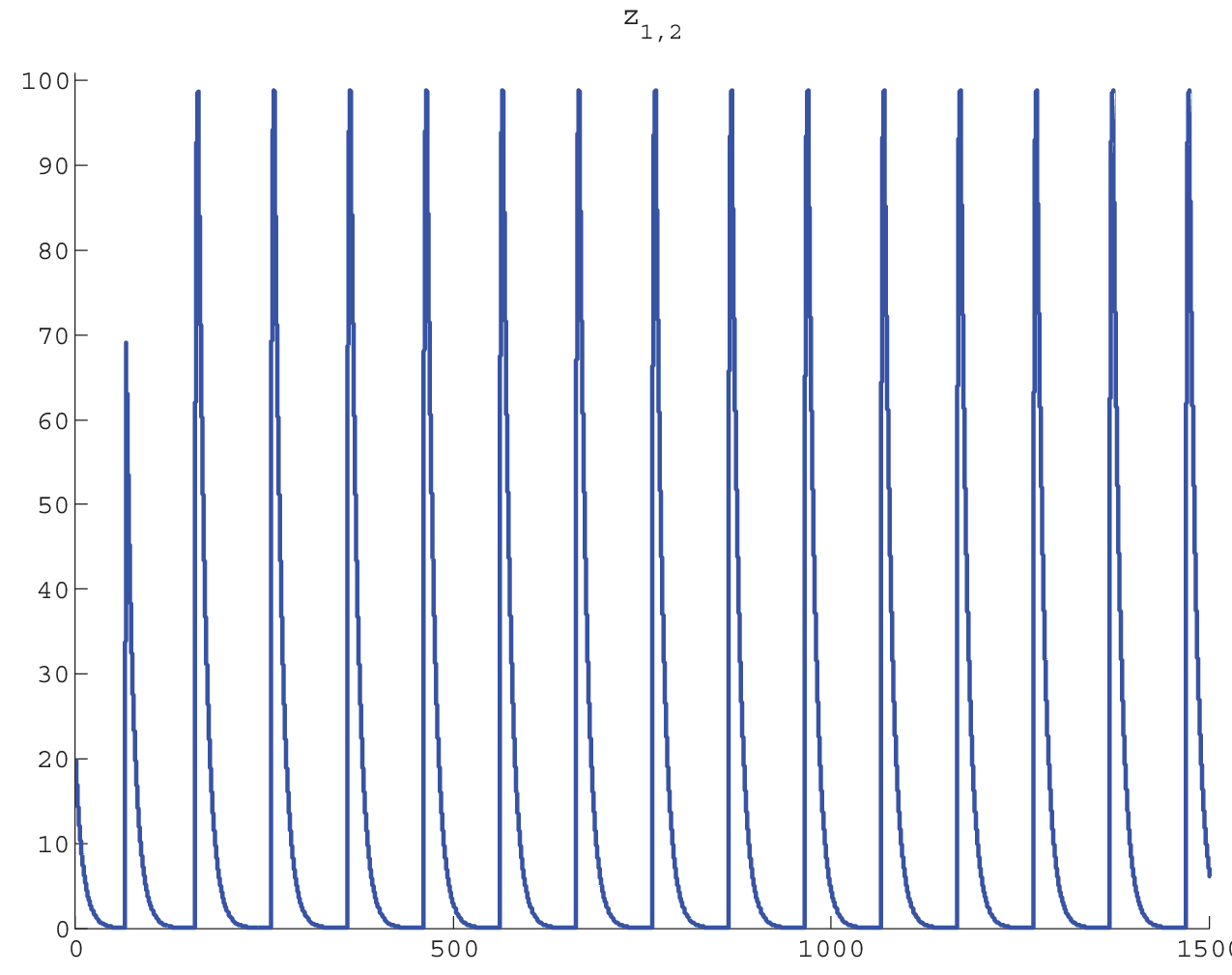}
      \caption{Oscillations $z_{1,2} (t)$ in the fluid model of the extreme example with abandonment.}
      \label{figZ12oscAbFluid}
    \end{center}
  \end{minipage}
  \hfill
   \begin{minipage}[t]{.4\textwidth}
    \begin{center}
      \includegraphics[width=6cm]{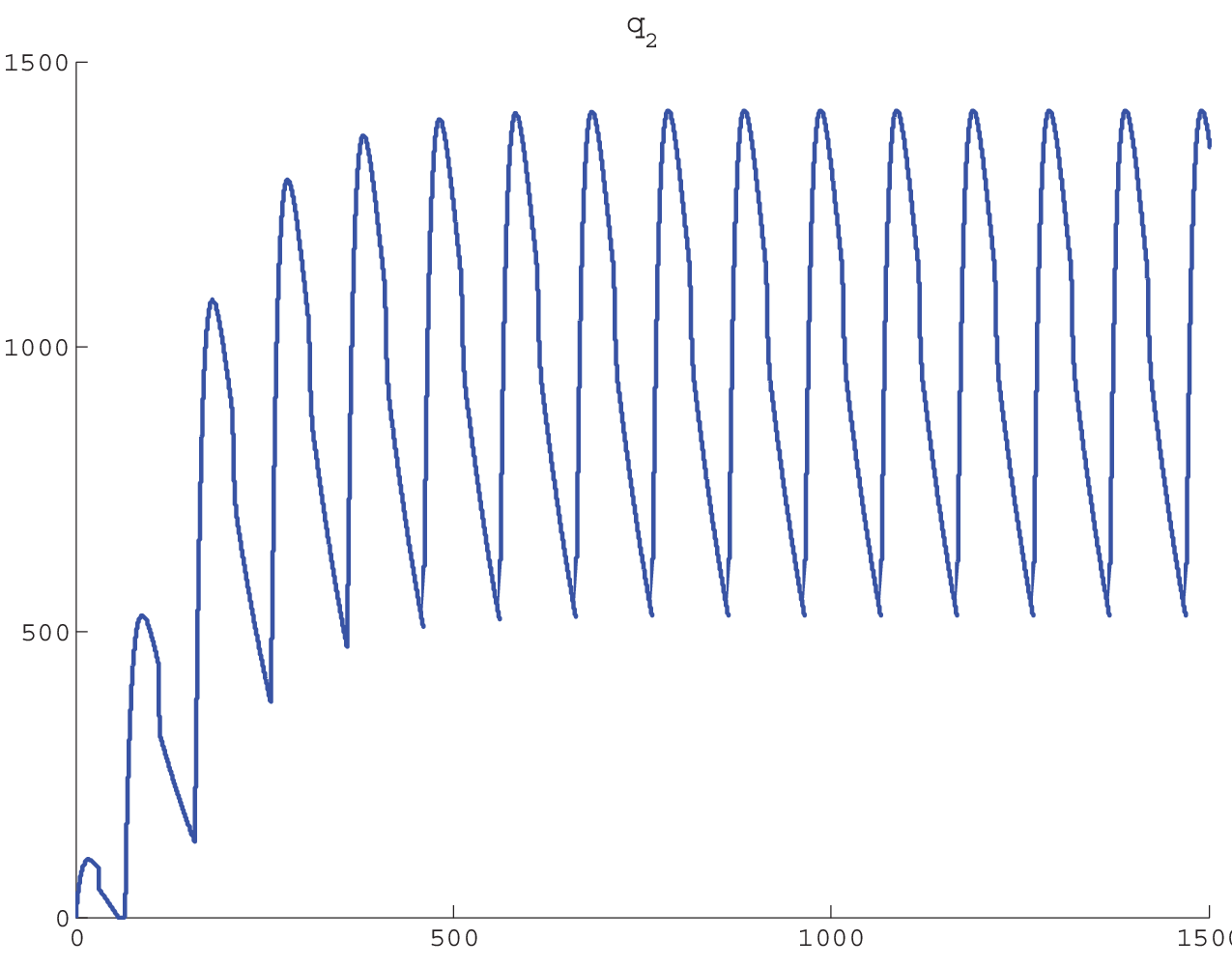}
      \caption{Oscillating growth of the content $q_{2} (t)$ in the fluid model of the extreme example with abandonment.}
      \label{figQ2oscAbFluid}
    \end{center}
  \end{minipage}
  \hfill
\end{figure}

Our final examples show that the fluid model can also predict
the bad oscillatory behavior. Here we consider the fluid model
of the examples shown in Figures
\ref{figZ12oscNo1}--\ref{figQ2oscAb} in \S \ref{secOSC}.
In particular, the parameters are $\mu_{1,1} = \mu_{2,2} = 1$,
$\mu_{1,2} = \mu_{2,1} = 0.1$, $\lm_1 = \lm_2 = 98$, $m_1 = m_2
= 100$ and $\tau_{i,j} = 0.01$ and $k_{i,j} = 10$, $i,j = 1,2$.
Figures \ref{figZ12oscF} and \ref{figQ2oscF} show the fluid
solution to the system with no abandonment (in \eqref{ODE-6} we
simply plug $\theta_1 = \theta_2 = 0$), whereas Figures
\ref{figZ12oscAbFluid} and \ref{figQ2oscAbFluid} show the fluid
solution to \eqref{ODE-6} with $\theta_1 = \theta_2 = 0.01$.

However, the initial conditions here are different than in
Figures \ref{figZ12oscNo1}--\ref{figQ2oscAb}. We now take
$z_{1,1}(0) = m_1 = 100$ and $z_{1,2}(0) = m_2 - z_{2,2}(0) =
20$. The reason is that, if the fluid is initialized with no
sharing and no queues, then its components $(q_1, q_2, z_{1,2},
z_{2,1})$ are fixed at $(0,0,0,0)$, i.e., there is never any
sharing, and the fluid queues are constant at zero. However, if
it is initialized at states with some sharing, then it may get
stuck at an oscillatory equilibrium, as shown in Figures
\ref{figZ12oscF} -- \ref{figQ2oscAbFluid}. In particular, this
is a numerical example that the fluid model may be {\em
bi-stable}, namely, have two very different stationary
behaviors. To which stationary behavior the fluid ends up
converging depends on the initial condition.

This fluid bi-stability property has two immediate implications
to the stochastic system. First, once an overload incident is
ending, with substantial sharing taking place, the system may
start to oscillate. Indeed, this is the case in the example
shown in Figures \ref{figZosc}.
The second implication is that the no-sharing equilibrium may
be unstable in practice, because stochastic noise can
eventually ``push'' the system out of this equilibrium, and
cause it to oscillate. That was demonstrated in Figures
\ref{figZ12oscNo1} and \ref{figQ2oscNo1} in \S \ref{secOSC}.
(Recall that the initial condition of the example in \S
\ref{secOSC} was of an empty system. In particular, with no
sharing initially.) Note also that the time scale in Figures
\ref{figZ12oscF} and \ref{figQ2oscF} is shorter than in Figures
\ref{figZ12oscAbFluid} and \ref{figQ2oscAbFluid}. As for the
corresponding figures in \S \ref{secOSC}, the time scale of the
second example is longer to make it clear that the system with
abandonment converges to an oscillatory equilibrium.

\section{Conclusions} \label{secConclude}

In this paper we studied a time-varying X model experiencing
periods of overloads. While our previous FQR-T control is
effective in automatically responding quickly to unexpected
overloads, the examples in \S \ref{secRelaxOneWay} and \S
\ref{secOSC} show that it needs to be modified to recover
rapidly after the overload is over, due to either a return to
normal loading or a sudden change in the direction of the
overload.  We thus proposed the {\em fixed-queue-ratio with
activation-and-release-thresholds} (FQR-ART) control. With
FQR-ART, the one-way sharing rule is relaxed by adding the
lower release thresholds. To avoid oscillations of the service
process, which in turn can cause congestion collapse, we
indicated that the activation thresholds also need to be
increased, being asymptotically of order $O(n)$ as in
\eqref{ThreshQ} instead of $o(n)$, as in \eqn{thresholds} with
FQR-T.

We then extended the fluid model developed in
\cite{PW11a,PW11b,PW13a} based on the stochastic averaging
principle to cover a more general time-varying environment. and
developed the corresponding algorithm to numerically compute
the performance functions in that fluid model. Simulation
experiments indicate that this fluid model captures the main
dynamics of the system, even in extreme cases, as the one
considered in \eqref{secSimNone}.  Thus the fluid model can be
used to ensure that the control parameters of FQR-ART are set
properly.

There are many directions for future research. First, it
remains to investigate the performance of FQR-ART in more
complex time-varying scenarios. Second, it remains to establish
theoretical properties of the new fluid model, paralleling
\cite{PW11b}. Third, it remains to establish many-server
heavy-traffic limits in this more general setting, paralleling
\cite{PW13a,PW13b}.

Fourth, and most important for engineering applications, it remains to extend the sharing
mechanism to more than two systems.  With more than two systems, there are more possible overload scenarios.
One scenario involves only a single system experiencing an overload.
For that scenario, assistance might be provided by several other systems, at less cost to each.
The previous results for multi-class multi-pool systems
in \cite{GW09a,GW09b,GW10} indicate that such an extension should be possible.
For example, one system might be judged to be
overloaded, and assistance from others might be activated, perhaps with help only provided to the one system experiencing the overload (the analog of one-way sharing),
if its queue length exceeds a specified proportion of the total queue length plus some activation threshold.
Then sharing, with help only provided to the one overloaded system, might aim to keep the queue length of the overloaded system close to its target proportion.
But then, as proposed here, evidently release thresholds should be used to ensure rapid recovery after the overload is over.

\noindent
\textbf{Acknowledgement} \\
The authors received support from NSF grants CMMI 1066372 and
1265070.

\bibliographystyle{nrl}

\begin{thebibliography}{99}

\bibitem{AAM07} Aksin, Z., M. Armony, V. Mehrotra. 2007. The
    modern call center: a multi-disciplinary perspective on
    operations management research. {\em Production Oper.
    Management}, {\bf 16} (6) 655--688.

\bibitem{ArmonyPatient} Armony, M., S. Israelit, A. Mandelbaum,
    Y. Marmor, Y. Tseytlin, G. Yom-tov. 2010. Patient flow in
    hospitals: a data-based queueing-science perspective. {\em
    Working paper}, New York University.


\bibitem{BFI05} Baier, V., R. F\"{o}disch, A. Ihring, E.
    Kessler, J. Lerchner, G. Wolf, J. M. K\"{o}hler, M.
    Nietzch, M. Krugel.  2006. Highly sensitive thermopile heat
    power sensor for micro-fluid calorimetry of biochemical
    processes. {\em Sensors and Actuators A} {\bf 123-124} (23)
    354--359.

\bibitem{BW98} Berger, A, W., W. Whitt. 1998. Effective
    bandwidths with priorities. {\em IEEE/ACM Transactions on
    Networking}, {\bf 6} (4), 447--460.



\bibitem{BBHA12} Boyle, A., K. Beniuk, I. Higginson, P.
    Atkinson. 2012. Emergency department crowding: Time for
    interventions and policy evaluations. {\em Emergency
    Medicine J.} {\bf 29} 460--466.

\bibitem{CLW95} Choudhury G. L., K. K. Leung, W. Whitt. 1995.
    Efficiently providing multiple grades of service with
    protection against overloads in shared resources. {\em
    AT\&T Technical Journal}, {\bf 74} (4), 50--63.

\bibitem{CAB11} Chan, C. W., M. Armony, N. Bambos. 2011.
    Fairness in overloaded parallel queues. {\em Working
    paper}, Columbia University, arXiv:1011.1237v2

\bibitem{ChanGalit} Chan, C. W., G. Yom-Tov, G. Escobar. 2011.
    When to use speedup: an examination of intensive care units
    with readmissions. {\em Working paper}, Columbia
    University.


\bibitem{DG11} Deo, S., I. Gurvich. 2011. Centralized versus
    decentralized ambulance diversion: a network perspective.
    {\em Mangement Sci.} {\bf 57} (7),  1300--1319.

\bibitem{DH86} Doshi, B., H. Heffes. 1986. Overload performance
    of several processor queueing disciplines for the M/M/1
    queue. {\em IEEE Transactions on Communications}, {\bf 34}
    (6), 538--546.

\bibitem{EF91} Erramilli, A., L. J. Forys. 1991. Oscillations
    and chaos in a flow model of a switching system. {\em IEEE
    journal on Selected Areas in Communications}, {\bf 9} (2),
    171--178.

\bibitem{FR94} Feinberg, E. A., M. I. Reiman. 1994. Optimality
    of randomized trunk reservation. {\em Prob Eng. Inf. Sci.}
    {\bf 8} (4) 463--489.

\bibitem{FMMW08} Feldman, Z., Mandelbaum, A., Massey, W. A.,
    Whitt, W. 2008. Staffing of time-varying queues to achieve
    time-stable performance. {\em Management Sci.} {\bf 54}
    (2), 324--338.

\bibitem{FF99} Floyd, S., K. Fall. 1999. Promoting the use of
    end-to-end congestion control in the Internet {\em IEEE/ACM
    Transactions on Networking}, {\bf 7} (4), 458--472.



\bibitem{GMR02}
Garnett, O., A. Mandelbaum, M. I. Reiman. 2002. Designing a
call center with impatient customers. {\em Manufacturing
Service Oper. Management}, {\bf 4} (3), 208--227.

\bibitem{GK04} Goldstein, M. A., K. A. Kavajecz. 2004. Trading
    strategies during circuit breakers and extreme market
    movements. {\em Journal of Financial Markets}, {\bf 7}
    301--333.


\bibitem{GurPer12} Gurvich, I., O. Perry. 2012. Overflow
    networks: Approximations and implications to call-center
    outsourcing. {\em Oper. Res.}, {\bf 60} (4), 996--1009.

\bibitem{GW09a} Gurvich, I., W. Whitt. 2009a.
    Queue-and-idleness-ratio controls in many-server service
    systems. {\em Math. Oper. Res.} {\bf 34} (2) 363--396.


\bibitem{GW09b} Gurvich, I., W. Whitt. 2009b. Scheduling
    flexible servers with convex delay costs in many-server
    service systems. {\em Manufacturing Service Oper.
    Management} {\bf 11} (2) 237--253.

\bibitem{GW10} Gurvich, I., W. Whitt. 2010. Service-level
    differentiation in many-server service systems via
    queue-ratio routing. {\em Oper. Res.} {\bf 58} (2)
    316--328.

\bibitem{HJM09} Hampshire, R. C., O. B. Jennings, W. A. Massey.
    2009. A time-varying call center design via Lagrangian
    mechanics. {\em Probability in the Engineering and
    Informational Sciences}, {\bf 23} (2), 231--259.

\bibitem{HMW09} Hampshire, R. C., W. A. Massey, Q. Wang. 2009.
    Dynamic pricing to control loss systems with quality of
    service targets. {\em Probability in the Engineering and
    Informational Sciences}, {\bf 23} (2), 357--383.

\bibitem{HZ05} Harrison J. M., A. Zeevi. 2005. A method for
    staffing large call centers using stochastic fluid models.
    {\em Manufacturing Service Oper. Management}, {\bf 7} (1),
    20--36.





\bibitem{K91} Kelly, F. P. 1991. Loss networks. {\em Ann. Appl.
    Probab.} {\bf 1} (3), 319--378.


\bibitem{K02} Khalil, H. K. 2002. {\em Nonlinear Systems.}
    Prentice Hall, New Jersey.


\bibitem{K72} Koopman, B. O. 1972. Air-terminal queues under
    time-dependent conditions. {\em Operations Research}, {\bf
    20} (6) 1089--1114.




\bibitem{Ko91} K\"{o}rner U. 1991. Overload control of SPC
    systems. {\em International Teletraffic Congress, ITC 13},
    Copenhagen, Denmark.

\bibitem{LW11} Liu, Y., W. Whitt. 2011. Nearly periodic
    behavior in the the overloaded G/D/S+GI Queue. {\em
    Stochastic Systems} {\bf 1} (2) 340--410.

\bibitem{LW12a} Liu, Y., W. Whitt. 2012a. The $G_t/GI/s_t+GI$
    many-server fluid queue. {\em Queueing Systems}, {\bf 71}
    (4) 405--444.

\bibitem{LW12b} Liu, Y., W. Whitt. 2012b. Stabilizing customer
    abandonment in many-server queues with time-varying
    arrivals. {\em Operations Research}, {\bf 60} (6)
    1551--1564.

\bibitem{LW12c} Liu, Y., W. Whitt. 2012c. A many-server fluid
    limit for the $G_t/GI/s_t+GI$ queueing model experiencing
    periods of overloading. {\em Operations Research Letters},
    {\bf 40} 307--312.

\bibitem{LPD02} Low, S. H., F. Paganini, J. C. Doyle. 2002.
    Internet congestion control. {\em Control Systems}, {\bf
    22} (1), 28--43.

\bibitem{N82} Newell, G. F. 1982. {\em Applications of Queueing
    Theory}, Chapman Hall, London.



\bibitem{PTW07} Pang, G., Talreja, R., Whitt, W., 2007.
    Martingale proofs of many-server heavy-traffic limits for
    Markovian queues. \textit{Probability Surveys}, \textbf{4},
    193--267.

\bibitem{PW09} Perry, O., W. Whitt. 2009. Responding to
    unexpected overloads in large-scale service systems. {\em
    Management Sci.}, {\bf 55} (8), 1353--1367.

\bibitem{PW11a} Perry, O., W. Whitt. 2011a. A fluid
    approximation for service systems responding to unexpected
    overloads. {\em Oper. Res.}, {\bf 59} (5), 1159--1170.


\bibitem{PW11b} Perry, O., W. Whitt. 2011b. An ODE for an
    overloaded X model involving a stochastic averaging
    principle. {\em Stochastic Systems}, {\bf 1} (1), 17--66.

\bibitem{PW13a} Perry, O., W. Whitt. 2013a. A fluid limit for an
    overloaded X model via a stochastic averaging principle. {\em Math,
    Oper. Res.} {\bf 13} (2), 294-349.


\bibitem{PW13b} Perry, O., W. Whitt. 2013b. Diffusion
    approximation for an overloaded X model via a stochastic
    averaging principle. {\em Queueing Systems}, published online May 24, 2013.


\bibitem{PW14} Perry, O., W. Whitt. 2014. The dynamics of a symmetric $X$ fluid model with inefficient sharing.
Paper in preparation.

\bibitem{PowEtal} Powell, E.S., R. K. Khare, A. K. Venkatesh,
    B. D. Van Roo, J. G. Adams, G. Reinhardt. 2012. The
    relationship between inpatient discharge timing and
    emergency department boarding. {\em Journal of Emergency
    Medicine}, {\bf 42} (2), 186--196.

\bibitem{SKT90} Schulzrinne, H., J. F. Kurose, D. Towsley.
    1990. Congestion control for real-time traffic in
    high-speed networks. {\em IEEE proceeding in Ninth Annual
    Joint Conference of the IEEE Computer and Communication
    Societies}, 543--550.



\bibitem{SW11} Shah, D., D. Wischik. 2011. Fluid models of
    congestion collapse in overloaded switched networks. {\em
    Queueing Systems} {\bf 69} 121--143.

\bibitem{Dai-ED} Shi, P., M. Chou, J. G. Dai, D. Ding, J. Sim.
    2012. Hospital Inpatient Operations: Mathematical Models
    and Managerial Insights. {\em Working paper}

\bibitem{S98} Sontag, E. D. 1998. {\em Mathematical Control
    Theory}, second edition, Springer, New York.


\bibitem{ST11} Stolyar, A. L., T. Tezcan. 2011. Shadow-routing
    based control of flexible multiserver pools in overload.
    {\em Oper. Res.} {\bf 59} (6), 1427--1444.

\bibitem{T09} Teschl, G. 2009. \emph{Ordinary Differential
    Equations and Dynamical Systems},
Universit\"{a}t Wien. 
Available online:
www.mat.univie.ac.at/$\sim$gerald/ftp/book-ode/ode.pdf

\bibitem{W64} Weber, J. H. 1964. A simulation study of routing
    control in communication networks. {\em Bell System Tech.
    J.} {\bf 43} 2639--2676.

\bibitem{WJLH06} Wei, D. X., C. Jin, S. H. Low, S. Hegde. 2006.
    FAST TCP: motivation, architecture, algorithms,
    performance. {\em IEEE/ACM Transactions on Networking},
    {\bf 14} (6), 1246--1259.

\bibitem{W02} 
Whitt, W. 2002. {\em Stochastic-Process Limits}, New York,
Springer, 2002.

\bibitem{W04} 
Whitt, W. 2004. Efficiency-driven heavy-traffic approximations
for many-server queues with abandonments.
{\em Management Sci.} {\bf 50} (10) 1449--1461.


\bibitem{W56} Wilkinson, R. I. 1956. Theory for toll traffic
    engineering in the U.S.A. {\em Bell System Tech. J.} {\bf
    35} 421--513.


\bibitem{YGGG10} Yankovic, N., S. Glied, L. V. Green, M. Grams.
    2010. The impact of ambulance diversion on heart attack
    deaths. {\em Inquiry} {\bf 47} (1) 81--91.

\end{thebibliography}

\end{document}